\documentclass[reqno]{amsart}
\usepackage{amsmath,amssymb}
\newtheorem{theorem}{Theorem}[section]
\newtheorem{proposition2}[theorem]{Proposition}
\newtheorem{definition2}[theorem]{Definition}
\newtheorem{lemma2}[theorem]{Lemma}
\newtheorem{corollary2}[theorem]{Corollary}
\newtheorem{remark2}[theorem]{Remark}
\newtheorem{example2}[theorem]{Example}
\numberwithin{equation}{section}
\numberwithin{theorem}{section}

\newcommand{\mc}[1]{{\mathcal #1}}
\newcommand{\bb}[1]{{\mathbb #1}}
\newcommand{\Glimsup}{\mathop{\textrm{$\Gamma\!$--$\varlimsup$}}\limits}
\newcommand{\Gliminf}{\mathop{\textrm{$\Gamma\!$--$\varliminf$}}\limits}
\newcommand{\Glim}{\mathop{\textrm{$\Gamma\!$--$\lim$}}\limits}
\newcommand{\res}{\mathop{\hbox{\vrule height 7pt width .5pt depth
               0pt\vrule height .5pt width 6pt depth 0pt}}\nolimits}
\newcommand{\upbar}[1]{\,\overline{\! #1}}
\newcommand{\id}{{1 \mskip -5mu {\rm I}}}

\begin{document}

\title{$\Gamma$-entropy cost for scalar conservation laws}

\author[G.\ Bellettini]{Giovanni Bellettini} 
\address{Giovanni
  Bellettini, Dipartimento di Matematica, Universit\`a di Roma `Tor
  Vergata', Via della Ricerca Scientifica 00133 Roma, Italy}
\email{Giovanni.Bellettini@lnf.infn.it}

\author[L.\ Bertini]{Lorenzo Bertini} 
\address{Lorenzo Bertini,
  Dipartimento di Matematica, Universit\`a di Roma `La Sa\-pien\-za',
  P.le Aldo Moro 2, 00185 Roma, Italy} 
\email{bertini@mat.uniroma1.it}

\author[M.\ Mariani]{Mauro Mariani} 
\address{Mauro Mariani,
  Dipartimento di Matematica, Universit\`a di Roma `La Sa\-pien\-za',
  P.le Aldo Moro 2, 00185 Roma, Italy \\
and CEREMADE, UMR-CNRS 7534, Universit\'e de Paris-Dauphine \\
  Place du Marechal de Lattre de TASSIGNY, F-75775 Paris Cedex 16}
\email{mariani@ceremade.dauphine.fr}

\author [M.\ Novaga]{Matteo Novaga} 
\address{Matteo Novaga,
  Dipartimento di Matematica, Universit\`a di Pisa, Largo Bruno
  Pontecorvo 5, 56127 Pisa, Italy} 
\email{novaga@dm.unipi.it}

\maketitle

\begin{abstract}
  We are concerned with a control problem related to the vanishing
  viscosity approximation to scalar conservation laws. We investigate
  the $\Gamma$-convergence of the control cost functional, as the
  viscosity coefficient tends to zero. A first order $\Gamma$-limit is
  established, which characterizes the measure-valued solutions to the
  conservation laws as the zeros of the $\Gamma$-limit. A second order
  $\Gamma$-limit is then investigated, providing a characterization of
  entropic solutions to conservation laws as the zeros of the
  $\Gamma$-limit.
\end{abstract}

\section{Introduction}
\label{s:1}
We are concerned with the scalar one-dimensional conservation law
\begin{equation}
  \label{e:1.1}
   u_t+f(u)_x=0 
\end{equation}
where, given $T>0$, $u= u(t,x)$, $(t,x)\in [0,T]\times \bb R$,
subscripts denote partial derivatives, and the \emph{flux} $f$ is a
Lipschitz function. As well known, even if the initial datum $u(0)=
u(0,\cdot)$ is smooth, the flow \eqref{e:1.1} may develop
singularities for some positive time. In general, these singularities
appear as discontinuities of $u$ and are called \emph{shocks}. It is
therefore natural to interpret \eqref{e:1.1} weakly; in the weak
formulation uniqueness is however lost, if no further conditions are
imposed. Given a function $\eta$, called \emph{entropy}, the
conjugated \emph{entropy flux} $q$ is defined up to an additive
constant as $q(u)=\int^u\!dv\,\eta'(v) \, f'(v)$. A weak solution to
\eqref{e:1.1} is called \emph{entropic} iff for each entropy --
entropy flux pair $(\eta,q)$ with $\eta$ convex, the inequality
$\eta(u)_t + q(u)_x \le 0$ holds in the sense of distributions. Note
that the entropy condition is always satisfied for smooth solutions to
\eqref{e:1.1}.  The classical theory, see e.g.\ \cite{Da,Se}, shows
existence and uniqueness in $C\big([0,T];L_{1,\mathrm{loc}}(\bb
R)\big)$ of the entropic solution to the Cauchy problem associated to
\eqref{e:1.1}. While the flow \eqref{e:1.1} is invariant w.r.t.\
$(t,x)\mapsto (-t,-x)$, the entropy condition breaks such invariance
and selects the ``physical'' direction of time.

In the conservation law \eqref{e:1.1} the viscosity effects are
neglected. This approximation is no longer valid if the gradients
become large as it happens when shocks appear. A more accurate
description is then given by the parabolic equation
\begin{equation}
  \label{e:1.2}
  u_t+f(u)_x= \frac{\varepsilon}2 \big( D(u)u_x \big)_x  
\end{equation}
in which $(t,x) \in [0,T]\times \bb R$, $D$, assumed uniformly
positive, is the diffusion coefficient and $\varepsilon>0$ is the
viscosity. In this context of scalar conservation laws, it is also
well known that, as $\varepsilon \to 0$, equibounded solutions to
\eqref{e:1.2} converge in $L_{1,\mathrm{loc}}(\bb [0,T]\times \bb R)$
to entropic solutions to \eqref{e:1.1}, see e.g.\ \cite{Da,Se}. This
approximation result shows that the entropy condition is relevant.

Perhaps less well known, at least in the hyperbolic literature, is the
fact that entropic solutions to \eqref{e:1.1} can be obtained as
scaling limit of discrete stochastic models of lattice gases, see
e.g.\ \cite[Ch.~8]{KL}. In a little more detail, consider particles
living on a one-dimensional lattice and randomly jumping to their
neighboring sites. It is then proven that, under hyperbolic scaling,
the empirical density of particles converges in probability to
entropic solutions to \eqref{e:1.1}.  A much studied example is the
\emph{totally asymmetric simple exclusion process}, where there is at
most one particle in each site and only jumps heading to the right are
allowed.  In this case, the empirical density takes values in $[0,1]$
and its scaling limit is given by \eqref{e:1.1} with flux
$f(u)=u(1-u)$. In this stochastic framework, it is also worth looking
at the large deviations asymptotic associated to the aforementioned
law of large numbers.  Basically, this amounts to estimate the
probability that the empirical density lies in a neighborhood of a
given trajectory. In general this probability is exponentially small,
and the corresponding decay rate is called the large deviations rate
functional. For the totally asymmetric simple exclusion process, this
issue has been analyzed in \cite{Je,Va}. It is there shown that the
large deviations rate functional is infinite off the set of weak
solutions to \eqref{e:1.1}; on such solutions the rate functional is
given by the total positive mass of the entropy production $h(u)_t +
g(u)_x$ where $h$ is the Bernoulli entropy, i.e.\ $h(u)=-u\log u -
(1-u)\log(1-u)$ and $g$ is its conjugated entropy flux.

A stochastic framework can also be naturally introduced in a PDE
setting by adding to \eqref{e:1.2} a random perturbation, namely
\begin{equation}
  \label{e:1.3}
  u_t+f(u)_x= \frac{\varepsilon}2 \big( D(u)u_x \big)_x
   + \sqrt{\gamma} \big(\sqrt{\sigma(u)}\alpha_\gamma \big)_x 
    \quad (t,x) \in (0,T)\times \bb R 
\end{equation}
where $\sigma(u)\ge 0$ is a conductivity coefficient and
$\alpha_\gamma$ is a Gaussian random forcing term white in time and
with spatial correlations on a scale much smaller than $\gamma$. Let
$u^{\varepsilon,\gamma}$ be the corresponding solution; if $\gamma \ll
\varepsilon$ then $u^{\varepsilon,\gamma}$ still converges in
probability to the entropic solution to \eqref{e:1.1} and the large
deviations asymptotic becomes a relevant issue. Referring to \cite{Ma}
for this analysis, here we formulate the problem from a purely
variational point of view quantifying, in terms of the parabolic
problem \eqref{e:1.2}, the asymptotic \emph{cost} of non-entropic
solutions to \eqref{e:1.1}. Introducing in \eqref{e:1.2} a
\emph{control} $E \equiv E(t,x)$ we get
\begin{equation}
  \label{e:1.4}
   u_t+f(u)_x= \frac{\varepsilon}2 \big( D(u)u_x \big)_x
   - \big( \sigma(u)E \big)_x \qquad (t,x) \in (0,T)\times \bb R 
\end{equation}
If we think of $u$ as a density of charge, then $E$ can be naturally
interpreted as the `controlling' external electric field and
$\sigma(u)\ge 0 $ as the conductivity. The flow \eqref{e:1.4}
conserves the \emph{total charge} $\int\!dx\, u(t,x)$, whenever it is
well defined.

The \emph{cost functional} $I_\varepsilon$ associated with
\eqref{e:1.2} can be now informally defined as the work done by the
optimal controlling field $E$ in \eqref{e:1.4}, namely
\begin{equation}
  \label{e:1.5}
I_\varepsilon(u) = \inf_E \frac 12 \int_{[0,T]}\!dt\,dx \, \sigma(u) E^2
                = \inf_E \frac 12 \int_{[0,T]}\!dt\, 
                         \big\| E \big\|^2_{ L_2(\bb R, \sigma(u) dx)}
\end{equation}
where the infimum is taken over the controls $E$ such that
\eqref{e:1.4} holds. For a suitable choice of the random perturbation
$\alpha_\gamma$, $I_\varepsilon$ is the large deviations rate
functional of the process $u^{\varepsilon,\gamma}$ solution to
\eqref{e:1.3}, when $\varepsilon$ is fixed and $\gamma\to 0$. To avoid
the technical problems connected to the possible unboundedness of the
density $u$, we assume that the conductivity $\sigma$ has compact
support. In this case, if $u$ is such that $I_\varepsilon(u)<+\infty$
then $u$ takes values in the support of $\sigma$, see
Proposition~\ref{t:boun} for the precise statement. For the sake of
simplicity, we assume that $\sigma$ is supported by $[0,1]$. The case
of strictly positive $\sigma$ also fits in the description below,
provided however that the analysis is a priori restricted to
equibounded densities $u$.

In this paper we analyze the variational convergence of
$I_\varepsilon$ as $\varepsilon\to 0$. Our first result holds for a
Lipschitz flux $f$, and identifies the so-called $\Gamma$-limit of
$I_\varepsilon$, which is naturally studied in a Young measures
setting.  The limiting cost of a Young measure $\mu \equiv
\mu_{t,x}(d\lambda)$ is
\begin{equation*}
\mc I (\mu) =  \frac 12 \int_{[0,T]}\!dt \,
               \Big\|  \big[ \mu(\lambda) \big]_t
                 +\big[ \mu( f(\lambda) )  \big]_x
               \Big\|_{H^{-1}(\bb R, \mu (\sigma(\lambda)) dx)}^2
\end{equation*}
where, given $F \in C([0,1])$ we set
$[\mu(F(\lambda))](t,x) = \int\!
\mu_{t,x}(d\lambda) \, F(\lambda)$ and, with a little abuse of
notation, $\| \varphi \|_{H^{-1}(\bb R, \mu_{t,\cdot}
  (\sigma(\lambda)) dx)}$ is the dual norm to $\big[\int\!dx\,
\mu_{t,x}(\sigma(\lambda)) \, \varphi_x^2\big]^{1/2}$.

Note that $\mc I(\mu)$ vanishes iff $\mu$ is a \emph{measure-valued}
solution to \eqref{e:1.1}. Hence we can obtain such solutions as
limits of solutions to \eqref{e:1.4} with a suitable sequence
$E_\varepsilon$ with vanishing cost.  On the other hand, if we set in
\eqref{e:1.4} $E=0$ we obtain, in the limit $\varepsilon \to 0$, an
entropic solution to \eqref{e:1.1}. If the flux $f$ is nonlinear, the
set of measure-valued solutions to \eqref{e:1.1} is larger than the
set of entropic solutions; it is thus natural to study the
$\Gamma$-convergence of the rescaled cost functional $H_\varepsilon :=
\varepsilon^{-1} I_\varepsilon$, which formally corresponds to the
scaling in \cite{Je,Va}.  Our second result concerns the
$\Gamma$-convergence of $H_\varepsilon$ which is studied under the
additional hypotheses that the flux $f$ is smooth and such that there
are no intervals in which $f$ is affine. A compensated compactness
argument shows that $H_\varepsilon$ has enough coercivity properties
to force its convergence in a functions setting and not in a Young
measures' one.

To informally define the candidate $\Gamma$-limit of $H_\varepsilon$,
we first introduce some preliminary notions. We say that a weak
solution $u$ to \eqref{e:1.1} is \emph{entropy-measure} iff for each
\emph{smooth} entropy $\eta$ the distribution $\eta(u)_t+q(u)_x$ is a
Radon measure on $(0,T)\times \bb R$. If $u$ is an entropy-measure
solution to \eqref{e:1.1}, then there exists a measurable map
$\varrho_u$ from $[0,1]$ to the set of Radon measures on $(0,T)\times
\bb R$, such that for each $\eta \in C^2([0,1])$ and $\varphi \in
C^\infty_{\mathrm{c}}\big((0,T)\times \bb R \big)$,
$-\int\!dt\,dx\,\big[\eta(u) \varphi_t+ q(u)\varphi_x\big] =
\int\!dv\,\varrho_u(v;dt,dx) \eta''(v) \varphi(t,x)$, see
Proposition~\ref{p:kin}. The candidate $\Gamma$-limit of
$H_\varepsilon$ is the functional $H$ defined as follows. If $u$ is
not an entropy-measure solution to \eqref{e:1.1} then $H(u)=+\infty$.
Otherwise $H(u)=\int\,dv\,\varrho_u^+(v;dt,dx) D(v)/\sigma(v)$, where
$\varrho_u^+$ denotes the positive part of $\varrho_u$. Note that
while $I_{\varepsilon}$ and $\mc I$ are nonlocal functionals, $H$ is
local. On the other hand, while $I_\varepsilon$, resp.\ $\mc I$,
quantifies in a suitable squared Hilbert norm the violation of
equation \eqref{e:1.2}, resp.\ \eqref{e:1.1}, this quadratic structure
is lost in $H$. In Proposition~\ref{p:H} we show that $H$ is a
coercive lower semicontinuous functional, this matching the necessary
properties for being the $\Gamma$-limit of a sequence of equicoercive
functionals. Note also that $H$ depends on the diffusion coefficient
$D$ and the conductivity coefficient $\sigma$ only through their
ratio, which is an expected property of well-behaving driven diffusive
systems, in hydrodynamical-like limits. We discuss this issue in
Remark~\ref{r:JV}, where a link between the functional $H$ and the
large deviations rate functional introduced in \cite{Je,Va} is also
investigated. In particular, $H$ comes as a natural generalization of
the functional introduced in \cite{Je,Va}, whenever the flux $f$ is
neither convex nor concave.

In this paper we prove that for each sequence $u^\varepsilon \to u$ in
$L_{1,\mathrm{loc}}([0,T]\times \bb R)$ we have
$\varliminf_\varepsilon H_\varepsilon(u^\varepsilon) \ge H(u)$, namely
$\Gliminf H_\varepsilon \ge H$.  Since the functional $H$ vanishes
only on entropic solutions to \eqref{e:1.1}, its zero-level set
coincides with the limit points of the minima of
$I_\varepsilon$. Concerning the $\Gamma$-limsup inequality, for each
weak solution $u$ to \eqref{e:1.1} in a suitable set $\mc S_\sigma$,
see Definition~\ref{d:splittable}, we construct a sequence
$u^\varepsilon \to u$ such that $H_\varepsilon(u^\varepsilon) \to
H(u)$. The above statements imply $(\Glim H_\varepsilon)(u) = H(u)$
for $u \in \mc S_\sigma$. To complete the proof of the
$\Gamma$-convergence of $H_\varepsilon$ to $H$ on the whole set of
entropy-measure solutions, an additional density argument is
needed. This seems to be a difficult problem, as Varadhan \cite{Va}
puts it: ``\ldots one does not see at the moment how to produce a
`general' non-entropic solution, partly because one does not know what
it is.''

The above results imply that if $u^\varepsilon$ solves
\eqref{e:1.4} for some control $E^\varepsilon$ such that
\begin{equation*}
\lim_{\varepsilon \to 0} \varepsilon^{-1} \int_{[0,T]}\!dt\,\|E^\varepsilon\|_{L_2(\bb
  R,\sigma(u^\varepsilon)dx)}^2=0
\end{equation*}
then any limit point of $u^\varepsilon$ is an entropic solution to
\eqref{e:1.1}. This statement is sharp in the sense that there are
sequences $\{E^\varepsilon\}$ with $\varliminf_{\varepsilon}
\varepsilon^{-1} \int_{[0,T]}\!dt\,\|E^\varepsilon\|_{L_2(\bb
  R,\sigma(u^\varepsilon)dx)}^2 > 0$ such that any limit point of the
corresponding $u^\varepsilon$ is not an entropic solutions to
\eqref{e:1.1}. More generally, the variational description of
conservation laws here introduced allows the following point of view.
Measure-valued solutions to \eqref{e:1.1} are the points in the
zero-level set of the $\Gamma$-limit of $I_\varepsilon$, while
entropic weak solutions are the points in the zero-level set of the
$\Gamma$-limit of $\varepsilon^{-1} I_\varepsilon$. In
Appendix~\ref{s:B} we introduce a sequence $\{J_\varepsilon\}$ of
functionals related to the viscous approximation of Hamilton-Jacobi
equations. In \cite{Po} a $\Gamma$-limsup inequality for a related
family of functionals has been independently investigated in a BV
setting. Following closely the proofs of the $\Gamma$-convergence of
$\{I_\varepsilon\}$, we establish the corresponding
$\Gamma$-convergence results, thus obtaining a variational
characterization of measure-valued and viscosity solutions to
Hamilton-Jacobi equations. Although this ``variational'' point of view
is consistent with the standard concepts of solution in the current
setting of scalar conservation laws and Hamilton-Jacobi equations, it
might be helpful for less understood model equations.

\section{Notation and results}
\label{s:2}
Hereafter in this paper, we assume that $f$ is a Lipschitz function on
$[0,1]$, $D$ and $\sigma$ are continuous functions on $[0,1]$, with
$D$ uniformly positive and $\sigma$ strictly positive on $(0,1)$. We
understand that these assumptions are supposed to hold in every
statement below.

We also let $\langle \cdot,\cdot\rangle$ denote the inner product in
$L_2(\bb R)$, for $T>0$ $\langle\langle \cdot,\cdot\rangle\rangle$
stands for the inner product in $L_2([0,T]\times\bb R)$, and for $O$
an open subset of $\bb R^n$, $C^\infty_{\mathrm{c}}(O)$ denotes the
collection of compactly supported infinitely differentiable functions
on $O$.

\medskip
\noindent\textit{Scalar conservation law} 

Our analysis will be restricted to equibounded densities $u$ that take
values in $[0,1]$. Let $U$ denote the compact separable metric space
of measurable functions $u:\bb R\to [0,1]$, equipped with the
following $H^{-1}_{\mathrm{loc}}$-like metric $d_U$.  For $L>0$, set
\begin{equation*}
\|u\|_{-1,L}:=\sup \big\{
       \langle u,\varphi\rangle,\,\varphi \in 
       C^\infty_{\mathrm{c}}\big((-L,L)\big),\, 
           \langle \varphi_x, \varphi_x\rangle =1 
            \big\}
\end{equation*}
and define the metric $d_U$ in $U$ by
\begin{equation}
  \label{e:2.1}
 d_U (u,v) := \sum_{N=1}^\infty 2^{-N} 
     \frac{\|u-v\|_{-1,N}}{1+\|u-v\|_{-1,N}}
\end{equation}
Given $T>0$, let $\mc U$ be the set $C\big([0,T]; U\big)$ endowed with
the uniform metric
\begin{equation}
\label{e:2.2}
  d_{\mc U}(u,v):= \sup_{t\in[0,T]} d_U\big( u(t),v(t) \big)
\end{equation}
An element $u\in \mc U$ is a \emph{weak solution} to \eqref{e:1.1} iff
for each $\varphi\in C^\infty_{\mathrm{c}}\big((0,T)\times \bb R\big)$
(in particular $\varphi(0)=\varphi(T)=0$) it satisfies
\begin{equation*}
  \langle\langle u,\varphi_t\rangle\rangle 
  +\langle\langle f(u),\varphi_x\rangle\rangle = 0
\end{equation*}

We also introduce a suitable space $\mc M$ of Young measures and
recall the notion of measure-valued solution to \eqref{e:1.1}.
Consider the set $\mc N$ of measurable maps $\mu$ from $[0,T]\times
\bb R$ to the set $\mc P([0,1])$ of Borel probability measures on
$[0,1]$. The set $\mc N$ can be identified with the set of positive
Radon measures $\mu$ on $[0,1]\times [0,T]\times \bb R$ such that
$\mu([0,1],\,dt,\, dx)=dt\,dx$. Indeed, by existence of a regular
version of conditional probabilities, for such measures $\mu$ there
exists a measurable kernel $\mu_{t,x}(d\lambda) \in \mc P([0,1])$ such
that $\mu(d\lambda,\,dt,\,dx)=dt\,dx\,\mu_{t,x}(d\lambda)$. For
$\imath : [0,1]\to [0,1]$ the identity map, we set
\begin{equation}
\label{e:2.3}
 \mc M := \big\{\mu \in \mc N\,:\:\text{the map
   $[0,T]\ni t \mapsto \mu_{t,\cdot}(\imath)$ is in 
             $\mc U$} \big\}
\end{equation}
in which, for a bounded measurable function $F:[0,1]\to \bb R$, the
notation $\mu_{t,x}(F)$ stands for $\int_{[0,1]} \! \mu_{t,x}(d\lambda)
F(\lambda)$. We endow $\mc M$ with the metric
\begin{equation}
  \label{e:2.4} 
  d_{\mc M} (\mu,\nu) := d_{\mathrm{w}}(\mu,\nu) 
   +d_{\mc U}\big(\mu(\imath),\nu(\imath)\big)
\end{equation}
where $d_{\mathrm{w}}$ is a distance generating the relative
topology on $\mc N$ regarded as a subset of the Radon measures
on $[0,1]\times [0,T]\times \bb R$ equipped with the weak* topology.
$(\mc M, d_{\mc M})$ is a complete separable metric space.

An element $\mu\in \mc M$ is a \emph{measure-valued solution} to
\eqref{e:1.1} iff for each $\varphi\in
C^\infty_{\mathrm{c}}((0,T)\times \bb R)$ it satisfies
\begin{equation*}
  \langle\langle \mu(\imath), \varphi_t \rangle\rangle
  + \langle\langle \mu(f), \varphi_x \rangle\rangle =0
\end{equation*}
If $u\in \mc U$ is a weak solution to \eqref{e:1.1}, then
$\delta_{u(t,x)}(d\lambda)\in \mc M$ is a measure-valued solution. On
the other hand, there exist measure-valued solutions which do not have
this form.

\medskip
\noindent\textit{Parabolic cost functional}

We next give the definition of the parabolic cost functional
informally introduced in \eqref{e:1.5}. Given $u \in \mc U$ we write
$u_x \in L_{2,\mathrm{loc}}([0,T]\times \bb R)$ iff $u$ admits a
locally square integrable weak $x$-derivative. For $\varepsilon>0$, $u
\in \mc U$ such that $u_x \in L_{2,\mathrm{loc}}([0,T]\times \bb R)$,
and $\varphi \in C^\infty_{\mathrm{c}}((0,T)\times \bb R)$ we set
\begin{equation}
\label{e:2.5}
\ell^u_\varepsilon(\varphi) :=
           - \langle \langle  u,\varphi_t\rangle \rangle
        -\langle \langle f(u),\varphi_x\rangle \rangle
        + \frac{\varepsilon}2 \langle \langle D(u) u_x, 
                           \varphi_x\rangle \rangle
\end{equation}
and define $I_\varepsilon: \mc U \to [0,+\infty]$ as follows.
If $u_x \in L_{2,\mathrm{loc}}([0,T]\times \bb R)$ we set 
\begin{equation}
\label{e:2.6}
I_\varepsilon (u):=
 \sup_{\varphi\in C^\infty_{\mathrm{c}}((0,T)\times \bb R)} 
  \Big[ 
   \ell_\varepsilon^u(\varphi)
   - \frac 12 \langle\langle \sigma(u)\varphi_x 
   \Big]
\end{equation}
letting $I_\varepsilon (u):= +\infty$ otherwise.
Note that $I_\varepsilon(u)$ vanishes iff $u \in \mc U$ is a weak
solution to \eqref{e:1.2}; more generally, by Riesz representation
theorem, it is not difficult to prove the connection of
$I_\varepsilon$ with the perturbed parabolic problem \eqref{e:1.4},
see Lemma~\ref{l:riesz} below for the precise statement.

In order to discuss the behavior of $I_\varepsilon$ as $\varepsilon\to
0$, we lift it to the space of Young measures $(\mc M, d_{\mc M})$.
We thus define $\mc I_\varepsilon:\mc M
\to [0,+\infty]$ by
\begin{equation}
\label{e:2.7}
   \mc I_\varepsilon(\mu):=
   \begin{cases}
     I_\varepsilon(u) & \text{ if } \quad
     \mu_{t,x} = \delta_{u(t,x)} \;\text{ for some }\; u\in\mc U \\
     +\infty  & \text{ otherwise }
   \end{cases}
\end{equation}

\medskip
\noindent\textit{Asymptotic parabolic cost}

As well known, a most useful notion of variational convergence is the
$\Gamma$-convergence which, together with some compactness
estimates, implies convergence of the minima. Let $X$ be a complete
separable metrizable space; recall that a sequence of functionals
$F_\varepsilon:X\to[-\infty,+\infty]$ is \emph{equicoercive} on $X$
iff for each $M >0$ there exists a compact set $K_M$ such that for any
$\varepsilon \in (0,1]$ we have $\{x \in X \,:\: F_\varepsilon(x) \le
M \} \subset K_M$. We briefly recall the basic definitions of the
$\Gamma$-convergence theory, see e.g.\ \cite{Br,DM}. Given $x \in X$
we define
\begin{eqnarray*}
\big( \Gliminf_{\varepsilon \to 0} F_\varepsilon \big)\, (x) &:= &
    \inf \big\{ \varliminf_{\varepsilon \to 0} F_\varepsilon(x^\varepsilon),\,
                   \{x^\varepsilon\} \subset X\,:\: x^\varepsilon \to x \big\}
\\
\big(\Glimsup_{\varepsilon \to 0} F_\varepsilon \big)\, (x) &:= & 
\inf \big\{\varlimsup_{\varepsilon \to 0}  F_\varepsilon(x^\varepsilon),\,
                   \{x^\varepsilon\} \subset X\,:\: x^\varepsilon \to x \big\}
\end{eqnarray*}
Whenever $\Gliminf_\varepsilon F_\varepsilon 
=\Glimsup_\varepsilon F_\varepsilon=F$ we say
that $F_\varepsilon$ $\Gamma$-converges to $F$ in $X$. Equivalently,
$F_\varepsilon$ $\Gamma$-converges to $F$ iff for each $x\in X$ we
have:
\begin{itemize}
\item[\rm{--}]{for any sequence $x^\varepsilon\to x$ we have
$\varliminf_\varepsilon F_\varepsilon(x^\varepsilon)\ge F(x)$ \
(\emph{$\Gamma$-liminf inequality});}
\item[\rm{--}]{ there exists a sequence $x^\varepsilon\to x$ such that
$\varlimsup_\varepsilon F_\varepsilon(x^\varepsilon)\le F(x)$ \
(\emph{$\Gamma$-limsup inequality}).}
\end{itemize}
Equicoercivity and $\Gamma$-convergence of a sequence
$\{F_\varepsilon\}$ imply an upper bound of infima over open sets, and
a lower bound of infima over closed sets, see e.g.\
\cite[Prop.~1.18]{Br}, and therefore it is the relevant notion of
variational convergence in the control setting introduced above.

\begin{theorem}
\label{t:2.1}
The sequence $\{\mc I_\varepsilon\}$ defined in \eqref{e:2.6},
\eqref{e:2.7} is equicoercive on $\mc M$ and, as $\varepsilon \to 0$,
$\Gamma$-converges in $\mc M$ to
\begin{equation}
  \label{e:2.8}
  \mc I (\mu) :=\sup_{\varphi \in C^\infty_{\mathrm{c}}((0,T) \times\bb R) } \Big\{ 
   - \langle\langle \mu(\imath), \varphi_t \rangle\rangle 
   - \langle\langle \mu (f), \varphi_x \rangle\rangle 
   - \frac 12\, \langle\langle \mu (\sigma) 
                  \varphi_x,\varphi_x \rangle\rangle 
     \Big\}
\end{equation}
\end{theorem}

Note that $\mc I (\mu)=0$ iff $\mu$ is a measure-valued solution to
\eqref{e:1.1}.
From Theorem~\ref{t:2.1} we deduce the $\Gamma$-limit of $I_\varepsilon$,
see \eqref{e:2.6}, on $\mc U$ by projection.

\begin{corollary2}
\label{c:rely}
The sequence of functionals $\{I_\varepsilon\}$ is equicoercive on $\mc
U$ and, as $\varepsilon \to 0$, $\Gamma$-converges in $\mc U$ to the
functional $I:\mc U \to [0,+\infty]$ defined by
\begin{equation*}
\begin{array}{lcl}
& & 
{\displaystyle
 I(u)  := 
  \inf \Big\{  
  \int \!dt \, dx\, 
  R_{f,\sigma}\big(u(t,x),\Phi(t,x)\big),
}
\\ & & \phantom{I(u)  := \inf \Big\{ }
   \Phi\in L_{2,\mathrm{loc}}([0,T]\times \bb R) \,:\: 
    \Phi_x=-u_t \; \mathrm{ weakly}
  \Big\}
\end{array}
\end{equation*}
where $R_{f,\sigma}:[0,1]\times \bb R \to [0,+\infty]$ is defined by
\begin{equation*}
  R_{f,\sigma}(w,c):=\inf \{\big(\nu(f)-c\big)^2 / \nu(\sigma),
            \,\nu   \in \mc P([0,1])\,:\: \nu(\imath)=w \}
\end{equation*}
in which we understand $(c-c)^2/0=0$.
\end{corollary2}
From the proof of Corollary~\ref{c:rely} it follows $I(\cdot) \le \mc
I (\delta_{\cdot})$, and the equality holds iff $f$ is linear. If we
restrict to \emph{stationary} $u$'s, namely to the case $u_t=0$,
Corollary~\ref{c:rely} can be regarded as a negative-Sobolev version
of classical relaxation results for integral functionals in weak
topology. More precisely, from the proofs of Theorem~\ref{t:lsce} and
Corollary~\ref{c:rely} it follows that if we define the functional
$\tilde{F}:U\to [0,+\infty]$ by
\begin{equation*}
\tilde{F}(u):= \inf_{c \in \bb R} \int \! dx\, 
                  \frac{\big[f(u(x))-c \big]^2}{\sigma(u(x))}
\end{equation*}
then its lower semicontinuous envelope w.r.t.\ the $d_{U}$-distance
\eqref{e:2.1} is given by
\begin{equation*}
F(u):= \inf_{c \in \bb R} \int \! dx \,
                      R_{f,\sigma}(u(x),c)
\end{equation*}
Note also that $R_{f,\sigma}$ can be explicitly calculated in some
cases. Let $\underline{f}, \overline{f}:[0,1] \to \bb R$ be
respectively the convex and concave envelope of $f$. Then, in the case
$\sigma=1$, we have $R_{f,1}(w,c) =
\big[\mathrm{distance}(c,[\underline{f}(w),\overline{f}(w)])\big]^2$. In
the case $f=\sigma$ (which includes the example mentioned in the
introduction $f(u)=\sigma(u)=u(1-u)$) then
\begin{equation*}
R_{f,f}(w,c)=
\begin{cases}
2 (|c|-c) &
      \text{if $ |c| \in [\underline{f}(w),\overline{f}(w)]$}
\\ \frac{(\overline{f}(w)-c)^2}{\overline{f}(w)} &
                           \text{if $|c|>\overline{f}(w)$}
\\ \frac{(\underline{f}(w)-c)^2}{\underline{f}(w)} &
     \text{if $|c|<\underline{f}(w)$}
\end{cases}
\end{equation*}

\medskip 
\noindent\textit{Entropy-measure solutions} 

Recalling \eqref{e:2.2}, we let $\mc X$ be the same set $C([0,T];U)$
endowed with the metric
\begin{equation}
  \label{e:2.9}
  d_{\mc X}(u,v):=  
      \sum_{N=1}^{\infty} \frac 1{2^N}  \|u-v\|_{L_1([0,T]\times[-N,N])}
     +d_{\mc U}(u, v)
\end{equation}
Convergence in $\mc X$ is equivalent to convergence in $\mc U$
\emph{and} in $L_{p,\mathrm{loc}}([0,T]\times \bb R)$ for $p\in
[1,+\infty)$.

Let $C^2([0,1])$ be the set of twice differentiable functions on
$(0,1)$ whose derivatives are continuous up to the boundary. A
function, resp. a convex function, $\eta \in C^2([0,1])$ is called an
\emph{entropy}, resp. a \emph{convex entropy}, and its
\emph{conjugated entropy flux} $q\in C([0,1])$ is defined up to a
constant by $q(u):=\int^u\!dv\,\eta'(v)f'(v)$. For $u$ a weak solution
to \eqref{e:1.1}, for $(\eta,q)$ an entropy -- entropy flux pair, the
\emph{$\eta$-entropy production} is the distribution $\wp_{\eta,u}$
acting on $C^\infty_{\mathrm{c}}\big((0,T)\times\bb R\big)$ as
\begin{equation}
  \label{e:2.10}
\wp_{\eta,u}(\varphi):= 
    - \langle\langle \eta(u) , \varphi_t \rangle\rangle
    - \langle\langle q(u) , \varphi_x \rangle\rangle
\end{equation}

Let $C^{2,\infty}_{\mathrm{c}}\big([0,1]\times (0,T)\times \bb R
\big)$ be the set of compactly supported maps
$\vartheta:[0,1]\times(0,T)\times \bb R \ni (v,t,x) \mapsto
\vartheta(v,t,x) \in \bb R$, that are twice differentiable in the $v$
variable, with derivatives continuous up to the boundary of
$[0,1]\times(0,T)\times \bb R$, and that are infinitely differentiable
in the $(t,x)$ variables. For $\vartheta \in
C^{2,\infty}_{\mathrm{c}}\big([0,1]\times (0,T)\times \bb R \big)$ we
denote by $\vartheta'$ and $\vartheta''$ its partial derivatives
w.r.t.\ the $v$ variable. We say that a function $\vartheta \in
C^{2,\infty}_{\mathrm{c}}\big([0,1]\times (0,T)\times \bb R \big)$ is
an \emph{entropy sampler}, and its \emph{conjugated entropy flux
  sampler} $Q:[0,1]\times (0,T)\times \bb R$ is defined up to an
additive function of $(t,x)$ by
$Q(u,t,x):=\int^u\!dv\,\vartheta'(v,t,x) f'(v)$. Finally, given a weak
solution $u$ to \eqref{e:1.1}, the \emph{$\vartheta$-sampled entropy
  production} $P_{\vartheta,u}$ is the real number
\begin{equation}
\label{e:2.11}
P_{\vartheta,u}
:=-\int \! dt\,dx\,
  \Big[\big(\partial_t \vartheta)\big(u(t,x),t,x \big) 
       + \big(\partial_x Q\big)\big(u(t,x),t,x \big)\Big]
\end{equation}
If $\vartheta(v,t,x)=\eta(v) \varphi(t,x)$ for some entropy $\eta$ and
some $\varphi \in C^{\infty}_{\mathrm{c}}\big((0,T)\times \bb R\big)$,
then $P_{\vartheta,u}=\wp_{\eta,u}(\varphi)$.

The next proposition introduces a suitable class of solutions to
\eqref{e:1.1} which will be needed in the following. We denote by
$M\big((0,T)\times \bb R\big)$ the set of Radon measures on $(0,T)
\times \bb R$ that we consider equipped with the weak* topology. In
the following, for $\varrho \in M\big((0,T)\times \bb R\big)$ we
denote by $\varrho^\pm$ the positive and negative part of
$\varrho$. For $u$ a weak solution to \eqref{e:1.1} and $\eta$ an
entropy, recalling \eqref{e:2.10} we set
\begin{equation}
\label{e:2.12}
\|\wp_{\eta,u}\|_{\mathrm{TV},L}:=
    \sup \big\{\wp_{\eta,u}(\varphi),\,\varphi \in
               C^\infty_{\mathrm{c}}\big((0,T)\times (-L,L) \big),\,
               |\varphi|\le 1
         \big\} 
\end{equation}
\begin{equation*}
\|\wp_{\eta,u}^+\|_{\mathrm{TV},L}:= 
    \sup \big\{\wp_{\eta,u}(\varphi),\,\varphi \in
               C^\infty_{\mathrm{c}}\big((0,T)\times (-L,L) \big),\,
               0 \le \varphi \le 1
         \big\}
\end{equation*}

\begin{proposition2}
\label{p:kin}
Let $u \in \mc X$ be a weak solution to \eqref{e:1.1}. The following
statements are equivalent:
\begin{itemize}
\item[{\rm (i)}] {there exists $c>0$ such that
    $\|\wp_{\eta,u}^+\|_{\mathrm{TV},L}<+\infty$ for any $L>0$ and
    $\eta \in C^2([0,1])$ with $0 \le \eta'' \le c$;}
\item[{\rm (ii)}]{for each entropy $\eta$, the $\eta$-entropy
    production $\wp_{\eta,u}$ can be extended to a Radon measure on
    $(0,T)\times \bb R$, namely
    $\|\wp_{\eta,u}\|_{\mathrm{TV},L}<+\infty$ for each $L>0$;}
\item[{\rm (iii)}]{there exists a bounded measurable map
    $\varrho_u:[0,1] \ni v \to \varrho_u(v;dt,dx) \in M\big((0,T)\times
    \bb R\big)$ such that for any entropy sampler $\vartheta$
\begin{equation}
\label{e:2.13}
P_{\vartheta,u} = \int \! dv\,\varrho_u(v;dt,dx)\,\vartheta''(v,t,x)
\end{equation}
}
\end{itemize}
\end{proposition2}

A weak solution $u \in \mc X$ that satisfies any of the equivalent
conditions in Proposition~\ref{p:kin} is called an
\emph{entropy-measure solution} to \eqref{e:1.1}. We denote by $\mc E
\subset \mc X$ the set of entropy-measure solutions to \eqref{e:1.1}.
Proposition~\ref{p:kin} establishes a so-called \emph{kinetic
  formulation} for entropy-measure solutions, see also
\cite[Prop.~3.1]{DOW} for a similar result. If $f \in C^2([0,1])$ is
such that there are no intervals in which $f$ is affine, using the
results in \cite{CR} we show that entropy-measure solutions have some
regularity properties, see Lemma~\ref{l:l1}.

A weak solution $u\in \mc X$ to \eqref{e:1.1} is called an
\emph{entropic solution} iff for each convex entropy $\eta$ the
inequality $\wp_{\eta,u} \le 0$ holds in distribution sense, namely
$\|\wp_{\eta,u}^+\|_{\mathrm{TV},L}=0$ for each $L>0$. In particular
entropic solutions are entropy-measure solutions such that
$\varrho_u(v;dt,dx)$ is a negative Radon measure for each $v \in
[0,1]$.  It is well known, see e.g.\ \cite{Da,Se}, that for each $u_0
\in U$ there exists a unique entropic solution $\bar u\in C([0,T];
L_{1,\mathrm{loc}}(\bb R))$ to \eqref{e:1.1} such that $\bar
u(0)=u_0$. Such a solution $\bar{u}$ is called the \emph{Kruzkov
  solution} with initial datum $u_0$.

\medskip
\noindent\emph{$\Gamma$-entropy cost of non-entropic solutions} 

We next introduce a rescaled cost functional and prove in particular
that entropic solutions are the only ones with vanishing rescaled
asymptotic cost. Recalling that $I_\varepsilon$ has been introduced in
\eqref{e:2.6}, the \emph{rescaled cost functional} $H_\varepsilon: \mc X
\to [0,+\infty]$ is defined by
\begin{equation}
\label{e:2.14}
H_\varepsilon(u) := \varepsilon^{-1}I_\varepsilon(u)
\end{equation}
In the $\Gamma$-convergence theory, the asymptotic behavior of the
rescaled functional $H_\varepsilon$ is usually referred to as the
development by $\Gamma$-convergence of $I_\varepsilon$, see e.g.\
\cite[\S 1.10]{Br}. In our case, while we lifted $I_\varepsilon$ to the
space of Young measures $\mc M$, we can consider the rescaled cost
functional $H_\varepsilon$ on $\mc X$. In fact, as shown below,
$H_\varepsilon$ has much better compactness properties than $I_\varepsilon$
and it is equicoercive on $\mc X$. Therefore the $\Gamma$-convergence
of the lift of $H_\varepsilon$ to $\mc M$ can be immediately retrieved
from the $\Gamma$-convergence of $H_\varepsilon$ on $\mc X$. Indeed,
since $\delta_{u_\varepsilon}\to \delta_u$ in $\mc M$ implies
$u_\varepsilon\to u$ in $\mc X$, the metric \eqref{e:2.9} generates the
relative topology of $\mc X$ regarded as a subset of $\mc M$.

Recall that $\mc E \subset \mc X$ denotes the set of entropy-measure
solutions to \eqref{e:1.1}, and that for $u \in \mc E$ there exists a
bounded measurable map $\varrho_u:[0,1]\to M\big((0,T)\times \bb
R\big)$ such that \eqref{e:2.13} holds. Let $\varrho_u^+$ be the
positive part of $\varrho_u$, and define $H:\mc X \to [0,+\infty]$ by
\begin{equation}
  \label{e:2.15}
  H(u):= 
  \begin{cases}
\displaystyle{
   \int \! dv\,\varrho_u^+(v;dt,dx)\, \frac{D(v)}{\sigma(v)}
             }
           & \text{if $u\in \mc E$}
 \\
+ \infty   & \text{ otherwise }
  \end{cases}
\end{equation}
As shown in the proof of Theorem~\ref{t:ecne}, if $u$ is a weak
solution to \eqref{e:1.1} and $H(u)<+\infty$, then
$H(u)=\sup_{\vartheta} P_{\vartheta,u}$, where the supremum is taken
over the entropy samplers $\vartheta$ such that $0 \le \sigma(v)
\vartheta''(v,t,x) \le D(v)$, for each $(v,t,x) \in [0,1]\times
[0,T]\times \bb R$.

\begin{definition2}
\label{d:splittable}
An entropy-measure solution $u \in \mc E$ is \emph{entropy-splittable}
iff there exist two closed sets $E^+,E^- \subset [0,T]\times \bb R$
such that
\begin{itemize}
\item[{\rm (i)}] {For a.e.\ $v \in [0,1]$, the support of
  $\varrho_u^+(v;dt,dx)$ is contained in $E^+$, and the support of
  $\varrho_u^-(v;dt,dx)$ is contained in $E^-$.}
\item[{\rm (ii)}] {For each $L>0$, the set $\Big\{ t \in [0,T]\,:\:
    \big(\{t\} \times [-L,L]\big) \cap E^+ \cap E^- \neq \emptyset
    \Big\}$ is nowhere dense in $[0,T]$.}
\end{itemize}
The set of entropy-splittable solutions to \eqref{e:1.1} is denoted by
$\mc S$. An entropy-splittable solution $u \in \mc S$ such that
$H(u)<+\infty$ and
\begin{itemize}
\item[{\rm (iii)}] {For each $L>0$ there exists $\delta_L>0$ such that
  $\sigma(u(t,x)) \ge \delta_L$ for a.e.\ $(t,x) \in [0,T]\times \bb
  [-L,L]$.}
\end{itemize}
is called \emph{nice w.r.t.\ $\sigma$}. The set of nice (w.r.t.\
$\sigma$) solutions to \eqref{e:1.1} is denoted by $\mc S_\sigma$.
\end{definition2}
Note that $\mc S_\sigma \subset \mc S \subset \mc E \subset \mc X$,
and that, if $\sigma$ is uniformly positive on $[0,1]$, then $\mc
S_\sigma=\mc S$. In Remark~\ref{r:splittable} we exhibit a few
classes of entropy-splittable solutions to \eqref{e:1.1}.

\begin{theorem}
 \label{t:ecne} 
Let $H_\varepsilon$ and $H$ be the functionals on $\mc X$ as
respectively defined in \eqref{e:2.14} and \eqref{e:2.15}. 
\begin{itemize}
\item[{\rm (i)}] {The sequence of functionals $\{H_\varepsilon\}$
    satisfies the $\Gamma$-liminf inequality
    $\Gamma$\textrm{-}$\varliminf_\varepsilon H_\varepsilon \ge H$ on $\mc X$.}

\item[{\rm (ii)}] {Assume that there is no interval where $f$ is
  affine. Then the sequence of functionals $\{H_\varepsilon\}$ is
  equicoercive on $\mc X$.}

\item[{\rm (iii)}] {Assume furthermore that $f \in C^2([0,1])$, and
  $D,\sigma \in C^{\alpha}([0,1])$ for some $\alpha>1/2$. Define
\begin{equation*}
  \upbar{H}(u) := 
            \inf \big\{\varliminf H(u_n),\,
               \{u_n\} \subset \mc S_\sigma\,:\:
                 u_n \to u \text{ in $\mc X$} \big\}
\end{equation*} 
Then the sequence of functionals $\{H_\varepsilon\}$ satisfies the
$\Gamma$-limsup inequality $\Gamma$\textrm{-}$\varlimsup_\varepsilon
H_\varepsilon \le \upbar{H}$ on $\mc X$.  }
\end{itemize}
\end{theorem}

From the lower semicontinuity of $H$ on $\mc X$, see
Proposition~\ref{p:H}, it follows that $\upbar{H} \ge H$ on $\mc X$
and $\upbar{H}= H$ on $\mc S_\sigma$, namely the $\Gamma$-convergence
of $H_\varepsilon$ to $H$ holds on $\mc S_\sigma$. To get the full
$\Gamma$-convergence on $\mc X$, the inequality $H(u)\ge \upbar{H}(u)$
is required also for $u \not\in \mc S_\sigma$. This amounts to show
that $\mc S_\sigma$ is $H$-dense in $\mc X$, namely that for $u \in
\mc X$ such that $H(u)<+\infty$ there exists a sequence $\{u^n\}
\subset \mc S_\sigma$ converging to $u$ in $\mc X$ such that
$H(u^n)\to H(u)$. As mentioned at the end of the introduction, this
appears to be a difficult problem. A preliminary step in this
direction is to obtain a chain rule formula for bounded vector fields
on $[0,T]\times \bb R$ the divergence of which is a Radon measure
(divergence-measure fields). This is a classical result for locally BV
fields \cite{AFP}. However, while entropic solutions to \eqref{e:1.1}
are in $BV_{\mathrm{loc}}([0,T]\times \bb R)$
\cite[Corollary~1.3]{AD} whenever $f$ is uniformly convex or concave, as shown
in Example~\ref{r:ex} below, the set $\{u \in \mc X\,:\:H(u)<+\infty\}$ is not
contained in $BV_{\mathrm{loc}}([0,T]\times \bb R)$ even under this assumptions
on $f$; see \cite{DW} for similar
examples including estimates in Besov norms. Chain rule formulae out
of the BV setting have been investigated in several recent papers; in
particular in \cite{DOW}, a chain rule formula for divergence-measure
fields is addressed, providing some partial results. In the remaining
of this section we discuss some properties of $H$, and some issues
related to the $H$-density of $\mc S_\sigma$.

In the following proposition we show that $H$ is lower semicontinuous,
and that it is coercive under the same hypotheses used for the
equicoercivity of $\{H_\varepsilon\}$. Moreover, we prove that the
minimizers of $H$ are limit points of the minimizers of
$I_\varepsilon$ as $\varepsilon \to 0$, so that no further rescaling
of $\{I_\varepsilon\}$ has to be investigated.

\begin{proposition2}
\label{p:H}
The functional $H$ is lower semicontinuous on $\mc X$ and $H(u)=0$ iff
$u$ is an entropic solution to \eqref{e:1.1}.  If furthermore there
are no intervals where $f$ is affine then $H$ is also coercive on $\mc
X$.
\end{proposition2}
From Proposition~\ref{l:l1} and the aforementioned regularity of entropy-measure solutions, see Lemma~\ref{l:l1}, 
it follows that if $f \in C^2([0,1])$ then the
zero-level set of $H$ coincides with the set of Kruzkov solutions to
\eqref{e:1.1}.

If $u \in \mc X$ is a weak solution with locally bounded variation,
Vol'pert chain rule, see \cite{AFP}, gives a formula for $H(u)$ in
terms of the normal traces of $u$ on its jump set, as shown in the following
remark.
\begin{remark2}
  \label{r:bv} 
  Let $u \in \mc X \cap BV_{\mathrm{loc}}([0,T]\times \bb R)$ be a
  weak solution to \eqref{e:1.1}. Denote by $J_u \subset [0,T]\times
  \bb R$ its jump set, by $\mc H^1 \res J_u$ the one-dimensional
  Hausdorff measure restricted to $J_u$, by $n=\big(n^{t},n^x\big)$ a
  unit normal to $J_u$ (which is well defined $\mc H^1 \res J_u$
  a.e.), and by $u^{\pm}$ the normal traces of $u$ on $J_u$ w.r.t.\
  $n$.  Then the Rankine-Hugoniot condition
  $(u^+-u^-)n^t+\big(f(u^+)-f(u^-)\big)n^x=0$ holds. In particular we
  can choose $n$ so that $n^x$ is uniformly positive, and thus $u^+$
  is the right trace of $u$ and $u^-$ is the left trace of $u$. Then
  $u \in \mc E$ and
\begin{equation*}
\varrho_u(v;dt,dx)= \frac{d \mc H^1 \res J_u}
             {\big\{(u^+-u^-)^2+[f(u^+)-f(u^-)]^2 \big\}^{1/2}}
             \rho(v,u^+,u^-)
\end{equation*} 
where, denoting by $u^- \wedge u^+$ and $u^- \lor u^+$ respectively
the minimum and maximum of $\{u^-,\,u^+\}$, $\rho: [0,1]^3 \to \bb R$
is defined by

\begin{eqnarray*}
\rho(v,u^+,u^-)&:=& \big[f(u^-)(u^+ -v) + f(u^+)(v-u^-)
                          -f(v)(u^+-u^-) \big]
\\                          
    && \times \id_{[u^- \wedge u^+,u^- \lor u^+]}(v)
\end{eqnarray*}
Hence, denoting by $\rho^+$ the positive part of $\rho$

\begin{eqnarray}
\label{e:2.16}
\nonumber
H(u)& =& \int_{J_u} \! \frac{d \mc H^1}
{\big\{(u^+-u^-)^2+[f(u^+)-f(u^-)]^2 \big\}^{1/2}}
\int \! dv\, \rho^+(v,u^+,u^-)\,\frac{D(v)}{\sigma(v)}
\\
& = & \int_{J_u} \! d \mc H^1 |n^x| 
   \int \! dv\, \frac{\rho^+(v,u^+,u^-)}{|u^+-u^-|} \frac{D(v)}{\sigma(v)}
\end{eqnarray}
\end{remark2}

Note $\rho(v,u^+,u^-) \le 0$ iff $\frac{f(v) -f(u^-)}{v-u^-} \ge
\frac{f(u^+)-f(v)}{u^+-v}$. This corresponds to the well known
geometrical secant condition for entropic solutions, see e.g.\
\cite{Da,Se}. Therefore $H(u)$ quantifies the violation of the entropy
condition along the non-entropic shocks of $u$.

In the following Example~\ref{r:ex} we show that neither the domain of
$H$, neither the $H$-closure of $\mc S_\sigma$ are contained in
$BV_{\mathrm{loc}}\big([0,T]\times \bb R \big)$.
\begin{figure}
\label{f:alb}
\begin{picture}(325,100)(0,0)
{
\thinlines
\put(0,0){\vector(1,0){300}}
\put(305,-3){$t=0$}
\put(0,100){\line(1,0){300}}
\put(305,97){$t=1$}
\put(99,-10){\footnotesize $0$}
\put(101.5,0){\line(-1,1){100}}
\put(45,25){{\footnotesize $\frac 12$}}
\put(140,-10){\footnotesize $b_1\!-\!b_2$}
\put(151.5,0){\line(-1,1){100}}
\put(49,60){{\footnotesize $\frac 12 +b_1$}}
\put(151.5,0){\line(-1,2){50}}
\put(90,75){{\footnotesize $\frac 12$}}
\put(178,-10){\footnotesize $b_1\!-\!b_3$}
\put(185.5,0){\line(-1,2){50}}
\put(125,60){{\footnotesize $\frac 12 +b_2$}}
\put(185.5,0){\line(-1,3){33.3}}
\put(148,80){{\footnotesize $\frac 12$}}
\put(219,-10){\footnotesize $b_1$}
\thicklines
\put(220,0){\line(0,1){100}}

\put(250,50){\footnotesize $\frac 12$}
\put(180,66){$\cdots$}
\put(190,33){$\cdots$}
\put(200,-2.9){$\cdots$}
}
\end{picture}
\smallskip
\caption{The values of $u$ in Example~\ref{r:ex} for $T=1$.}
\end{figure}
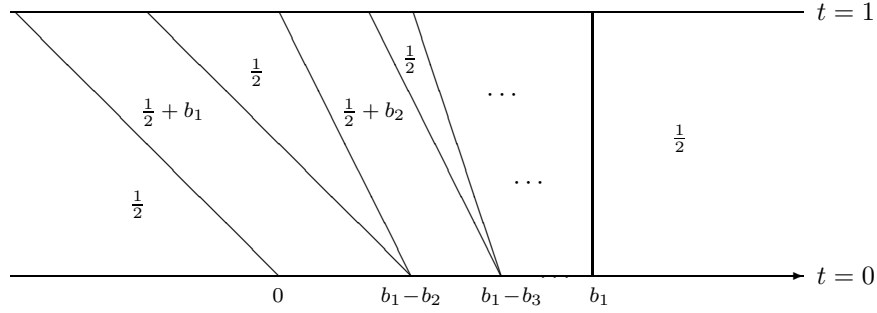

\begin{example2}
\label{r:ex}
Let $f(u)=u(1-u)$ and pick a decreasing sequence $\{b_i\}$ of positive
reals such that $b_1 <1/2$, $\sum_i b_i=+\infty$ and $\sum_i b_i^3
<+\infty$. Let $u$ be defined by
\begin{equation*}
u(t,x):=
\begin{cases}
1/2+b_i & 
        \text{if $T(b_1-b_i)< x+b_i\,t 
            <T\big(b_1-b_{i+1}\big)$ for some $i$}
\\
1/2    & \text{otherwise}
\end{cases}
\end{equation*}
Then $H(u)= \frac T2 \sum_i \int_{[0,b_i]}\!dv\,
\frac{D(1/2+v)}{\sigma(1/2+v)} v(b_i-v) <+\infty$. Note that, even if
the initial datum is in $BV(\bb R)$ and $f$ is concave, $u \not\in
BV_{\mathrm{loc}}([0,T]\times \bb R)$. However
$H(u)=\upbar{H}(u)$. Indeed the sequence $\{u^n\} \subset \mc
S_\sigma$ defined by
\begin{equation*}
u^n(t,x):=
\begin{cases}
u(t,x) & \text{if $x+b_n\,t < T(b_1-b_{n+1})$}
\\
1/2 & \text{otherwise}
\end{cases}
\end{equation*}
is such that $u^n \to u$ in $\mc X$ and $\lim_n H(u^n)=H(u)$.
\end{example2}

In the following remarks we identify some classes of
entropy-splittable solutions to \eqref{e:1.1}, see
Definition~\ref{d:splittable}.

\begin{remark2}
\label{r:splittable}
Weak solutions to \eqref{e:1.1} such that, for each convex entropy
$\eta$, $\wp_{\eta,u} \le 0$ (entropic solutions) or $\wp_{\eta,u} \ge
0$ (anti-entropic solutions) are entropy-splittable. Indeed they are
entropy-measure solutions (see Proposition~\ref{p:H}) and they fit in
Definition~\ref{d:splittable} with the choice $E^-=[0,T]\times \bb R$
and $E^+= \emptyset$ (for entropic solutions), and respectively
$E^+=[0,T]\times \bb R$ and $E^-=\emptyset$ (for anti-entropic
solutions).

Let $u \in BV_{\mathrm{loc}}\big([0,T]\times \bb R\big)$ be a weak
solution to \eqref{e:1.1}.  In the same setting of Remark~\ref{r:bv},
let us define $J_u^\pm:=\mathrm{Closure}\big(\{(t,x)\in J_u\,:\:
\exists v \in [0,1]\,:\:\pm \varrho(v;u^+,u^-) > 0\}\big)$. Suppose
that for each $L>0$ the set $\{t \in [0,T]\,:\: \big(\{t\}\times
[-L,L] \big) \cap J_u^+ \cap J_u^-\}$ is nowhere dense in
$[0,T]$. Then $u$ is an entropy-splittable solution. If $f$ is convex
or concave the sign of $\rho(v,u^+,u^-)$ does not depend on $v \in
[u^- \wedge u^+, u^-\vee u^+]$. Therefore, under this convexity
hypothesis, weak solutions to \eqref{e:1.1} with locally bounded
variations and with a jump set $J_u$ consisting of a locally finite
number of Lipschitz curves, intersecting each other at a locally
finite number of points are entropy splittable.
\end{remark2}

For a general (possibly neither convex nor concave) flux $f$, even
piecewise constant solutions to \eqref{e:1.1} may fail to be
entropy-splittable. However, in the following Example~\ref{r:ex2} we
introduce a family of weak solutions $u$ to \eqref{e:1.1} that are not
entropy-splittable, and show that they are in the $H$-closure of $\mc
S_\sigma$, and thus $\upbar{H}(u)=H(u)$. However, while
Example~\ref{r:ex2} can be widely generalized to prove
$\upbar{H}(u)=H(u)$ for $u$ in suitable classes of piecewise smooth
solutions, it does not seem that the ideas suggested by this example
may work in the general setting of entropy-measure solutions $u \in
\mc E$.
\begin{example2}
\label{r:ex2}
Let $\gamma:[0,T]\to \bb R$ be a Lipschitz map, let $u$ be a weak
solution of bounded variation to \eqref{e:1.1}, and suppose that the
jump set of $u$ coincides with $\gamma$. Let $u^-\equiv u^-(t)$ and
$u^+\equiv u^+(t)$ be the traces of $u$ on $\gamma$, and suppose that
there exists $u^0\in (0,1)$ such that $u^-(t)<u^0<u^+(t)$ for each $t$
and $\frac{f(v) -f(u^-)}{v-u^-} \ge \frac{f(u^+)-f(v)}{u^+-v}$ for $v
\in [u^-,u^0]$ and $\frac{f(v) -f(u^-)}{v-u^-} \le
\frac{f(u^+)-f(v)}{u^+-v}$ for $v \in[u^-,u^0]$. Then, if these
inequalities are strict at some $v$ and $t$, $u$ is not
entropy-splittable. However defining $u^n \in \mc X$ by
\begin{equation*}
u^n(t,x):=
\begin{cases}
u(t,x+n^{-1}) & \text{if $x \le \gamma(t)-n^{-1}$}
\\
u^0 & \text{if $\gamma(t)-n^{-1} < x < \gamma(t)+n^{-1} $}
\\
u(t,x-n^{-1}) &  \text{if $x \le \gamma(t)+n^{-1}$}
\end{cases}
\end{equation*}
we have that $u^n \in \mc S$, $u^n \to u$ in $\mc X$ and
$H(u^n)=H(u)$. In particular, if $\sigma(u)$ is uniformly positive on
compact subsets of $[0,T]\times \bb R$, then $\upbar{H}(u)=H(u)$. It
is easy to extend this example to the case in which the jump set of
$u$ consists of a locally finite number of Lipschitz curves
non-intersecting each other, provided that on each curve the quantity
$\frac{f(v) -f(u^-)}{v-u^-} - \frac{f(u^+)-f(v)}{u^+-v}$ changes its
sign a finite number of times for $v \in [u^+ \wedge u^-,u^+\vee
u^-]$.
\end{example2}

We next discuss the link between this paper and \cite{Je,Va}. In the
introduction we informally described the connection between the
problem \eqref{e:1.4} and stochastic particles systems under Euler
scaling. It is interesting to note that such a quantitative connection
can also be established for the limiting functionals. The key point is
that we expect the functional $H$ defined in \eqref{e:2.15} to
coincide with the large deviations rate functional introduced in
\cite{Je,Va}, provided the functions $f$, $D$ and $\sigma$ are chosen
correspondingly. Unfortunately, we cannot establish such an
identification off the set of weak solutions to \eqref{e:1.1} with
locally bounded variation.
\begin{remark2}
\label{r:JV}
Let $H^\prime:\mc X \to [0,+\infty]$ be defined as follows.
If $u \in \mc E$ we set 
\begin{equation*}
H^\prime(u):=
\sup \big\{ \| \wp_{\eta,u}^+\|_{\mathrm{TV},L},\,L>0,\, \eta \in
C^2([0,1])\,:\: 0\le \sigma\,\eta'' \le D \big\} 
\end{equation*}
letting $H^\prime(u):=+\infty$ otherwise. 
Then $H \ge H^\prime$ and $H(u)=H^\prime(u)$ whenever there exists a
Borel set $E^+ \subset [0,T]\times \bb R$ such that 
for a.e.\ $v\in [0,1]$ the measure
$\varrho_u^+(v;dt,dx)$ is concentrated on $E^+$ and
$\varrho_u^-(v;dt,dx)=0$ on $E^+$. In particular if
$f$ is convex or concave and $u \in BV_{\mathrm{loc}}([0,T]\times \bb
R)$, then $H(u)=H^\prime(u)$. If $f$ is neither convex nor concave,
then there exists $u \in \mc X$ such that $H(u)>H^\prime(u)$.
\end{remark2}

A general connection between dynamical transport coefficients and
thermodynamic potentials in driven diffusive systems is the so-called
\emph{Einstein relation}, see e.g.\ \cite[II.2.5]{Sp}. For a physical
model described by \eqref{e:1.4}, this relation states that the
\emph{Einstein entropy} $h \in C^2((0,1))\cap C([0,1];[0,+\infty])$
defined by
\begin{equation*}
\sigma(v) h''(v)=D(v)\qquad v \in (0,1)
\end{equation*}
is a physically relevant entropy in the limit $\varepsilon \to 0$. 
We let $g$ be the conjugated flux to $h$, i.e.\
$g(u):=\int_{1/2}^u\!dv\, h'(v)f'(v)$.
Note that $h$, $g$ may be unbounded if $\sigma$ vanishes at the
boundary of $[0,1]$ and that $g\le C_1+C_2\,h$ for some constants
$C_1,C_2 \ge 0$. If $u$ is a weak solution to \eqref{e:1.1} such that
$h(u) \in L_{1,\mathrm{loc}}([0,T]\times \bb R)$ and such that the
distribution $h(u)_t+g(u)_x$ acts as a Radon measure on $(0,T)\times
\bb R$, we let $\|\wp_{h,u}^+\|_{\mathrm{TV}}$ be the total variation
of the positive part of such a measure. By monotone convergence
$H^\prime(u) \ge \|\wp_{h,u}^+\|_{\mathrm{TV}}$ for such a $u$, and if
$f$ is convex or concave and $u$ has locally bounded variation, then
indeed $H^\prime(u)=\|\wp_{h,u}^+\|_{\mathrm{TV}}$. If $f$ is convex or
concave, we do not know whether
$H(u)=H^\prime(u)=\|\wp_{h,u}^+\|_{\mathrm{TV}}$ for all $u \in \mc
X$, since a chain rule formula for divergence-measure fields is
missing.

The problem investigated in \cite{Je,Va} formally corresponds to the
case $f(u)=\sigma(u)=u(1-u)$ and $D(u)=1$, so that the Einstein
entropy $h$ coincides with the Bernoulli entropy $h(u)=-u\log u
-(1-u)\log(1-u)$. The (candidate) large deviations rate functional
$H^{JV}$ introduced in \cite{Je,Va} is defined as $+\infty$ off the
set of weak solutions to \eqref{e:1.1}, while
$H^{JV}(u)=\|\wp_{h,u}^+\|_{\mathrm{TV}}$ for $u$ a weak solution
(this is well defined, since $h$ is bounded). We thus have $H \ge
H^{JV}$, and in view of the $\Gamma$-liminf inequality, $H$ comes as a
natural generalization of $H^{JV}$ for diffusive systems with no
convexity assumptions on the the flux $f$.

\medskip
\noindent\emph{Outline of the proofs}

Standard parabolic a priori estimates on $u$ in terms of
$I_\varepsilon(u)$ imply equicoercivity of $\mc I_\varepsilon$ on $\mc M$.
Equicoercivity of $H_\varepsilon$ on $\mc X$ is obtained by the same
bounds and a classical compensated compactness argument.

The $\Gamma$-liminf inequality in Theorem~\ref{t:2.1} follows from the
variational definition \eqref{e:2.6} of $I_\varepsilon$. The
$\Gamma$-liminf inequality in Theorem~\ref{t:ecne} still follows from
\eqref{e:2.6} by choosing test functions of the form $\varepsilon
\vartheta(u^\varepsilon(t,x),t,x)$, with $\sigma \vartheta'' \le D$.

The $\Gamma$-limsup inequality in Theorem~\ref{t:2.1} is not difficult
if $\mu_{t,x}=\delta_{u(t,x)}$ for some smooth $u$; the general result
is obtained by taking the lower semicontinuous envelope. The
$\Gamma$-limsup statement in Theorem~\ref{t:ecne} is proven by
building, for each $u \in \mc S_\sigma$, a recovery sequence
$\{u^\varepsilon\}$ such that \emph{a priori} $H_\varepsilon (u^\varepsilon)
\to H(u)$. The convergence $u^\varepsilon \to u$ is then obtained by a
stability analysis of the parabolic equation \eqref{e:1.4} w.r.t.\ small
variations of the control $E$.

Eventually, in Appendix~\ref{s:B} we apply our results to
Hamilton-Jacobi equations.

\section{Representation of $I_\varepsilon$ and a priori bounds}
\label{s:3}

Given a bounded measurable function $a \ge 0$ on $[0,T]\times \bb R$
let $\mc D^1_a$ be the Hilbert space obtained by identifying and
completing the functions $\varphi\in C^\infty_{\mathrm{c}}([0,T]\times
\bb R)$ w.r.t.\ the seminorm $\langle \langle \varphi_x, a\, \varphi_x
\rangle \rangle^{1/2}$. Let $\mc D^{-1}_a$ be its dual space. The
corresponding norms are denoted respectively by $\|\cdot \|_{\mc
  D^1_a}$ and $\|\cdot \|_{\mc D^{-1}_a}$.

We first establish the connection between the cost functional
$I_\varepsilon$ and the perturbed parabolic problem \eqref{e:1.4}. The
following lemma is a standard tool in large deviations theory, see
e.g.\ \cite[Lemma~10.5.3]{KL}. We however detail its proof for sake of
completeness.

\begin{lemma2}
  \label{l:riesz}
  Fix $\varepsilon>0$ and let $u\in \mc U$. Then $I_\varepsilon(u) <
  +\infty$ iff there exists $\Psi^{\varepsilon,u} \in \mc
  D^1_{\sigma(u)}$ such that $u$ is a weak solution to \eqref{e:1.4}
  with $E = \Psi^{\varepsilon,u}_x$, namely for each $\varphi \in
  C^\infty_{\mathrm{c}}([0,T]\times\bb R)$
\begin{eqnarray}
\label{e:3.1}
\nonumber
&& 
\langle u(T), \varphi(T)\rangle  - \langle u(0), \varphi(0) \rangle 
\\
&&\qquad\quad
- \Big[ 
  \langle\langle u,\varphi_t\rangle \rangle
  + \big\langle\big\langle f(u) - \frac {\varepsilon}2  D(u)u_x 
  + \sigma(u) \Psi^{\varepsilon,u}_x,\varphi_x \big\rangle\big\rangle
 \Big] =0 \quad
\end{eqnarray}
In such a case $\Psi^{\varepsilon,u}$ is unique and
\begin{equation}
\label{e:3.2}
 I_\varepsilon(u) =  \frac 12 \Big\| u_t + f(u)_x 
         -\frac{\varepsilon}2 \big( D(u) u_x\big)_x
          \Big\|_{\mc D^{-1}_{\sigma(u)}}^2
   =\frac 12 \, \|\Psi^{\varepsilon,u}\|_{\mc D^{1}_{\sigma(u)}}^2
\end{equation}
\end{lemma2}

\proof
Fix $\varepsilon>0$ and $u \in \mc U$ such that
$I_\varepsilon(u)<+\infty$. The functional $\ell^\varepsilon_u$
defined in \eqref{e:2.5} can be extended to a linear functional on
$C^\infty_{\mathrm{c}}([0,T]\times \bb R)$ by setting
\begin{eqnarray}
\label{e:3.3}
\nonumber
\ell^u_\varepsilon(\varphi)& = &
        \langle u(T),\varphi(T)\rangle
        -\langle u(0),\varphi(0)\rangle 
     - \langle \langle  u,\varphi_t\rangle \rangle
        -\langle \langle f(u),\varphi_x\rangle \rangle
\\
&&
       + \: \frac{\varepsilon}2 \langle \langle D(u) u_x, 
                            \varphi_x\rangle \rangle
\end{eqnarray}
Since for any $\varphi \in C^\infty_{\mathrm{c}}([0,T]\times \bb R)$
the map $[0,T]\ni t \mapsto \langle u(t),\varphi(t)\rangle \in \bb R$
is continuous, it is easily seen that
\begin{equation*}
  I_\varepsilon(u) = \sup_{\varphi \in C^\infty_{\mathrm{c}}([0,T]\times \bb R)} \Big\{ 
  \ell^u_\varepsilon(\varphi) -  \frac 12  
      \langle \langle \sigma(u) \varphi_x, \varphi_x\rangle \rangle\Big\}
\end{equation*}
We claim that $\ell^u_\varepsilon$ defines a bounded linear functional on
$\mc D^1_{\sigma(u)}$. Indeed, since $I_\varepsilon(u)<+\infty$
\begin{equation*}
 \ell^u_\varepsilon(\varphi)
    \le I_\varepsilon(u)+  \frac 12 \langle \langle \sigma(u)\,\varphi_x,\varphi_x \rangle \rangle
 = I_\varepsilon(u)+ \frac 12 \|\varphi\|_{\mc D^1_{\sigma(u)}}^2
\end{equation*}
which shows that $\ell^u_\varepsilon(\varphi)=0$ whenever
$\langle \langle \sigma(u)\,\varphi_x,\varphi_x \rangle \rangle=0$ (as $\ell^u_\varepsilon(\cdot)$ is
$1$-homogeneous), namely $\ell^u_\varepsilon$ is compatible with the
identification in the definition of $D^1_{\sigma(u)}$ above. We also get that
$\ell^u_\varepsilon(\varphi)$ is bounded by the $\mc
D^1_{\sigma(u)}$-norm of $\varphi$ (up to a multiplicative constant), and it can therefore be extended by compatibility
and density to a continuous linear functional on $\mc D^1_{\sigma(u)}$. Still denoting by $\ell^u_\varepsilon$ such a
functional we get
\begin{equation}
\label{e:3.4}
I_\varepsilon(u)= \sup_{\varphi \in \mc D^1_{\sigma(u)}} \Big\{
              \ell^u_\varepsilon(\varphi)
             -\frac 12 \langle \langle \sigma(u)\varphi_x , 
                   \varphi_x\rangle \rangle\ \Big\}
\end{equation}
which is equivalent to the first equality in \eqref{e:3.2}. By Riesz
representation theorem we now get existence and uniqueness of
$\Psi^{\varepsilon,u} \in \mc D^1_{\sigma(u)}$ such that
$\ell^u_\varepsilon(\varphi)=\big(\Psi^{\varepsilon,u},\varphi\big)_{\mc
 D^1_{\sigma(u)}}$ for any $\varphi \in \mc D^1_{\sigma(u)}$, which
implies \eqref{e:3.1}. Riesz representation also yields
$I_\varepsilon(u)=\frac 12 \|\Psi^{\varepsilon,u}\|^2_{\mc
  D^1_{\sigma(u)}}$. The converse statements are obvious.
\qed

\smallskip
In the following lemma we give some regularity results for $u \in \mc
U$ with finite cost, and we prove some a priori bounds.

\begin{lemma2}
 \label{l:aprbound}
 Let $\varepsilon>0$ and $u\in \mc U$ be such that
 $I_\varepsilon(u)<+\infty$. Then $u \in
 C\big([0,T];L_{1,\mathrm{loc}}(\bb R) \big)$. Moreover for each
 entropy -- entropy flux pair $(\eta, q)$, each $\varphi \in
 C^\infty_{\mathrm{c}}([0,T] \times \bb R)$, and each $t \in [0,T]$
\begin{eqnarray}
\label{e:3.5}
\nonumber
&&
\langle \eta(u(t)),\varphi(t)\rangle 
    - \langle \eta(u(0)),\varphi(0)\rangle
    - \int_{[0,t]}\!ds\, \big[ \langle \eta(u), \varphi_s \rangle 
              +\langle q(u), \varphi_x \rangle \big]
\\
&&
\nonumber
\qquad
 =  -\frac{\varepsilon}2  \int_{[0,t]}\!ds\, \Big[
  \langle \eta''(u) D(u)u_x, \varphi\,u_x \rangle 
         +\langle
          \eta'(u)D(u) u_x,\varphi_x \rangle   
          \Big]
\\
&&   
\qquad \quad \,
+\int_{[0,t]}\!ds \Big[
        \langle \eta''(u) \sigma(u)\,u_x, 
          \Psi^{\varepsilon,u}_x\,\varphi \rangle 
       +\langle \eta'(u)\sigma(u)\, 
        \Psi^{\varepsilon,u}_x, \varphi_x \rangle \Big]
\qquad
\end{eqnarray}
where $\Psi^{\varepsilon,u}$ is as in Lemma~\ref{l:riesz}. Finally, there
exists a constant $C>0$ depending only on $f$, $D$ and $\sigma$ such
that for any $\varepsilon,\,L>0$
\begin{equation}
\label{e:3.6}
\varepsilon \int \! dt \int_{[-L,L]}\!dx \,u_x^2 \le
                   C \big[\varepsilon^{-1}I_\varepsilon(u)+L+1\big]
\end{equation}
\end{lemma2}

\proof
Recall that the linear functional $\ell^u_\varepsilon$ on $\mc
D^1_{\sigma(u)}$ is defined as the extension of \eqref{e:3.3}. Let
$\theta:=-f(u)+\frac{\varepsilon}{2}
D(u)\,u_x-\sigma(u)\,\Psi^{\varepsilon,u}_x \in
L_{2,\mathrm{loc}}\big([0,T]\times \bb R\big)$; by \eqref{e:3.1}
$u_t=\theta_x$ holds weakly. Since $I_\varepsilon(u)<+\infty$ we also
have $u_x \in L_{2,\mathrm{loc}} \big([0,T]\times \bb R \big)$, so
that $u \in C\big([0,T];L_{2,\mathrm{loc}}(\bb R) \big)$ by standard
interpolations arguments, see e.g.\ \cite{LM}. Since $u$ is bounded,
this is equivalent to the statement $u \in
C\big([0,T];L_{1,\mathrm{loc}}(\bb R) \big)$.

This fact implies that integrations by parts are allowed in the first
line on the r.h.s.\ of \eqref{e:3.3}, namely for each measurable
compactly supported $\phi:[0,T]\times \bb R \to \bb R$ with $\phi_x
\in L_2\big( [0,T]\times \bb R \big)$
\begin{equation}
\label{e:3.7}
\ell^u_\varepsilon(\phi)=
      \langle \langle u_t,\phi \rangle \rangle 
     + \langle \langle f(u)_x,\phi\rangle\rangle 
     + \frac{\varepsilon}2 \langle \langle D(u) u_x, \phi_x\rangle\rangle
\end{equation}
where indeed we understand $\langle \langle u_t,\phi \rangle \rangle \equiv - \langle \langle \theta,\phi_x \rangle \rangle$.
Since $u_x$ is locally square integrable, if $\eta \in C^{2}([0,1])$
and $\varphi \in C^\infty_{\mathrm{c}}([0,T]\times \bb R)$, then
$\eta'(u) \varphi$ has compact support and its weak $x$-derivative is
square integrable. We can thus evaluate \eqref{e:3.7} with $\phi$
replaced by $\eta'(u)\varphi$; since
$\ell^u_\varepsilon(\eta'(u)\varphi) =
\big(\Psi^{\varepsilon,u},\eta'(u)\varphi\big)_{\mc D^1_{\sigma(u)}}$
and $u \in C\big([0,T];L_{2,\mathrm{loc}}(\bb R) \big)$ we get
\eqref{e:3.5}.

To prove the last statement, consider an entropy -- entropy flux pair
$(\eta,q)$. By \eqref{e:3.4} and \eqref{e:3.7}
\begin{eqnarray*}
&&\!\!\!\!
I_\varepsilon(u) \ge \ell^u_\varepsilon (\varepsilon\,\eta'(u)\,\varphi) 
        -\frac{\varepsilon^2}{2} \big\langle \big\langle \big(\eta'(u)
         \varphi\big)_x , 
         \sigma(u) \big( \eta'(u)\varphi \big)_x \big\rangle \big\rangle
\\
&& \; 
 = \varepsilon \langle\eta(u(T)),\varphi(T)\rangle
      - \varepsilon \langle \eta(u(0)),\varphi(0)\rangle 
      -\varepsilon \big[\langle \langle \eta(u), \varphi_t \rangle \rangle 
      +\langle \langle q(u), \varphi_x \rangle\rangle \big] 
\\
&& \;
\phantom{=}
     + \frac{\varepsilon^2}2  \Big[ 
        \langle \langle D(u)\eta''(u)u_x^2, \varphi \rangle \rangle
     +  \langle \langle\eta'(u)D(u)u_x, \varphi_x \rangle \rangle
\\
&&  \;
\phantom{= + \frac{\varepsilon^2}2  \Big[}
-  \langle \langle \sigma(u) \eta''(u)^2 u_x^2,\varphi^2 \rangle \rangle
- \langle \langle \sigma(u) \eta'(u)^2\varphi_x, \varphi_x \rangle\rangle
\\
&& \;
\phantom{= + \frac{\varepsilon^2}2  \Big[}
- 2\langle \langle \sigma(u) \eta''(u) \eta'(u) u_x, 
\varphi\,\varphi_x \rangle \rangle
\Big]
\end{eqnarray*}
We now choose $\eta \ge 0$, uniformly convex and such that $\sigma
\eta'' \le D$, and for such a $\eta$ we let $\alpha:= \max_{v}\big[
D(v)\eta'(v)^2 / \eta''(v)\big]$, so that $\sigma\, (\eta')^2 \le
\alpha$. By Cauchy-Schwarz inequality
\begin{eqnarray*}
&&
2\big|   \langle \langle \sigma(u) \eta''(u) \eta'(u) u_x, 
                    \varphi\,\varphi_x \rangle \rangle \big|
\\ 
&&
\qquad\qquad
 \le  
    \langle \langle \sigma(u) \eta''(u)^2 u_x^2,
           \varphi^2\rangle \rangle
 +  \langle \langle \sigma(u) \eta'(u)^2, \varphi_x
           \varphi_x \rangle \rangle
\\
&&
\qquad\qquad
\le  
\langle \langle D(u) \eta''(u) u_x^2,
           \varphi^2\rangle \rangle
 + \alpha  \langle \langle \varphi_x,
           \varphi_x \rangle \rangle
\end{eqnarray*}
Letting $\zeta:[0,1]\to\bb R$ be such that $\zeta'= \eta' \,D$, and
integrating by parts we get $\langle\eta'(u)D(u)u_x, \varphi_x \rangle
= - \langle \zeta(u), \varphi_{x x} \rangle$. Collecting all the bounds
\begin{eqnarray*} 
&& 
\langle\eta(u(T)),\varphi(T)\rangle
      + \frac{\varepsilon}2 \langle \langle 
       D(u)\eta''(u)u_x^2, \varphi-2 \varphi^2 \rangle \rangle
\\
&& \qquad\qquad
 \le \varepsilon^{-1} I_\varepsilon(u)  
    + \langle \eta(u(0)),\varphi(0)\rangle
    +\langle \langle \eta(u), \varphi_t \rangle \rangle
    +\langle \langle q(u), \varphi_x \rangle \rangle
\\
&& \qquad\qquad \phantom{ \le}
    +\frac{\varepsilon}2 \langle \langle \zeta(u),\varphi_{x x}\rangle \rangle
    +\varepsilon\,\alpha\, \langle \langle \varphi_x, \varphi_x \rangle \rangle
\end{eqnarray*}
We now choose $\varphi$ independent of $t$ and such that $\varphi(x)
= 1/4$ for $|x|\le L$, $0\le \varphi(x) \le 1/4$ for $L\le |x| \le
L+1$, $\varphi(x) =0$ for $|x| \ge L+1$, and $\langle \varphi_x,
\varphi_x\rangle +\langle\varphi_{x x}, \varphi_{x x} \rangle \le 2$.
Since $q$, $\zeta$ are bounded and $\eta \ge 0$, estimate
\eqref{e:3.6} easily follows.
\qed

\begin{lemma2}
\label{l:equic}
The sequence of functionals $\{I_\varepsilon\}$ is equicoercive on $\mc
U$.
\end{lemma2}

\proof
Let $u \in \mc U$ be such that $I_\varepsilon(u)<+\infty$ and
$\Psi^{\varepsilon,u}$ be as in Lemma~\ref{l:riesz}. By \eqref{e:3.1},
\eqref{e:3.2} and the bound \eqref{e:3.6}, for each $s,t \in [0,T]$,
each $L>0$, each $\varphi \in C^\infty_{\mathrm{c}}(\bb R)$ supported
by $[-L,L]$
\begin{eqnarray*} 
&&
|\langle u(t)-u(s),\varphi\rangle|  = 
      \Big|\int_{[s,t]} \!dr\, \big\langle f(u) -\frac{\varepsilon}{2} D(u)u_x
      +\sigma(u)\Psi^{\varepsilon,u}_x,\varphi_x \big\rangle \Big| 
\\
&&\qquad  \le 
         \Big\{ 2 \int_{[s,t]\times [-L,L]}\!dr\, dx\, \Big[f(u)^2
                + \frac{\varepsilon^2}{4} D(u)^2 u_x^2 \Big] \Big\}^{1/2} 
           \big[|t-s| \langle \varphi_x,\varphi_x \rangle \big]^{1/2}
\\   
& &\qquad \quad     +\Big[ \int_{[s,t]}\!dr \langle \sigma(u) 
                                           \Psi^{\varepsilon,u}_x, 
               \Psi^{\varepsilon,u}_x \rangle\Big]^{1/2} 
    \big[|t-s| \langle \sigma(u)\varphi_x,\varphi_x \rangle
    \big]^{1/2}
\\
&&
\qquad    \le C\Big[1+L+I_\varepsilon(u)\Big]^{1/2} |t-s|^{1/2} \langle
          \varphi_x,\varphi_x \rangle^{1/2}
\end{eqnarray*}
for a suitable constant $C$ depending only on $f$, $D$, and $\sigma$.
Since $(U,d_U)$ is compact, see \eqref{e:2.1}, recalling \eqref{e:2.2}
and the Ascoli-Arzel\'{a} theorem, the equicoercivity of $\{I_\varepsilon\}$
on $\mc U$ follows.
\qed

\smallskip
As mentioned in the introduction, the assumption that $\sigma$ is
supported by $[0,1]$ allows us to consider only functions $u$ that
take values in $[0,1]$. More precisely, consider a cost functional
$\hat{I}_\varepsilon$ analogous to $I_\varepsilon$ but defined on
$L_{1,\mathrm{loc}}([0,T]\times \bb R)$. We next prove that, if $u \in
L_{1,\mathrm{loc}}([0,T]\times \bb R)$ is such that
$\hat{I}_\varepsilon(u)<+\infty$ and satisfies some growth conditions,
then $u$ takes values in $[0,1]$.

\begin{proposition2}
\label{t:boun}
Let $f,\,D,\,\sigma:\bb R \to \bb R$; assume $f$ Lipschitz, $\sigma$
and $D$ continuous and bounded, with $\sigma \ge 0$ and $D$ uniformly
positive.  Let ${\hat{I}}_\varepsilon: L_{1,\mathrm{loc}}([0,T] \times
\bb R) \to [0,+\infty]$ be defined as follows. If $f(u) \in
L_{2,\mathrm{loc}}([0,T]\times \bb R)$, we define
$\hat{I}_\varepsilon(u)$ as in \eqref{e:2.6}, and we set
$\hat{I}_\varepsilon(u)=+\infty$ otherwise. Suppose that $u \in
L_{1,\mathrm{loc}}([0,T] \times \bb R)$ is such that
$\hat{I}_\varepsilon(u)<+\infty$. Then $u \in
C\big([0,T];L_{1,\mathrm{loc}}(\bb R)\big)$. Moreover, if $\sigma$ is
supported by $[0,1]$, and $u$ is such that $u(0) \in U$ and $\int \!dt
\,dx\,|u(t,x)| e^{-r|x|}<+\infty$ for some $r>0$, then $u$ takes
values in $[0,1]$, hence $u\in \mc U$.
\end{proposition2}

\proof
Let $u \in L_{1,\mathrm{loc}}([0,T]\times \bb R)$ be such that
$\hat{I}_\varepsilon(u) <+\infty$. By the same arguments of
Lemma~\ref{l:aprbound}, since $f(u)\in L_{2,\mathrm{loc}}([0,T]\times
\bb R)$, $u_t=\theta_x$ holds weakly for some $\theta \in
L_{2,\mathrm{loc}}([0,T]\times \bb R)$. Hence, as in
Lemma~\ref{l:aprbound}, $u\in C\big([0,T];L_{1,\mathrm{loc}}(\bb
R)\big)$. Suppose now that $\sigma$ is supported by $[0,1]$. Pick a
sequence of strictly convex, strictly positive entropies $\eta_n\in
C^2(\bb R)$ such that: $\big|\eta'_n(u)\big|,\,\eta''_n(u) \le C_n$
for some $C_n>0$; for $u\in (0,1)$, $\eta_n(u)$ does not depend on $n$
and satisfies $0< c\le \eta''_n(u) \le D(u)/\sigma(u)$; $\eta_n$ is
decreasing for $u < 0$ and increasing for $u>1$; for $u \not\in [0,1]$
the sequence $\{\eta_n(u)\}$ increases pointwise to $+\infty$ as $n
\to \infty$.  Still following the proof of Lemma~\ref{l:aprbound}, for
$t \in [0,T]$ and $\varphi \in C^\infty_{\mathrm{c}}(\bb R)$
\begin{eqnarray*} 
&& \langle\eta_n(u(t)),\varphi \rangle
      + \frac{\varepsilon}2 \int_{[0,t]} \!ds\, \langle 
       D(u)\eta_n''(u)u_x^2, \varphi-2 \varphi^2 \rangle 
     \le \varepsilon^{-1} \hat{I}_\varepsilon(u)
\\ 
&&\qquad
+ \langle \eta_n(u(0)),\varphi\rangle 
+ \int_{[0,t]}\!ds \Big[\langle q_n(u), \varphi_x \rangle
+\frac{\varepsilon}2 \langle \zeta_n(u),\varphi_{x x}\rangle 
+\varepsilon\,\alpha\, \langle \varphi_x, \varphi_x \rangle
\Big]
\end{eqnarray*}
where $q_n$ and $\zeta_n$ are defined (up to a constant) by $q_n(v)=
\int^v\!dw\,\eta_n'(w)\,f'(w)$ and $\zeta_n'=\eta_n'\,D$, and
$\alpha:= \max_{u \in [0,1]} D(u)\eta_n'(u)^2 / \eta_n''(u)$ is a
constant independent of $n$, since $\sigma$ is supported by
$[0,1]$. Since $f$ is Lipschitz and $D$ is bounded, it is possible to
choose the arbitrary constants in the definition of $q_n$ and
$\zeta_n$ such that $|q_n|,\,|\zeta_n| \le C \eta_n$ for some constant
$C>0$ independent of $n$. In particular $\zeta_n,\,q_n \in
L_{1,\mathrm{loc}}([0,T]\times \bb R)$; for each $\varphi\in
C^\infty_{\mathrm{c}}(\bb R)$ such that $0\le \varphi(x) \le 1/2$
\begin{eqnarray*}
 && \langle \langle\eta_n(u),\varphi\rangle\rangle
   \le 
  T\,\varepsilon^{-1} \hat{I}_\varepsilon(u)
      +T\,\langle \eta_n(u(0)),\varphi\rangle
\\ 
&&\qquad\qquad
  +\int_{[0,T]}\!dt \int_{[0,t]}\!ds \,\Big[\langle q_n(u), \varphi_x \rangle
  +\frac{\varepsilon}2 \langle \zeta_n(u),\varphi_{x x}\rangle 
  +\varepsilon\,\alpha\, \langle \varphi_x, \varphi_x \rangle \Big]
\end{eqnarray*}
Let now $r$ be such that $\int \! dt \, dx\,e^{-r |x|} |u(t,x)|
<+\infty$. By a limiting procedure, the above bound holds for any
$\varphi \in C^\infty(\bb R)$ such that $0\le \varphi\le 1/2$ and $
\sup_{x\in \bb R} \, e^{ r|x|} \, \big[ |\varphi(x)|+|\varphi_x(x)| +
|\varphi_{x x}(x)| \big] <+\infty$. For such $\varphi$, by the choice
of $q_n,\,\zeta_n$
\begin{eqnarray*}
\frac 1T \langle \langle\eta_n(u),\varphi \rangle\rangle
   &   \le &  \varepsilon^{-1} \hat{I}_\varepsilon(u)
              +  \langle \eta_n(u(0)),\varphi\rangle     
\\ & &
        +\varepsilon \frac{\alpha T }{2} 
                 \langle \varphi_x, \varphi_x \rangle     
      +C \langle \langle \eta_n(u), 
      |\varphi_x| +\frac{\varepsilon}{2} |\varphi_{x x}|\rangle \rangle
\end{eqnarray*}
It is easy to verify that, given $L>0$ large enough, we can choose
$\varphi$ such that $\varphi(x)=1/2$ for $|x| \le L$,
$\varphi(x)=\frac 12 e^{-r |x-L|}$ for $|x|>2\,L$ and $ |\varphi_{x
  x}(x)|\le r |\varphi_x(x)| \le r^2 \, \varphi(x) \le r^2/2$ for
$|x|>L$. Moreover, with no loss of generality, we can assume that
$\frac 1T - C \,(r+\frac{\varepsilon}{2}r^2)>0$, otherwise we can
suppose $T$ small enough and iterate this proof. Therefore
\begin{eqnarray*}
&& \Big[\frac 1T - C \big(r+\frac{\varepsilon}{2}r^2\big)\Big] 
      \int_{[0,T] \times [-L,L]}\!dt\,dx\, \eta_n(u)
\\
&&\qquad\qquad
\le \varepsilon^{-1} \hat{I}_\varepsilon(u)
      +\langle \eta_n(u(0)),\varphi\rangle 
      +\varepsilon \frac{\alpha T}{2} \langle \varphi_x, \varphi_x \rangle     
\end{eqnarray*}
If $u(0)\in U$ the r.h.s.\ of this formula is finite and independent
of $n$, and therefore the l.h.s.\ is bounded uniformly in $n$. Taking
the limit $n \to \infty$, by the choice of $\eta_n$ necessarily
$u(t,x) \in [0,1]$ for a.e.\ $(t,x) \in[0,T] \times \bb R$.
\qed

\smallskip
The following result is not used in the sequel, but together with
Lem\-ma~\ref{l:riesz} and Proposition~\ref{t:boun}, motivates the choice
of $I_\varepsilon$ as the cost functional related to \eqref{e:1.2}.
\begin{proposition2}
  For each $\varepsilon>0$ the functional $I_\varepsilon: \mc U\to
  [0,+\infty]$ is lower semicontinuous.
\end{proposition2}

\proof
Let $\{u^n\}\subset \mc U$ be a sequence converging to $u$ in $\mc U$,
and such that $I_\varepsilon(u^n)$ is bounded uniformly in $n$. By
\eqref{e:3.6}, for each $L>0$ we have that $\int_{[0,T]\times [-L,L]}
\!dt\, dx\, (u^n_x)^2$ is also bounded uniformly in $n$. Therefore,
recalling definition \eqref{e:2.6}, the lower semicontinuity of
$I_\varepsilon$ is established once we show that $u^n$ converges to
$u$ strongly in $L_{1,\mathrm{loc}}([0,T]\times \bb R)$. Fix $L>0$ and
pick $\chi_L\in C^\infty_{\mathrm{c}} (\bb R)$ such that $0\le \chi_L
\le 1$ with $\chi_L(x)=1$ for $x \in [-L,L]$. We show that $u^{n,L}:=
u^n \, \chi_L$ converges to $u^{L}:= u \, \chi_L$ in $L_2([0,T]\times
\bb R)$. Choose a sequence of mollifiers $\jmath_k : \bb R \to \bb
R^+$ with $\int\!dx \,\jmath_k(x)=1$, then
\begin{eqnarray*}
& &\big\| u^{n,L} - u^{L} \big\|_{L_2([0,T]\times \bb R)} \le
      \big\| u^{n,L} - \jmath_k * u^{n,L} 
                 \big\|_{L_2([0,T]\times \bb R)}
\\ & &
\qquad \qquad 
      +\big\|\jmath_k * u^{n,L} -  \jmath_k * u^{L}
                 \big\|_{L_2([0,T]\times \bb R)}
      +\big\|\jmath_k * u^{L} - u^{L}
                 \big\|_{L_2([0,T]\times \bb R)}
\end{eqnarray*}
where the convolution is only in the space variable. For each $k$ the
second term on the r.h.s. above vanishes as $n\to \infty$ by the
convergence $u^n\to u$ in $\mc U$. Since the third term vanishes as
$k\to \infty$ it remains to show that the first one vanishes as $k\to
\infty$ uniformly in $n$. Integration by parts and Young inequality
for convolutions yield
\begin{eqnarray*}
&&\big\| u^{n,L} - \jmath_k * u^{n,L} \big\|_{L_2([0,T]\times\bb R)}  
\\
&&\qquad\qquad
\le  \Big\| \id_{[0,+\infty)} 
     - \int_{-\infty}^{\,\cdot} \!dy\,j_k(y)  \Big\|_{L_1(\bb R)}
   \: \big\|  u^{n,L}_x \big\|_{L_2([0,T]\times\bb R)} 
\end{eqnarray*}
The uniform boundedness of $\int_{[0,T]\times [-L,L]} \!dt \, dx\,
(u^n_x)^2$, \eqref{e:3.6} and the choice of $\chi_L$ imply that the
second term on the r.h.s.\ is bounded uniformly in $n$, while the
first term vanishes as $k \to \infty$.
\qed

\section{$\Gamma$-convergence of $\mc I_\varepsilon$}
\label{s:4}
In this section we prove the $\Gamma$-convergence of the parabolic
cost functional $\mc I_\varepsilon$ as $\varepsilon \to 0$, see
Theorem~\ref{t:2.1}. Some technical steps are postponed in
Appendix~\ref{s:A}.

\noindent
\textbf{Proof of Theorem~\ref{t:2.1}: equicoercivity of 
$\mc I_\varepsilon$.} 
Recall that $(\mc M,d_{\mc M})$ has been de\-fi\-ned in \eqref{e:2.3},
\eqref{e:2.4} and note that $(\mc N,d_{\mathrm{w}})$ is compact. By
Lemma~\ref{l:equic}, for each $C>0$ there exists a compact $K_C
\subset \mc U$, such that for any $\varepsilon$ small enough $\{\mu
\in \mc M\,:\:\mc I_\varepsilon(\mu)\le C\}\subset \{\mu \in \mc
M\,:\:\mu_{t,x}=\delta_{u(t,x)} \text{for some $u \in K_C$}\} =:\mc
K_C$. In order to prove that $\mc K_C$ is compact in $(\mc
M,d_{\mc M})$, consider a sequence $\{\mu^n=\delta_{u^n}\} \subset \mc
K_C$.  Then there exists a subsequence $\{\mu^{n_j}\}$ such that,
for some $\mu \in \mc N$ and $u\in \mc U$, $\mu^{n_j} \to \mu$ in
$(\mc N,d_{\mathrm{w}})$, and $\mu^{n_j}(\imath)=u^{n_j}\to u$ in
$\mc U$, hence $\mu(\imath)=u$. Therefore $\mu \in \mc M$ and
$\mu^{n_j}\to \mu$ in $(\mc M,d_{\mc M})$.
\qed

\smallskip
\noindent
\textbf{Proof of Theorem~\ref{t:2.1}: $\Gamma$-liminf inequality.}
Let $\{\mu^\varepsilon\} \subset \mc M$ be a sequence converging to
$\mu$ in $\mc M$. In order to prove $\varliminf_{\varepsilon \to 0}
\mc I_\varepsilon(\mu^\varepsilon) \ge \mc I(\mu)$, it is not
restrictive to assume $\mc I_\varepsilon(\mu^\varepsilon)<+\infty$,
and therefore $\mu^\varepsilon_{t,x}= \delta_{u^\varepsilon(t,x)}$ for
some $u^\varepsilon \in \mc U$. For each $\varphi \in
C^\infty_{\mathrm{c}}((0,T)\times \bb R)$, recalling definition
\eqref{e:2.6}
\begin{eqnarray*}
  &&\!\!
  \mc I_\varepsilon(\mu^\varepsilon) \ge 
  \ell^{u^\varepsilon}_\varepsilon(\varphi) 
  -\frac 12 \|\varphi\|_{\mc D^1_{\sigma(u^\varepsilon)}}^2
  \\
  &&\; = 
  - \langle\langle \mu^\varepsilon(\imath),\varphi_t \rangle \rangle
  -\langle\langle\mu^\varepsilon(f), \varphi_x\rangle \rangle
  -\frac 12\, \langle\langle\mu^\varepsilon(\sigma)\varphi_x, 
  \varphi_x\rangle \rangle
  + \frac{\varepsilon}2 \langle \langle D(u^\varepsilon) 
  u^\varepsilon_x,\varphi_x \rangle \rangle
\end{eqnarray*}
Let $d \in C^1([0,1])$ be such that $d'(u)=D(u)$. Then $D(u^\varepsilon)
u^\varepsilon_x=d(u^\varepsilon)_x$, and an integration by parts shows that
the last term on the r.h.s.\ of the previous formula vanishes as
$\varepsilon \to 0$. Hence
\begin{equation*}
\varliminf_{\varepsilon \to 0} \mc I_\varepsilon(\mu^\varepsilon) \ge  
- \langle\langle \mu(\imath),\varphi_t \rangle \rangle
   -\langle\langle \mu(f), \varphi_x\rangle \rangle
   -\frac 12\, \langle\langle\mu(\sigma)\varphi_x, 
                             \varphi_x\rangle \rangle
\end{equation*}
By optimizing over $\varphi \in C^\infty_{\mathrm{c}}((0,T)\times \bb
R)$ the $\Gamma$-liminf inequality follows.
\qed

\smallskip
\noindent
\textbf{Proof of Theorem~\ref{t:2.1}: $\Gamma$-limsup inequality.}
  Let
\begin{eqnarray}
\label{e:mg}
&&\!\!\!\!\!\!\!\!\!\!
\begin{array}{lcl}
\!\!
\mc M_g & \!:= &\!
\Big\{\mu \in \mc M\,:\: \mc I(\mu)<+\infty,\,
\exists r,\,L>0,\, \exists \mu_\infty \in \mc P([0,1]) 
         \; \text{such that}
\\ & & 
\phantom{\Big\{\mu \in \mc M\,:\:\; }
\mu(\imath),\, \mu(\sigma) \ge r,\;
\mu_{t,x}=\mu_\infty\,\text{for $|x|>L$} \Big\}
\end{array}
\\
&&\!\!\!\!\!\!\!\!\!\!
\label{e:m0}
\mc M_0 :=
\Big\{\mu \in \mc M_g\,:\: \mu =\delta_u \text{ for some 
         $u \in  C^1\big([0,T]\times \bb R;[0,1]\big)$} 
\Big\}
\end{eqnarray}
and define $\tilde{\mc I}:\mc M \to [0,+\infty]$ by
\begin{equation}
\label{e:itilde}
\widetilde{\mc I}(\mu):=
\begin{cases}
\mc I(\mu)  &  \text{if $\mu \in \mc M_0$}
 \\
+\infty     &  \text{otherwise}
\end{cases}
\end{equation}
We claim that for $\mu \in \mc M_0$, a recovery sequence is simply
given by $\mu^\varepsilon=\mu$. Indeed, if $\mu=\delta_{u}$ for some
$u\in C^1\big([0,T]\times \bb R;[0,1]\big)$, we have
\begin{equation*}
\begin{array}{lcl}
\mc I_\varepsilon(\mu^\varepsilon) & = & 
 I_\varepsilon(u) =
       \frac 12 \Big\| u_t + f(u)_x - \frac{\varepsilon}2 
             \big( D(u) u_x\big)_x  \Big\|_{\mc D^{-1}_{\sigma(u)}}^2
\\
& \le&
 \frac{1 + \varepsilon}{2}
            \Big\| u_t + f(u)_x \Big\|_{\mc D^{-1}_{\sigma(u)}}^2
    + \frac{1+\varepsilon^{-1}}{2}
      \Big\|\frac{\varepsilon}{2} \big( D(u) u_x\big)_x 
        \Big\|_{\mc D^{-1}_{\sigma(u)}}^2
\end{array}
\end{equation*}
As $\mu \in \mc M_g$, $u$ is constant for $|x|$ large enough, in
particular $u_x \in L_2([0,T]\times \bb R)$. Since we have also
$\sigma(u)\ge r>0$, the last term in the above formula vanishes as
$\varepsilon \to 0$. Hence $\Gamma$\textrm{-}$\varlimsup_\varepsilon \mc
I_\varepsilon \le \tilde{\mc I}$. As well known, see e.g.
\cite[Prop.~1.28]{Br}, any $\Gamma$-limsup is lower semicontinuous;
the proof is then completed by Theorem~\ref{t:lsce} below.
\qed

\smallskip
The relaxation of the functional $\widetilde{\mc I}$ on $\mc M$
defined in \eqref{e:itilde} might have an independent interest; in the
following result we show it coincides with $\mc I$, as defined in
\eqref{e:2.8}.
\begin{theorem}
\label{t:lsce}
$\mc I$ is the lower semicontinuous envelope of $\widetilde{\mc I}$.
\end{theorem}

The following representation of $\mc I$ is proven similarly to
Lemma~\ref{l:riesz}.
\begin{lemma2}
\label{l:riesz2}
Let $\mu \in \mc M$. Then $\mc I(\mu) < +\infty$ iff
there exists $\Psi^\mu \in \mc D^1_{\mu(\sigma)}$ such that $\mu$ is a
measure-valued solution to $u_t+f(u)_x = -
\big(\sigma(u)\Psi^{\mu}_x\big)_x$, namely
\begin{equation}
\label{e:riesz2}
\mu(\imath)_t+\mu(f)_x =-\big(\mu(\sigma)\Psi^\mu_x\big)_x
\end{equation}
holds weakly. In such a case $\Psi^\mu$ is unique and
\begin{equation*}
\mc I(\mu)
  = \frac 12 \Big\|\mu(\imath)_t+\mu(f)_x\Big\|_{\mc D^{-1}_{\mu(\sigma)}}^2
  = \frac 12 \|\Psi^\mu\|_{\mc D^1_{\mu(\sigma)}}^2
\end{equation*}
Furthermore, suppose that $\mu(\sigma)\ge r$ for some constant $r>0$.
Then $\mc I(\mu)<+\infty$ iff there exists $G^\mu \in L_2([0,T]\times
\bb R)$ such that weakly 
\begin{equation}
\label{e:measg}
\mu(\imath)_t+\mu(f)_x=-G^\mu_x
\end{equation}
In such a case $\Psi^{\mu}_x$ can be identified with a
function in $L_2([0,T]\times \bb R)$, and
\begin{equation}
\label{e:Gdef}
G^\mu=\mu(\sigma)\Psi^\mu_x,
 \qquad
\mc I(\mu)=\frac 12
\int \!dt\,dx \frac{\big(G^\mu(t,x)\big)^2}{\mu_{t,x}(\sigma)}
\end{equation}
\end{lemma2}

The following remark is a consequence of Lemma~\ref{l:riesz2}.
\begin{remark2}
\label{l:iconv}
Let $\{\mu^k\}\subset \mc M$ be such that $\mu^k \to \mu$ in $\mc M$,
$\mc I(\mu^k)<+\infty$ and $\mu^k(\sigma)\ge r$ for some $r>0$.
Let also $G^{\mu^k}$ be defined as in Lemma~\ref{l:riesz2}. If
$\mu^k(\sigma)\to \mu(\sigma)$ strongly in
$L_{1,\mathrm{loc}}\big([0,T]\times \bb R\big)$ and $\{G^{\mu^k}\}$ is
strongly precompact in $L_2\big([0,T]\times \bb R\big)$, then $\mc
I(\mu^k) \to \mc I(\mu)$.
\end{remark2}

Throughout the proof of Theorem~\ref{t:lsce}, approximation of Young
measures by piecewise smooth measures is a much used procedure. In
particular we will refer repeatedly to the following result, which is
a simple restatement of the Rankine-Hugoniot condition for the
divergence-free vector field $(\mu(\imath),\, \mu(f)+G^\mu)$ on
$(0,T)\times \bb R$.

\begin{lemma2}
\label{l:ruho}
Let $\gamma:(0,T)\to \bb R$ be a Lipschitz map with a.e.\ derivative
$\dot{\gamma}$, and let $O^\mp \subset (0,T)\times \bb R$ be a left,
resp.\ a right, open neighborhood of the graph of $\gamma$; namely
$\mathrm{Graph}(\gamma) \subset \mathrm{Closure}(O^-)\cap
\mathrm{Closure}(O^+)$, and for all $(t,x) \in O^-$, resp.\ $(t,x)\in
O^+$, the inequality $x<\gamma(t)$, resp.\ $x>\gamma(t)$, holds. Let
also $O:=O^+ \cup O^- \cup \mathrm{Graph}(\gamma)$. Suppose that a
Young measure $\mu \in \mc M$ is such that, for each continuous
function $F\in C([0,1])$ the map $(t,x)\mapsto \mu_{t,x}(F)$ is
continuously differentiable in $O^- \cup O^+$, and such that there
exist the respective traces $\mu^\mp(F)$ of $\mu(F)$ on the graph of
$\gamma$. Then there exists a map $G:O\to \bb R$, defined up to an
additive measurable function of the $t$ variable, which is continuous
in $O^- \cup O^+$, admits traces $G^\mp$ on the graph of $\gamma$, and
is such that \eqref{e:measg} holds weakly in $O$.  Moreover the
Rankine-Hugoniot condition holds for a.e.\ $t\in [0,T]$, namely
\begin{equation}
\label{e:ruho}
G^+-G^-= \big[\mu(\imath)^+ - \mu(\imath)^-\big]\dot{\gamma} 
       - \big[\mu(f)^+ - \mu(f)^-\big]
\end{equation}
\end{lemma2}

\noindent\textbf{Proof of Theorem~\ref{t:lsce}.}
Since $\mc I$ is lower semicontinuous, it is enough to prove that $\mc
M_0$, as defined in \eqref{e:m0}, is $\mc I$-dense in $\mc M$, namely
that for each $\mu \in \mc M$ with $\mc I(\mu)<+\infty$, there exists
a sequence $\{\mu^k\} \subset \mc M_0$ such that $\mu^k \to \mu$ in
$\mc M$ and $\varlimsup_k \mc I(\mu^k) \le \mc I(\mu)$ (we will also
say that $\mu^k$ $\mc I$-converges to $\mu$). We split the proof in
several steps.

\smallskip
\noindent\textit{Step 1.} Here we show that $\mc M_0$ is $\mc I$-dense in the
set of Young measures which are a finite convex combination of Dirac
masses for a.e. $(t,x)$. More precisely, recalling definition
\eqref{e:mg}, we set
\begin{eqnarray*}
\begin{array}{lcl}
& &
\mc M_1^n := \Big\{ 
                \mu \in \mc M_g\, :\: 
                   \mu=\sum_{i=1}^n \alpha^i\delta_{u^i}
     \text{ for some } \alpha^i \in L_\infty\big([0,T]
                              \times \bb R;[0,1]\big)
\\
& & \phantom{\mc M_1^n := \Big\{ } \;
     \text{with } \sum_{i=1}^n \alpha^i=1
     \text{ and } u^i \in 
            L_\infty\big([0,T]\times \bb R; [0,1]\big)  
     \Big\}
\end{array}
\end{eqnarray*}
and
\begin{equation*}
\mc M_1 := 
        \bigcup_{n=1}^\infty \mc M_1^n
\end{equation*}
In this step, we prove that $\mc M_0$ is $\mc I$-dense in $\mc M_1$.
We proceed by induction on $n$; to this aim, for $n\ge 1$, we
introduce the auxiliary sets
\begin{eqnarray*}
&&
\begin{array}{lcl}
&&
\upbar{\mc M}_1^n  := \Big\{
     \mu \in \mc M_g\, :\:    
                        \exists r>0  \text{ such that }  
         \mu = \sum_{i=1}^n \alpha^i \delta_{u^i},
           
\\ & & \phantom{\upbar{\mc M}_1^n  := \Big\{ }
       \text{ for some }  \alpha^i \in L_\infty\big([0,T]\times \bb R;[r,1]\big)
    \text{ with }  \sum_{i=1}^n \alpha^i=1

\\ & & \phantom{\upbar{\mc M}_1^n  := \Big\{ }
    \text{ and } u^i \in C^0\big([0,T]\times \bb R; [0,1]\big)
 \Big\}
\end{array}
\\
&&
\begin{array}{lcl}
  & &\widetilde{\mc M}_1^n := \Big\{
  \mu \in \mc M_g\, :\:
  \exists r>0 \text{ such that }  
  \mu=\sum_{i=1}^n \alpha^i \delta_{u^i}
  \\
  & &
  \phantom{\widetilde{\mc M}_1^n := \Big\{ }
 \text{ for some }
  \alpha^i \in C^1\big([0,T]\times \bb R;[r,1]\big)
  \text{ with }  \sum_{i=1}^n \alpha^i=1
  \\
  & &
  \phantom{\widetilde{\mc M}_1^n := \Big\{ }
  \text{ and } u^i \in C^1\big([0,T]\times \bb R; [r,1-r]\big)
  \text{ with }
  u^{i+1}\ge u^i +r
  \Big\}
\end{array}
\end{eqnarray*}
Note that $\widetilde{\mc M}_1^n \subset \upbar{\mc M}_1^n \subset \mc
M_1^n$, and $\widetilde{\mc M}_1^1 \subset \mc M_0$. We claim that for
each $n \ge 1$, $\widetilde{\mc M}_1^n$ is $\mc I$-dense in
$\upbar{\mc M}_1^n$, that $\upbar{\mc M}_1^n$ is $\mc I$-dense in $\mc
M_1^n$, and that $\mc M_1^n$ is $\mc I$-dense in $\widetilde{\mc
  M}_1^{n+1}$. The $\mc I$-density of $\mc M_0$ in $\mc M_1$ then
follows by induction. The previous claims are proven in
Appendix~\ref{s:A}.

\smallskip
\noindent\textit{Step 2.} 
In this step we prove that $\mc M_1$ is $\mc I$-dense in $\mc M_g$,
see \eqref{e:mg}. We use the following elementary extension of the
mean value theorem.

\begin{lemma2}
\label{l:subsimp}
Let $X$ be a connected compact separable metric space, $F_1$,$\ldots$,
$F_d \in C(X)$ be continuous functions on $X$, and $\bb P \in \mc
P(X)$ be a Borel probability measure on $X$.  Then there exist
$\alpha^1,\ldots,\alpha^d \ge 0$ with $\sum_i \alpha^i=1$,
$x^1,\ldots,\,x^d \in X$ such that $\bb P(F^i)=\sum_{j=1}^d \alpha^j
F^i(x^j)$, $i=1,\ldots,\,d$.  Furthermore there exists a sequence
$\{\bb P^n\} \subset \mc P(X)$ converging weakly* to $\bb P$, such
that each $\bb P^n$ is a finite convex combination of Dirac masses,
$\bb P^n(F^i)=\bb P(F^i)$ for $i=1,\ldots,d$, and for each $n$ the map
$\mc P(X) \ni \bb P \mapsto \bb P^n \in \mc P(X)$ is Borel measurable
w.r.t.\ the weak* topology.
\end{lemma2}

\proof 
It is easy to see that the point $\bb P(F):= \big(\bb
P(F_1),\ldots,\,\bb P(F_d)\big)\in \bb R^d$ belongs to the closed
convex hull of the set $B:=\{\big(F_1(x),\ldots,\,F_d(x)\big),\,x\in
X\} \subset \bb R^d$. Since $B$ is compact and connected, Caratheodory
theorem implies that $\bb P(F)$ is a convex combination of at most $d$
points in $B$, namely the first statement of the lemma holds. Since
$X$ is compact, for each integer $n \ge 1$, there exist an integer
$k=k(n)$ and pairwise disjoint measurable sets $A_1^n,\ldots,\,A_k^n
\subset X$, such that $\bb P(X \setminus \cup_{l=1}^k A_l^n)=0$, $\bb
P(A_l^n)>0$, and $\mathrm{diameter}(A_l^n) \le n^{-1}$,
$l=1,\ldots,\,k$. For $l=1,\ldots,k$, let $\bb P(\cdot|A_l^n) \in \mc
P(X)$ be defined by $\bb P(B|A_l^n):=\bb P(A_l^n\cap B)/\bb P(A_l^n)$
for any Borel set $B \subset X$.  By the first part of the lemma,
there exists a probability measure $\bb P^n_l \in \mc P(X)$, which is
a convex combination of $d$ Dirac masses, such that $\bb
P^n_l(F_i)=\bb P(F_i|A_l^n)$. The sequence $\{\bb P^n\}$ defined as
$\bb P^n(\cdot):=\sum_{l=1}^k \bb P(A_l^n) \bb P^n_l(\cdot)$ satisfies
the requirements of the lemma.
\qed

\smallskip
Let $\mu \in \mc M_g$. By Lemma~\ref{l:subsimp}, there exists a
sequence $\{\mu^n\}\subset \mc M$ converging to $\mu$ in $\mc M$ such
that $\mu_{t,x}$ is a convex combination of Dirac masses $(t,x)$ for
a.e.\ $(t,x)$, and
$\mu^n(\imath)=\mu(\imath)$,
$\mu^n(f)=\mu(f)$, $\mu^n(\sigma)=\mu(\sigma)$. Hence $\mc I(\mu^n) =
\mc I(\mu)$ and $\mu^n \in \mc M_1$.

\smallskip
\noindent \textit{Step 3.} Recall Lemma~\ref{l:riesz2} and set
\begin{equation*}
\begin{array}{lcl}
& & \mc M_3 := \Big\{
\mu \in \mc M\,:\: \mc
I(\mu)<+\infty, \,\exists r>0 \text{ such that } 
\mu(\imath),\,\mu(\sigma) \ge r,
\\ & &
\phantom{ \mc M_3 := \Big\{ }
G^\mu \in C^1\big([0,T]\times \bb R\big) \cap
L_{\infty}\big([0,T]\times \bb R\big),
\\ & &
\phantom{ \mc M_3 := \Big\{ }
\text{ for each $F \in  C([0,1])$ }\;
 \mu(F) \in C^1\big([0,T]\times \bb R\big)
\Big\} 
\end{array}
\end{equation*}
In this step we prove that $\mc M_g$ is $\mc I$-dense in $\mc M_3$.

Let $\mu \in \mc M_3$, and choose a constant $u_\infty >0$ such that
$\mu(\imath)-u_\infty > \delta$ for some $\delta >0$. Define the maps
$\gamma^k_{\pm} \in C([0,T])\cap C^1((0,T))$ as the solutions to the
Cauchy problems
\begin{equation*}
\begin{cases}
{\displaystyle
\dot{\gamma}(t)=\frac{G^{\mu}(t,\gamma(t))+\mu_{t,\gamma(t)}(f) 
                  - f(u_\infty)}{\mu_{t,\gamma(t)}(\imath) - u_\infty} 
}
\\ %
\gamma(0)=\pm k
\end{cases}
\end{equation*}
$\gamma^k_\pm$ are well-defined by the smoothness hypotheses on $\mu$
and $G^\mu$. On the other hand, since we assumed $G^\mu$ to be
uniformly bounded, $|\gamma^k_{\pm}(t) \mp k| \le C$, for some
constant $C>0$ not depending on $k$. We define, for $k > C$, $\mu^k$
by $\mu^k_{t,x}=\mu_{t,x}$ if $\gamma^k_-(t) <x<\gamma^k_+(t) $ and
$\mu^k_{t,x}=\delta_{u_\infty}$ otherwise. Clearly $\mu^k \to \mu$ in
$\mc M$ as $k \to \infty$. We also let $G^{\mu^k}(t,x)=G^{\mu}(t,x)$
if $\gamma^k_-(t) <x<\gamma^k_+(t)$, and $G^{\mu^k}(t,x)=0$ otherwise.
By \eqref{e:ruho} and the definition of $\gamma^k_{\pm}$, the equation
$\mu^k(\imath)_t +\mu^k(f)_x=-G^{\mu^k}_x$ holds weakly in
$(0,T)\times \bb R$. In particular, by Lemma~\ref{l:riesz2}, $\mc
I(\mu^k) \le \mc I(\mu)$.

\smallskip
\noindent \textit{Step 4.} Here we prove that $\mc M_3$ is $\mc
I$-dense in
\begin{equation*}
\mc M_4:= \{\mu \in \mc M\,:\:I(\mu)<+\infty,\,\exists r>0 
 \text{ such that }\mu(\imath),\,\mu(\sigma)
 \ge r\}
\end{equation*}
Let $\mu \in \mc M_4$ and $\{\jmath^k\}_{k \ge 1}\subset
C^\infty_{\mathrm{c}}(\bb R \times \bb R)$ be a sequence of smooth
mollifiers supported by $[-T/k,T/k]\times [-1,1]$. For $k \ge 1$, let
us define the rescaled time-space variables $b^k:[0,T]\times \bb R \to
\bb R \times \bb R$ by
\begin{equation}
\label{e:b}
b^k(t,x):=\Big(\frac{t+T/k}{1+2T/k} ,\frac{x}{1+2T/k}\Big)
\end{equation}
For $k\ge 1$ we also define the Young measure $\mu^k$ by setting for
$F \in C([0,1])$ and $(t,x) \in [0,T]\times \bb R$
\begin{equation*}
\mu^{k}_{t,x}(F)  :=  
      \int\!dy\,ds\,\jmath^k(t-s,x-y) \mu_{b^k(s,y)}(F)
\end{equation*}
It is immediate to see that $\mu^k \in \mc M_3$. Moreover, as $k \to
\infty$, $\mu^k \to \mu$ in $\mc M$ and $\mu^k(F) \to \mu(F)$ strongly
in $L_{1,\mathrm{loc}}([0,T]\times \bb R)$ for each $F \in C([0,1])$.

Let us also define $G^{\mu^k} \in L_2([0,T]\times \bb R)$ by
\begin{equation*}
G^{\mu^{k}}_{t,x}  :=  
      \int \! dy \,ds\, \jmath^k(t-s,x-y) 
                 G^{\mu}\big(b^k(s,y)\big)
\end{equation*}
Then  $\mu^k(\imath)_t
+\mu^k(f)_x=-G^{\mu^k}_x$ holds weakly, and $G^{\mu^k} \to G^{\mu}$
in $L_2([0,T]\times \bb R)$ as $k \to \infty$. The proof is then
achieved by  Remark~\ref{l:iconv}.

\smallskip
\noindent \textit{Step 5.} $\mc M_4$ is $\mc I$-dense in $\mc M$. For
$\mu \in \mc M$ with $\mc I(\mu)<+\infty$, we define
$\mu^k:=(1-k^{-1})\mu+k^{-1}\delta_{1/2}$.  Clearly $\mu^k \to \mu$ in
$\mc M$, and $\mu^k(\imath) \ge k^{-1}/2$, $\mu^k(\sigma) \ge k^{-1}
\sigma(1/2)$. Therefore $\mu^k \in \mc M_4$.  From \eqref{e:2.8} it
follows that $\mc I$ is convex, and since $\mc I(\delta_{1/2})=0$, we
have $\mc I(\mu^k)\le (1-k^{-1})\mc I(\mu)$.
\qed

\smallskip
The following proposition is easily proven, and will be used in the
proof of Corollary~\ref{c:rely}.

\begin{proposition2}
\label{p:contraction}
Let $X$, $Y$ be complete separable metrizable spaces, and let
$\omega:X \to Y$ be continuous. Let also $\{\mc F_\varepsilon\}$ be a
family of functionals $\mc F_\varepsilon:X \to [-\infty,+\infty]$. Let us
define $F_\varepsilon:Y \to [-\infty,+\infty]$ by
\begin{equation*}
  F_\varepsilon(y) = \inf_{x \in \omega^{-1}(y)} \mc F_\varepsilon(x)
\end{equation*}
Then
\begin{equation*}
  \big(\Glimsup_{\varepsilon\to 0}  F_\varepsilon\big)(y) \le 
  \inf_{ x \in \omega^{-1}(y)} 
       \big(\Glimsup_{\varepsilon\to 0}  \mc F_\varepsilon\big)(x)
\end{equation*}
Furthermore if $\{\mc F_\varepsilon\}$ is equicoercive on $X$ then
$\{F_\varepsilon\}$ is equicoercive on $Y$. In such a case
\begin{equation*}
  \big(\Gliminf_{\varepsilon\to 0}  F_\varepsilon\big)(y) \ge 
  \inf_{ x \in \omega^{-1}(y)} 
       \big(\Gliminf_{\varepsilon\to 0}  \mc F_\varepsilon\big)(x)
\end{equation*}
\end{proposition2}

\noindent
\textbf{Proof of Corollary~\ref{c:rely}.}
Since the map $\mc M \ni \mu \mapsto \mu(\imath) \in \mc U$ is
continuous, by Proposition~\ref{p:contraction} we have that
$I_\varepsilon$ is equicoercive on $\mc U$ (which we already knew from
Lemma~\ref{l:equic}) and $\Gamma$-converges to $I:\mc U \to
[0,+\infty]$ defined by
\begin{equation*}
I(u)= \inf_{\mu\in \mc M\,:\: \mu(\imath)=u} \mc I(\mu)
\end{equation*}
Recall that, if $\mc I(\mu)<+\infty$, $\Psi^\mu_x$ has been defined in
Lemma~\ref{l:riesz2}. Equality \eqref{e:riesz2} yields
\begin{eqnarray*}
&& 
I(u)=  \inf \Big\{\langle \langle \mu(\sigma) \Psi^\mu_x , 
                               \Psi^\mu_x \rangle \rangle,\, 
          \Phi \in L_{2,\mathrm{loc}}\big([0,T]\times \bb R\big),\,
  \mu \in \mc M\,:
\\ 
&& 
\qquad\quad
  \mc I(\mu)<+\infty,\,
 \mu(\imath)=u,\,
  \Phi_x=-\mu(\imath)_t \text{ weakly},\,
  \mu(\sigma) \Psi^{\mu}_x = \Phi-\mu(f)
  \Big\}
\end{eqnarray*}
The corollary then follows by direct computations.
\qed

\section{$\Gamma$-convergence of $H_\varepsilon$}
\label{s:5}

\noindent\textbf{Proof of Proposition~\ref{p:kin}.} 

\smallskip \noindent 
(i) $\Rightarrow$ (ii). 
We first show that $\|\wp_{\eta,u}\|_{\mathrm{TV},L}$ is finite for
each $\eta$ such that $0 \le \eta'' \le c$. It is easily seen that for
each $\varphi \in C^\infty_{\mathrm{c}}\big((0,T)\times
(-L,L);[0,1]\big)$ there exists $\bar{\varphi}\in
C^\infty_{\mathrm{c}}\big((0,T)\times (-L,L);[0,1]\big)$ such that
$\bar{\varphi} \ge \varphi$ and
$\||\bar{\varphi}_t|+|\bar{\varphi}_x|\|_{L_1} \le 2(2L+T)$.
Therefore
\begin{eqnarray*}
 \wp_{\eta,u}(-\varphi) & = &
  \wp_{\eta,u}(\bar{\varphi}-\varphi)-\wp_{\eta,u}(\bar{\varphi})
\\
  & \le &
  \|\wp_{\eta,u}^+\|_{\mathrm{TV},L} 
  + \langle \langle \eta(u),\bar{\varphi}_t \rangle \rangle
  + \langle \langle q(u),\bar{\varphi}_x \rangle \rangle 
\\ & \le  &
  \|\wp_{\eta,u}^+\|_{\mathrm{TV},L} 
  + 2 (\|\eta\|_{\infty}+\|q\|_{\infty})(2L+T)
\end{eqnarray*}
and thus $\|\wp_{\eta,u}\|_{{\mathrm{TV},L}} \le 2
\|\wp_{\eta,u}^+\|_{\mathrm{TV},L} + 2
(\|\eta\|_{\infty}+\|q\|_{\infty})(2L+T)$.

Let now $\tilde{\eta}(v):=c\,v^2/2$, and for $\eta \in C^2([0,1])$
arbitrary, let $\alpha:=c^{-1} \max_v |\eta''(v)|$. Then
$\wp_{\eta,u}= -\alpha
\wp_{\tilde{\eta}-\eta/\alpha,u}+\alpha\wp_{\tilde{\eta},u}$. Since both
$\tilde{\eta}-\eta/\alpha$ and $\tilde{\eta}$ are convex with second
derivative bounded by $c$, $\wp_{\eta,u}$ is a linear combination of
Radon measures, and thus a Radon measure itself.

\smallskip \noindent 
(ii) $\Rightarrow$ (iii). Throughout this proof, we say that
$\eta_1,\,\eta_2 \in C^2([0,1])$ are equivalent, and we write $\eta_1
\sim \eta_2$, iff $\eta_1''=\eta_2''$. We identify $C^2([0,1])/\sim$
with $C([0,1])$, which we equip with the topology of uniform
convergence. For $u \in \mc X$ a weak solution to \eqref{e:1.1}, for
$\varphi \in C^\infty_{\mathrm{c}}\big((0,T)\times \bb R\big)$, the
linear mapping $C^2([0,1]) \ni \eta \mapsto \wp_{\eta,u}(\varphi) \in
\bb R$ is compatible with $\sim$, and it thus defines a linear mapping
$P_{\varphi,u}:C([0,1]) \to \bb R$. It is immediate to see that
$P_{\varphi,u}$ is continuous, and by (ii) for each $\eta \in
C^2([0,1])$ and $L>0$
\begin{equation*}
\sup \big\{ P_{\varphi,u}(\eta''),
   \,\varphi \in C^\infty_{\mathrm{c}}\big((0,T)\times
               (-L,L)\big),\, |\varphi|\le 1 
     \big\} 
=\|\wp_{\eta,u}\|_{TV,L} <+\infty
\end{equation*}
By Banach-Steinhaus theorem
\begin{eqnarray*}
&& \sup \Big\{
P_{\varphi,u}(e)\,,
   \:\varphi \in C^\infty_{\mathrm{c}}\big((0,T)\times
               (-L,L)\big),\,|\varphi| \le 1,
\\
&&\phantom{ \sup \Big\{ P_{\varphi,u}(e)\,,}
          e \in C([0,1]),\, |e|\le 1
     \Big\}<+\infty
\end{eqnarray*}
Therefore the linear mapping $P_u^L: C([0,1])\times
C^\infty_{\mathrm{c}}\big((0,T)\times (-L,L)\big) \to \bb R$,
$P_u^L(e,\varphi):=P_{\varphi,u}(e)$ can be extended to a finite Borel
measure on $[0,1]\times \bb (0,T)\times (-L,L)$. The collection
$\{P_u^L\}_L$ defines a unique Radon measure $P_u$ on $[0,1]\times \bb
(0,T)\times \bb R$, since two elements of this collection coincide on
the intersection of their domains.  Recalling \eqref{e:2.10}, we thus
gather for each $\eta \in C^2([0,1])$, for each $\varphi \in
C^\infty_{\mathrm{c}}\big((0,T)\times \bb R\big)$ and for some
constant $C>0$ depending only on $f$
\begin{equation*}
  \Big|\int\!P_u(dv,dt,dx)\, \eta''(v)\varphi(t,x)\Big|
  =\Big|\wp_{\eta,u}(\varphi)\Big| 
  \le C \|\varphi\|_{C^1\big((0,T)\times \bb R\big)}
\int \!dv\,|\eta''(v)|
\end{equation*}
$P_u$ thus defines a linear continuous functional on $L_1([0,1])
\times C^1_{\mathrm{c}}\big((0,T)\times \bb R\big)$. This implies that
the Radon measure $P_u$ can be disintegrated as
$P_u=dv\,\varrho_u(v;dt,dx)$, for some bounded measurable map
$\varrho_u:[0,1]\to M\big((0,T)\times \bb R \big)$. From the
definition of $P_u$, we obtain for $\eta \in C^2([0,1])$, $\varphi \in
C^\infty_{\mathrm{c}}\big((0,T)\times \bb R \big)$ and
$\vartheta(v,t,x)=\eta(v) \varphi(t,x)$
\begin{eqnarray*}
  P_{\vartheta,u} &=& \wp_{\eta,u}(\varphi)
  = \int\!P_u(dv,dt,dx) \eta''(v) \varphi(t,x)
\\ 
&=& \int\!dv\,\varrho_u(v;dt,dx)\, \vartheta''(v,t,x)
\end{eqnarray*}
By linearity and density \eqref{e:2.13} holds for each entropy sampler
$\vartheta$.

\smallskip
\noindent (iii) $\Rightarrow$ (i).
It follows by choosing $\vartheta(v,t,x)= \eta(v) \varphi(t,x)$ in
equation \eqref{e:2.13} for $\varphi \in
C^\infty_{\mathrm{c}}\big((0,T)\times \bb R;[0,1] \big) $ and $\eta
\in C^2([0,1])$ with $0 \le \eta'' \le c$ for an arbitrary $c>0$.
\qed

\smallskip
\noindent
\textbf{Proof of Theorem~\ref{t:ecne}, item (ii): equicoercivity
  of $H_\varepsilon$.} 
The equicoercivity of $H_\varepsilon$ w.r.t. the topology generated by
the $d_{\mc U}$-distance \eqref{e:2.2} follows from
Lemma~\ref{l:equic}. It remains to show that, if $u^\varepsilon$ is
such that $H_\varepsilon(u^\varepsilon)$ is bounded uniformly in
$\varepsilon$, then $\{u^\varepsilon\}$ is precompact in
$L_{1,\mathrm{loc}}([0,T]\times \bb R)$. By equicoercivity of $\{\mc
I_\varepsilon\}$, the sequence $\{\mu^\varepsilon\}$ defined by
$\mu_{t,x}^\varepsilon=\delta_{u^\varepsilon(t,x)}$ is precompact in
$\mc M$. Therefore we have only to show that any limit point $\mu \in
\mc M$ of $\{\mu^\varepsilon\}$ has the form
$\mu_{t,x}=\delta_{u(t,x)}$ for some $u \in \mc X$, to obtain the
existence of limit points for $\{u^\varepsilon\}$ in $\mc X$. This is
implied by a compensated compactness argument due to Tartar, see
\cite[Ch. 9]{Se}, provided that there is no interval where $f$ is
affine, and that, for any entropy - entropy flux pair $(\eta,q)$, the
sequence $\{\eta(u^\varepsilon)_t+q(u^\varepsilon)_x\}$ is precompact
in $H^{-1}_{\mathrm{loc}}\big([0,T]\times \bb R\big)$.  Let us show
the latter. By \eqref{e:3.5}, there exists $C>0$ such that for each
$\varphi \in C^\infty_{\mathrm{c}}\big((0,T)\times (-L,L)\big)$
\begin{eqnarray*}
&&
\big|\langle \langle \eta(u^\varepsilon)_t 
                   +q(u^\varepsilon)_x,\varphi\rangle\rangle\big|
\\
&& 
\qquad\quad 
\le 
\frac{\varepsilon}2 \Big|\langle \langle \eta''(u^\varepsilon) 
                D(u^\varepsilon)u^\varepsilon_x, \varphi\,u^\varepsilon_x
\rangle \rangle\Big|
+\frac{\varepsilon}2 \Big|\langle \langle
\eta'(u^\varepsilon)D(u^\varepsilon) 
                u^\varepsilon_x,\varphi_x \rangle \rangle\Big|
\\
&& \qquad\quad  \phantom{\le}
  +\Big|\langle \langle \eta''(u^\varepsilon)
                               \sigma(u^\varepsilon)\,u^\varepsilon_x, 
              \Psi^{\varepsilon,u^\varepsilon}_x \varphi \rangle \rangle\Big|
    +\Big|\langle \langle \eta'(u^\varepsilon)\sigma(u^\varepsilon)\, 
          \Psi^{\varepsilon,u^\varepsilon}_x, \varphi_x \rangle \rangle\Big|
\\
&&
\qquad\quad 
\le   C\big[1+H_\varepsilon(u^\varepsilon)\big]
     \Big[\varepsilon \int_{[0,T]\times [-L,L]}\!dt \, dx \,
                                          (u^\varepsilon_x)^2\Big] 
\: \|\varphi\|_{L_{\infty}([0,T]\times \bb R)}
\\
&& \qquad\quad  \phantom{\le}
+C\Big[\varepsilon H_\varepsilon(u^\varepsilon)
           +\varepsilon^2\int_{[0,T]\times [-L,L]} \! dt \,dx\,  
      (u^\varepsilon_x)^2 \Big]^{1/2}
\: \|\varphi_x\|_{L_2([0,T]\times \bb R)}
\end{eqnarray*}
By the bound \eqref{e:3.6}, $\eta(u^\varepsilon)_t+q(u^\varepsilon)_x$
is the sum of a term bounded in $L_{1,\mathrm{loc}}([0,T]\times \bb
R)$ and a term vanishing in $H^{-1}_{\mathrm{loc}}([0,T]\times \bb R)$
as $\varepsilon \to 0$. By Sobolev compact embedding and boundedness
of $\eta$, $q$, the sequence
$\{\eta(u^\varepsilon)_t+q(u^\varepsilon)_x\}$ is compact in
$H^{-1}_{\mathrm{loc}}([0,T]\times \bb R)$.
\qed

\smallskip
\noindent
\textbf{Proof of Theorem~\ref{t:ecne}, item (i): $\Gamma$-liminf
  inequality.}
Let $\{u^\varepsilon\}$ be a sequence converging to $u$ in $\mc X$. If
$u$ is not a weak solution to \eqref{e:1.1}, by Theorem~\ref{t:2.1} we
have $\varliminf_{\varepsilon \to 0} I_\varepsilon(u^\varepsilon) \ge
\mc I(\delta_u)>0$, and therefore $\varliminf_{\varepsilon \to 0}
H_\varepsilon(u^\varepsilon)=+\infty$. Let now $u$ be a weak solution
to \eqref{e:1.1}. With no loss of generality we can suppose
$H_\varepsilon(u^\varepsilon) \le C_H$. We now consider an entropy
sampler -- entropy sampler flux pair $(\vartheta,Q)$ such that
\begin{equation}
\label{e:thetaineq}
0\le \sigma(v) \vartheta''(v,t,x) \le D(v),
\qquad (v,t,x) \in [0,1] \times (0,T)\times \bb R
\end{equation}
We also let $\varphi^\varepsilon(t,x)=\varepsilon
\vartheta'(u^\varepsilon(t,x),t,x)$, and introduce the short hand
notation $\big(\vartheta'(u^\varepsilon)\big)(t,x)\equiv
\vartheta'(u^\varepsilon(t,x),t,x)$,
$\big(\vartheta''(u^\varepsilon)\big)(t,x)\equiv
\vartheta''(u^\varepsilon(t,x),t,x)$, $\big((\partial_x
\vartheta')(u^\varepsilon)\big)(t,x) \equiv \big(\partial_x
\vartheta'\big)(u^\varepsilon(t,x),t,x)$. As we assumed
$H_\varepsilon(u^\varepsilon)<+\infty$, $u^\varepsilon_x$ is locally square
integrable, see \eqref{e:2.6}, and since $\vartheta$ is compactly
supported we have $\varphi^\varepsilon_x= \varepsilon
\vartheta''(u^\varepsilon)\, u^\varepsilon_x +\varepsilon (\partial_x
\vartheta')(u^\varepsilon) \in L_2([0,T]\times \bb R)$. The
representation \eqref{e:3.7} of
$\ell^{u^\varepsilon}_\varepsilon(\varphi^\varepsilon)$ thus holds, and
recalling \eqref{e:2.11} we get
\begin{eqnarray*}
&& \!\!
H_\varepsilon(u^\varepsilon) \ge \varepsilon^{-1} 
         \ell^{u^\varepsilon}_\varepsilon(\varphi^\varepsilon)
        - \frac{\varepsilon^{-1}}{2}
         \|\varphi^\varepsilon\|_{\mc D^1_{\sigma(u^\varepsilon)}}^2 
\\
&&\;
=  \langle \langle u^\varepsilon_t, 
            \vartheta'(u^\varepsilon)\rangle \rangle
      +\langle \langle f(u^\varepsilon)_x,
             \vartheta'(u^\varepsilon) \rangle \rangle
     + \frac{\varepsilon}{2} \langle \langle
          D(u^\varepsilon) u^\varepsilon_x,
       \vartheta''(u^\varepsilon) u^\varepsilon_x \rangle \rangle
\\
&& \; \phantom{=}
+ \frac{\varepsilon}{2} \langle \langle
          D(u^\varepsilon) u^\varepsilon_x,
       \big(\partial_x \vartheta'\big)(u^\varepsilon) \rangle \rangle
-\frac{\varepsilon}{2} \langle \langle 
            \sigma(u^\varepsilon) \vartheta''(u^\varepsilon) u^\varepsilon_x,
               \vartheta''(u^\varepsilon) u^\varepsilon_x \rangle \rangle
\\ 
&&\; \phantom{=}
- \varepsilon \langle \langle 
            \sigma(u^\varepsilon) \vartheta''(u^\varepsilon) u^\varepsilon_x,
  \big(\partial_x \vartheta'\big)(u^\varepsilon) \rangle \rangle
-\frac{\varepsilon}{2} \langle \langle 
            \sigma(u^\varepsilon)\big(\partial_x
            \vartheta'\big)(u^\varepsilon),
     \big(\partial_x \vartheta'\big)(u^\varepsilon) \rangle \rangle
\\ 
&& \;
= -\int \! dt\,dx\,
  \Big[\big(\partial_t \vartheta)\big(u^\varepsilon(t,x),t,x \big) 
       + \big(\partial_x Q\big)\big(u^\varepsilon(t,x),t,x \big)\Big]
\\ 
&&\; \phantom{=}
+\frac{\varepsilon}{2} \langle \langle
          D(u^\varepsilon) -\sigma(u^\varepsilon)\vartheta''(u^\varepsilon),
       \vartheta''(u^\varepsilon) (u^\varepsilon_x)^2 \rangle \rangle
+ \frac{\varepsilon}{2} \langle \langle
          D(u^\varepsilon) u^\varepsilon_x,
       \big(\partial_x \vartheta'\big)(u^\varepsilon) \rangle \rangle
\\ 
&&\; \phantom{=}
- \varepsilon \langle \langle 
            \sigma(u^\varepsilon) \vartheta''(u^\varepsilon) u^\varepsilon_x,
  \big(\partial_x \vartheta'\big)(u^\varepsilon) \rangle \rangle
-\frac{\varepsilon}{2} \langle \langle 
            \sigma(u^\varepsilon)\big(\partial_x
            \vartheta'\big)(u^\varepsilon),
     \big(\partial_x \vartheta'\big)(u^\varepsilon) \rangle \rangle
\end{eqnarray*}
By the bound \eqref{e:3.6}, the last three terms in the above formula
vanish as $\varepsilon \to 0$, while $\langle \langle [
D(u^\varepsilon) -\sigma(u^\varepsilon)\vartheta''(u^\varepsilon)]
u^\varepsilon_x, \vartheta''(u^\varepsilon) u^\varepsilon_x \rangle
\rangle \ge 0$ for each entropy sampler $\vartheta$ satisfying
\eqref{e:thetaineq}. Therefore, taking the limit $\varepsilon \to 0$
and optimizing over $\vartheta$
\begin{eqnarray*}
&&  \varliminf_{\varepsilon \to 0} H_\varepsilon(u^\varepsilon) 
\\
&&\qquad
 \ge \sup_{\vartheta}  \varliminf_{\varepsilon \to 0} -\int \! dt\,dx\,
  \Big[\big(\partial_t \vartheta)\big(u^\varepsilon(t,x),t,x \big) 
  + \big(\partial_x Q\big)\big(u^\varepsilon(t,x),t,x \big)\Big]
\quad
\\
&&\qquad
  =  \sup_{\vartheta} P_{\vartheta,u}
\end{eqnarray*}
where the supremum is taken on the $\vartheta \in
C^{2,\infty}_{\mathrm{c}}\big([0,1]\times (0,T)\times \bb R\big)$
satisfying \eqref{e:thetaineq}. Recalling that we assumed the l.h.s.\
of this formula to be finite, we next show that this inequality
implies that $u \in \mc E$, and that the r.h.s.\ is equal to $H(u)$.
By taking $\vartheta(v,t,x) = \eta(v) \varphi(t,x)$ for some $\varphi
\in C^\infty_{\mathrm{c}}\big([0,T]\times \bb R;[0,1]\big)$ and
entropy $\eta$ such that $0\le \sigma(v) \eta''(v) \le D(v)$, we get
$\wp_{\eta,u}(\varphi) \le \varliminf_{\varepsilon}
H_\varepsilon(u^\varepsilon)$. Optimizing over $\varphi$ it follows that $u$
fulfills condition (i) in Proposition~\ref{p:kin} with $c=\min_{v}
D(v)/\sigma(v)>0$, and thus $u \in \mc E$. By (iii) in
Proposition~\ref{p:kin} and monotone convergence we then get
\begin{eqnarray}
\label{e:Hsup}
\nonumber
\varliminf_{\varepsilon \to 0} H_\varepsilon(u^\varepsilon) 
&  \ge &
 \sup_{\vartheta} P_{\vartheta,u} 
= \sup_{\vartheta} \int\! dv\,\varrho_u(v;dt,dx)\,
            \vartheta''(v,t,x)
\\
&=&
 \int\! dv\,\varrho_u^+(v;dt,dx)\, \frac{D(v)}{\sigma(v)} =H(u)
\end{eqnarray}
which concludes the proof. \qed

\smallskip
\begin{lemma2}
\label{l:l1}
Let $f \in C^2([0,1])$ and assume that there is no interval where $f$
is affine. Then entropy-measure solutions to \eqref{e:1.1} 
belong to the space
$C\big([0,T];L_{1,\mathrm{loc}}(\bb R)\big)$.
Let furthermore
\begin{equation}
\label{e:vpm}
V_f^+:=\max_{v \in [0,1]} f'(v)  \qquad V_f^-:=\min_{v \in [0,1]} f'(v)
\end{equation}
Then for each $u \in \mc E$, $x \in \bb R$, $V>V_f^+$ or $V<V_f^-$
\begin{equation}
\label{e:xcont}
\lim_{\zeta \to 0} \int \!dt 
     \big|u\big(t, x+\zeta+V\,t\big)-u\big(t, x+V\,t\big) \big|=0
\end{equation}
\end{lemma2}

\proof 
With the same hypotheses of this lemma, in \cite[Sect.~4]{CR} it is
shown that if a weak solution $u$ to \eqref{e:1.1} is such that
$\wp_{f,u}$ is a Radon measure, then, for each $L>0$ and $t \in
[0,T)$, $\lim_{s \downarrow t} \int_{[-L,L]}\!|u(s,x)-u(t,x)|dx = 0$.
Therefore, by item (ii) in Proposition~\ref{p:kin}, entropy-measure
solutions enjoy this property. Since the set $\mc E$ of
entropy-measure solutions is invariant under the symmetry $(t,x)
\mapsto (-t,-x)$, the same holds true also for $s \uparrow t$, and
thus $\mc E \subset C\big([0,T];L_{1,\mathrm{loc}}(\bb R)\big)$.

If $u$ is an entropy-measure solution to the conservation law
\eqref{e:1.1}, then $u^{V,\pm}(t,x):=u(\pm t,x\pm Vt)$ is an
entropy-measure solution to the conservation law with flux $f^\pm$,
where $f^\pm (w)=f(w)\mp V w$. With no loss of generality, we can thus
prove \eqref{e:xcont} only in the case $V=0$ with the assumption
$V_f^->0$. In this case $f$ is invertible on its range $[a,b]$, and we
let $g \in C^2([a,b])$ be its inverse. We define $v:\bb R \times [0,T]
\mapsto [a,b]$ by $v(x,t)=f(u(t,x))$. Then $v$ satisfies
\begin{equation}
\label{e:consg}
v_x+g(v)_t=0
\end{equation}
Furthermore, if $l,m \in C^2([a,b])$ satisfy $m'=l'g'$, then by chain
rule $l(v)_x+m(v)_t= \wp_{\eta,u}$, where $\eta(w):=\int^w\!
dz\,l'(f(z))$. Therefore $v$ is an entropy-measure solution to
\eqref{e:consg}, and by the first part of this lemma
\begin{equation*}
  \lim_{\zeta \to 0} \int \! ds\, \big| v(x+\zeta,s)-v(x,s)\big|=0
\end{equation*}
The result then follows by recalling $u(t,x)=g(v(x,t))$.
\qed

\smallskip
\noindent
\textbf{Proof of Theorem~\ref{t:ecne}, item (iii): 
$\Gamma$-limsup inequality.}
Given an nice (w.r.t.\ $\sigma$) solution $\tilde{u} \in \mc
S_\sigma$, let $E^{\pm}$ be as in Definition~\ref{d:splittable}. We
want to construct a recovery sequence $\{u^\varepsilon\} \subset \mc
X$ that converges to $\tilde{u}$ in $\mc X$ as $\varepsilon \to 0$,
and such that $\varlimsup_{\varepsilon} H_\varepsilon(u^\varepsilon)
\le H(\tilde{u})$. We split the proof in four steps.  In
\textit{Step~1} we build a suitable family of rectangles contained in
$[0,T]\times \bb R$. In \textit{Step~2}, for $\varepsilon,\,\delta,\,L
\ge 1$, we introduce two collections
$\{v^{\varepsilon,\delta,L,\pm}\}$ of auxiliary functions on
$[0,T]\times \bb R$. In \textit{Step~3}, for $N \in \bb N$ we define a
collection $\{u^{\varepsilon,\delta,N,L}\} \subset \mc X$, and we
prove
\begin{equation}
\label{e:ucost}
\varlimsup_{\delta \to 0} \varlimsup_{\varepsilon \to 0}
H_\varepsilon(u^{\varepsilon,\delta,N,L}) \le H(\tilde{u})
\end{equation}
In particular $\{u^{\varepsilon,\delta,N,L}\}$ is precompact in $\mc
X$. In \textit{Step~4} we show that any limit point of
$\{u^{\varepsilon,\delta,N,L}\}$ coincides with $\tilde{u}$ in $\mc X$,
provided we consider the limit in $\varepsilon,\,\delta,\,N,\,L$ in a
suitable order. More precisely we show
\begin{equation}
\label{e:uconv}
\lim_{L \to \infty} \lim_{N \to \infty} \lim_{\delta \to 0}
\lim_{\varepsilon \to 0} u^{\varepsilon,\delta,N,L} = \tilde{u}
\end{equation}
By \eqref{e:ucost} and \eqref{e:uconv} it follows that there exist
subsequences $\{\delta^\varepsilon\},\,\{L^\varepsilon\} \subset
(0,+\infty)$ and $\{N^\varepsilon\} \subset \bb N$ such that
$u^\varepsilon:= u^{\varepsilon,\delta^\varepsilon,N^\varepsilon,L^\varepsilon}$
provides the required recovery sequence for $\tilde{u}$.

Throughout this proof, we assume $f'$ to be uniformly positive in
$[0,1]$, namely that $V_f^-$, as defined in \eqref{e:vpm}, is
positive. As noted in the proof of Lemma~\ref{l:l1}, this assumption
is not restrictive. Note also that the calculations carried out below
make sense also if $E^+=\emptyset$ or $E^-=\emptyset$.

\smallskip
\noindent
\textit{Step 1.}
For each $t$ such that $\big(\{t\} \times [-L,L]\big) \cap E^+ \cap
E^- = \emptyset $, the compact sets $\big(\{t\}\times [-L,L]\big) \cap
E^{\pm}$ are disjoint, hence strictly separated. By (ii) in
Definition~\ref{d:splittable}, there exists a countable collection of
pairwise disjoint time intervals $\{(s_i^L,t_i^L)\}_{i \in \bb N}$,
with $(s_i^L,t_i^L) \subset (0,T)$ such that $\tau^L:=\cup_i
(s_i^L,t_i^L)$ is dense in $[0,T]$, and for each $i \in \bb N$ the two
sets $E^{L,\pm}_i:= \big((s_i^L,t_i^L) \times [-L,L]\big) \cap E^\pm $
are strictly separated. By splitting each of these intervals in a
finite number of intervals, with no loss of generality we can assume
\begin{equation}
\label{e:tsmall}
t_i^L-s_i^L <
\frac{1}{10+8\,V_f^+}\,\mathrm{distance}\big(E^{L,+}_i,E^{L,-}_i \big)
\end{equation}
where $V_f^+$ is defined in \eqref{e:vpm}, and it coincides with the
Lipschitz constant of $f$ since we supposed $V_f^->0$.

For $i \in \bb N$ let $n_i^L \in \bb N$ be such that
\begin{equation}
\label{e:nil}
\frac{L}{n_i^L} \le \frac{1}{10} \min \big\{1,
\mathrm{distance}\big(E^{L,+}_i,E^{L,-}_i \big) \big\}
\end{equation}
and consider the rectangles $R_{i,j}^L:= (s_i^L,t_i^L) \times
(\frac{j}{n_i^L} L,\frac{j+1}{n_i^L} L )$, for
$j=-n_i^L$, $-n_i^L+1$, $\ldots$, $n_i^L-1$. By the definition
\eqref{e:nil} of $n_i^L$ and condition \eqref{e:tsmall}, for $j=-n_i^L+1,-n_i^L+2,\ldots, n_i^L-2$
\begin{equation}
\label{e:xsmall}
\mathrm{diameter}( R_{i,j-1}^L \cup  R_{i,j}^L \cup R_{i,j+1}^L \cup R_{i,j+2}^L ) 
 < \frac{1}{2} \mathrm{distance}\big(E^{L,+}_i,E^{L,-}_i \big)
\end{equation}
In particular each $R_{i,j}^L$ has nonempty intersection with at most
one of the sets $E^+$, $E^-$. We define
\begin{equation}
\label{e:ril}
R_i^{L,\pm} := \bigcup_{\substack{
                         j = -n_i^L+1, \\
                   j\,:\: (R_{i,j-1}^L \cup  R_{i,j}^L \cup
                   R_{i,j+1}^L \cup R_{i,j+2}^L) \cap E^\mp = \emptyset } }^{n_i^L-2} R_{i,j}^L
\end{equation}
and for $N \in \bb N$
\begin{equation}
\label{e:rnl}
R^{N,L,\pm}:= \cup_{i=1}^N R_i^{L,\pm} 
\qquad \qquad 
R^{L,\pm}:=\cup_N R^{N,L,\pm}
\end{equation}

Note that by \eqref{e:tsmall} and \eqref{e:nil}
\begin{eqnarray}
\label{e:rilkruz}
\nonumber
&& R_{i,j}^L \subset 
 \Big\{ (r,x)\,:\:  s_i^L<r<t_i^L,
\\
&&\phantom{ R_{i,j}^L \subset  \Big\{ }
         \frac{j-1}{n_i^L}L + V_f^+ (r-s_i^L) 
                       \le  x \le
        \frac{j+2}{n_i^L}L - V_f^+ (r-s_i^L)
\Big\}
\qquad\quad
\end{eqnarray}
and by \eqref{e:xsmall}
\begin{equation}
\label{e:rlarge}
R^{L,+}\cup R^{L,-}= \bigcup_i \bigcup_{j=-n_i^L+1}^{n_i^L-2} R_{i,j}^L
\end{equation}

\smallskip
\noindent
\textit{Step 2.}
For $L \ge 1$ and $\delta \in (0,1/2)$, let $\tilde{u}^{\delta,L} \in \mc
X$ be defined by
\begin{equation}
\label{e:udl}
\tilde{u}^{\delta,L}(t,x):=
\begin{cases}
\tilde{u}(t,x) 
   & \text{if $|x|\le L$ and $\tilde{u}(t,x) \in
     [\delta,1-\delta]$}
\\
\delta
   & \text{if $|x|\le L$ and $\tilde{u}(t,x) \le \delta$}
\\
1-\delta 
   & \text{if $|x|\le L$ and $\tilde{u}(t,x) \ge 1-\delta$}
\\
1/2 & \text{if $|x|>L$}
\end{cases}
\end{equation}
For $\varepsilon>0$, $i \in \bb N$, we define
$v^{\varepsilon,\delta,L,-}_{i}:(s_i^L,t_i^L) \times \bb R \to \bb R$ as
the solution to the forward-parabolic Cauchy problem
\begin{equation}
\label{e:v-}
\begin{cases}
{\displaystyle
v_t+ f(v)_x=\frac{\varepsilon}2 \big(D(v) v_x\big)_x
}
\\
{\displaystyle
v(s_i^L)=\tilde{u}^{\delta,L}(s_i^L)
}
\end{cases}
\end{equation}
and $v^{\varepsilon,\delta,L,+}_{i}:(s_i^L,t_i^L) \times \bb R \to \bb R$
as the solution to the backward-parabolic Cauchy problem
\begin{equation}
\label{e:v+}
\begin{cases}
{\displaystyle
v_t+ f(v)_x=-\frac{\varepsilon}2 \big(D(v) v_x\big)_x
}
\\
{\displaystyle
v(t_i^L)=\tilde{u}^{\delta,L}(t_i^L)
}
\end{cases}
\end{equation}
We define $v^{\varepsilon,\delta,L,\pm}:\tau^L \times \bb R \to
\bb R$ by requiring
$v^{\varepsilon,\delta,L,\pm}(r,x)=v^{\varepsilon,\delta,L,\pm}_{i}(r,x)$
for $r \in (s_i^L,t_i^L)$. Note that $v^{\varepsilon,\delta,L,\pm} \in
C\big(\tau^L;U\big)$ and $v^{\varepsilon,\delta,L,\pm}(t,x) \in
[\delta,1-\delta]$ by maximum principle.  Furthermore
$v^{\varepsilon,\delta,L,\pm}_x \in L_{2,\mathrm{loc}}(\tau^L \times
\bb R)$, and indeed by standard parabolic estimates
\begin{eqnarray}
\label{e:vest}
\nonumber
&& \varepsilon \int_{R^{N,L,\pm}}\!dr\,dx\, 
\big( v^{\varepsilon,\delta,L,\pm}_x(r,x)\big)^2
\\
&&
\qquad\qquad
\le
\sum_{i=1}^N \varepsilon \int_{[s_i^L,t_i^L]\times [-L,L]} \!dr \, dx\,
        \big( v^{\varepsilon,\delta,L,\pm}_x(r,x)\big)^2 \le C^{N,L}
\qquad
\end{eqnarray}
for some constant $C^{N,L}>0$ independent of $\varepsilon$ and $\delta$.

We claim
\begin{equation}
\label{e:vconv}
\lim_{\delta \to 0} \lim_{\varepsilon \to 0} \int_{R^{N,L,\pm}}\!dr\,dx\,
\big|v^{\varepsilon,\delta,L,\pm}(r,x)-\tilde{u}(r,x)\big|=0
\end{equation}
We show \eqref{e:vconv} for $v^{\varepsilon,\delta,L,-}$. The analogous
statement for $v^{\varepsilon,\delta,L,+}$ follows by the fact that the
set $\mc S_\sigma$ is invariant w.r.t.\ the symmetry $(t,x) \mapsto
(-t,-x)$, while the supports of $\varrho^{\pm}_u$ are exchanged under
this symmetry. By the well known results of convergence of the
vanishing viscosity approximations to conservation laws (and as it
also follows from the $\Gamma$-liminf inequality in
Theorem~\ref{t:ecne} item (i))
\begin{equation}
\label{e:vv}
  \lim_{\varepsilon \to 0} \int_{[s_i^L,t_i^L]\times [-L,L]}\!dr \,dx\, 
  \big|v^{\varepsilon,\delta,L,-}(r,x)-\bar{u}_i^{\delta,L}(r-s_i^L,x)\big|=0
\end{equation}
where $\bar{u}_i^{\delta,L}$ is the Kruzkov solution to \eqref{e:1.1}
with initial condition
$\bar{u}_i^{\delta,L}(0,\cdot)=\tilde{u}^{\delta,L}(s_i^L,\cdot)$. On
the other hand, by the definition \eqref{e:ril} of $R_i^{L,-}$, if $j$
is such that $R_{i,j}^L \subset R_i^{L,-}$, then $\tilde{u}$ is
entropic in the rectangle $(s_i^L,t_i^L) \times (\frac{j-1}{n_i^L} L
,\frac{j+2}{n_i^L} L )$, namely $\wp_{\eta,\tilde{u}}(\varphi) \le 0$
for each convex entropy $\eta$ and each positive test function
$\varphi$ compactly supported in $(s_i^L,t_i^L) \times
(\frac{j-1}{n_i^L} L ,\frac{j+2}{n_i^L} L )$. Therefore, by Kruzkov
theorem \cite{Se}
\begin{eqnarray*}
&&
\varlimsup_{\delta \to 0} \sup_{s_i^L \le r \le t_i^L} 
    \int_{\frac{j-1}{n_i^L}L + V_f^+ (r-s_i^L)}^{\frac{j+2}{n_i^L}L-V_f^+
(r-s_i^L)}\!\!dx 
             \big|\bar{u}_i^{\delta,L}(r-s_i^L,x)-\tilde{u}(r,x) \big|
\\ 
&& \qquad\quad
\le 
  \varlimsup_{\delta \to 0} 
     \int_{\frac{j-1}{n_i^L}L}^{\frac{j+2}{n_i^L}L} \!\!\!dx 
      \big| \bar{u}_i^{\delta,L}(0,x) - \tilde{u}(s_i^L,x) \big|
\\
&& \qquad\quad
= \varlimsup_{\delta \to 0} 
    \int_{\frac{j-1}{n_i^L}L}^{\frac{j+2}{n_i^L}L} \!\!\!dx 
      \big| \tilde{u}^{\delta,L}(s_i^L,x) - \tilde{u}(s_i^L,x) \big|
=0
\end{eqnarray*}
and thus, fixed $N \in \bb N$, by \eqref{e:rilkruz} the convergence
claimed in \eqref{e:vconv} holds on each $R_{i,j}^L$ for each $i \le
N$ and each $j$ such that $R_{i,j}^L \subset R_i^{L,-}$, and therefore
on $R^{N,L,-}$ itself.

Next we claim that for each $L \ge 1$, $N \in \bb N$ and $\varphi \in
C^\infty_{\mathrm{c}}\big(R^{N,L,+};[0,1]\big)$
\begin{equation}
\label{e:vcost}
\varlimsup_{\delta \to 0} \varlimsup_{\varepsilon \to 0}
\frac{\varepsilon}{2} \Big\langle\Big\langle 
   D(v^{\varepsilon,\delta,L,+}) v^{\varepsilon,\delta,L,+}_x, 
   \varphi
     \frac{D(v^{\varepsilon,\delta,L,+})}{\sigma(v^{\varepsilon,\delta,L,+})} 
      v^{\varepsilon,\delta,L,+}_x 
  \Big\rangle\Big\rangle \le H(\tilde{u})
\end{equation}
Note that the l.h.s.\ of this formula is well defined, since $\delta
\le v^{\varepsilon,\delta,L,+} \le 1-\delta$ and thus $\sigma(
v^{\varepsilon,\delta,L,+})$ is uniformly positive. For each $\varphi \in
C^\infty_{\mathrm{c}}\big(([0,T]\times \bb R) \setminus
E^-;[0,1]\big)$ and $\eta \in C^2\big([0,1]\big)$ such that $\sigma
\eta'' \le D$ we have
\begin{equation}
\label{e:Hlarger}
H(\tilde{u}) \ge \int\!dw\, \varrho_{\tilde{u}}(w;dt,dx)\, \eta''(w)
\varphi(t,x) = \wp_{\eta,\tilde{u}}(\varphi)
\end{equation}
By \eqref{e:v+} and \eqref{e:vest} for each $\eta \in C^2([0,1])$, $N
\in \bb N$ and $\varphi \in C^\infty_{\mathrm{c}}\big(R^{N,L,+}\big)$
\begin{equation}
\label{e:etaconv}
\lim_{\delta \to 0} \lim_{\varepsilon \to 0} 
   \frac{\varepsilon}{2} \big\langle \big\langle 
   D(v^{\varepsilon,\delta,L,+}) v^{\varepsilon,\delta,L,+}_x, 
   \varphi \eta''(v^{\varepsilon,\delta,L,+})v^{\varepsilon,\delta,L,+}_x 
  \big\rangle\big\rangle 
= \wp_{\eta,\tilde{u}}(\varphi)
\end{equation}
This implies \eqref{e:vcost} if $\sigma$ is uniformly positive on
$[0,1]$, since we can evaluate \eqref{e:etaconv} on an entropy $\eta$
such that $\eta''=D/\sigma$ and use the trivial bound
\eqref{e:Hlarger}. On the other hand, if $\sigma(0)=0$, resp.\ if
$\sigma(1)=0$, then by condition (iii) in
Definition~\ref{d:splittable}, we have that $\tilde{u}(t,x) \ge
\zeta_L$, resp.\ $\tilde{u}(t,x) \le 1-\zeta_L$, for a.e.\ $(t,x) \in
(0,T)\times (-L,L)$ and for some $\zeta_L>0$. By the definition of
$\tilde{u}^{\delta,L}$ and maximum principle, we have also
$v^{\varepsilon,\delta,L,+} \ge \zeta_L$, resp.\ $v^{\varepsilon,\delta,L,+}
\le 1-\zeta_L$, and thus \eqref{e:vconv} follows by evaluating
\eqref{e:etaconv} on an entropy $\eta$ such that
$\eta''(w)=D(w)/\sigma(w)$ for all $w \ge \zeta_L$, resp.\ $w \le
1-\zeta_L$.

\smallskip
\noindent\textit{Step 3.}  In this step, with a little abuse of
notation, we denote by $f$ and $D$ two bounded continuous functions on
$\bb R$, such they their restrictions to $[0,1]$ coincide with $f$ and
$D$, and $f$ is uniformly Lipschitz and $D$ uniformly positive.  We
also let $\sigma^{\delta} \in C^\alpha([0,1])$ be such that
$\sigma^\delta(w)=\sigma(w)$ for $w \in [\delta,1-\delta]$,
$\sigma^\delta(w) \le \sigma(w)$ for $w \in [0,1]$, and
$\sigma^\delta(w)=0$ for $w \le \delta/2$ or $w \ge 1-\delta/2$.

For $L \ge 1$ and $N \in \bb N$, let $\Xi^{N,L} \in
C^\infty_{\mathrm{c}}\big(R^{N,L,+};[0,1]\big)$, and define
\begin{equation}
\begin{array}{lcl}
\label{e:fnl}
& &
P^{N,L,+}:= \mathrm{Interior} \Big( 
    \big\{ (t,x) \in R^{N,L,+}\,:\: \Xi^{N,L}(t,x)=1 \big\} \Big)
\\
& &
P^{N,L,-}:=  \mathrm{Interior} \Big( 
    \big\{ (t,x) \in R^{N,L,-}\,:\: \Xi^{N,L}(t,x)=0 \big\} \Big)
\end{array}
\end{equation}
For each fixed $L \ge 1$, we require the sequence $\{\Xi^{N,L}\}$ to
be increasing in $N$ and such that
\begin{equation}
\label{e:xilarge}
\cup_{N} P^{N,L,+} =  R^{L,+}
\end{equation}

For $\delta,\,L \ge 1$ and $N\in \bb N$ define
$u^{\varepsilon,\delta,N,L}:[0,T]\times \bb R \to \bb R$ as the
solution to the Cauchy problem
\begin{equation}
\label{e:udef}
\begin{cases}
{\displaystyle
  u_t+f(u)_x=\frac{\varepsilon}{2} \big(D(u)u_x \big)_x
}
\\
{\displaystyle
\phantom{u_t+f(u)_x=}
   -\varepsilon \Big[\Xi^{N,L}
     \frac{\sqrt{\sigma^\delta(u)}}{\sqrt{\sigma(v^{\varepsilon,\delta,L,+})}}
     D(v^{\varepsilon,\delta,L,+}) v^{\varepsilon,\delta,L,+}_x \Big]_x
}
\\
{\displaystyle
u(0,x)=\tilde{u}^{\delta,L}(0,x) \qquad x\in \bb R
}
\end{cases}
\end{equation}
Note that the term in square brackets in \eqref{e:udef} is
well-defined since $v^{\varepsilon,\delta,L,+}$ is well-defined on the
support of $\Xi^{N,L}$, and since $\delta \le v^{\varepsilon,\delta,L,+}
\le 1-\delta$, $\sigma(v^{\varepsilon,\delta,L,+})$ is uniformly
positive.

It is easily seen that the problem \eqref{e:udef} admits at least a
solution $u^{\varepsilon,\delta,N,L} \in L_{\infty}\big([0,T]\times\bb
R\big)$ with $u^{\varepsilon,\delta,N,L}_x \in
L_{2,\mathrm{loc}}\big([0,T]\times\bb R\big)$. By \eqref{e:udef} we
also gather
\begin{eqnarray*}
&& 
\big\| u^{\varepsilon,\delta,N,L}_t+f(u^{\varepsilon,\delta,N,L})_x -
\frac{\varepsilon}{2}
\big(D(u^{\varepsilon,\delta,N,L})u^{\varepsilon,\delta,N,L}_x \big)_x
\big\|_{\mc D^{-1}_{\sigma^\delta(u)}}^2
\\
&&\qquad\quad
= \varepsilon^2 \Big\langle \Big\langle
 D(v^{\varepsilon,\delta,L,+}) v^{\varepsilon,\delta,L,+}_x, 
   (\Xi^{N,L})^2
     \frac{D(v^{\varepsilon,\delta,L,+})}{\sigma(v^{\varepsilon,\delta,L,+})} 
      v^{\varepsilon,\delta,L,+}_x 
  \Big\rangle \Big\rangle <+\infty
\end{eqnarray*}
Therefore, replacing $\sigma$ with $\sigma^\delta$ in the statement of
Proposition~\ref{t:boun}, we have $\delta \le
u^{\varepsilon,\delta,N,L} \le 1-\delta$ and
$u^{\varepsilon,\delta,N,L} \in \mc X$. Since $(\Xi^{N,L})^2 \in
C^\infty_{\mathrm{c}}(R^{N,L,+};[0,1])$, by the same estimate and
\eqref{e:vcost}
\begin{eqnarray*}
&&
\varlimsup_{\delta} \varlimsup_\varepsilon
H_\varepsilon(u^{\varepsilon,\delta,N,L}) 
\\
&&
\qquad
=\varlimsup_\delta \varlimsup_\varepsilon \frac{\varepsilon}{2} 
\Big\langle\Big\langle
 D(v^{\varepsilon,\delta,L,+}) v^{\varepsilon,\delta,L,+}_x, 
   (\Xi^{N,L})^2 
\\
&&\qquad
\phantom{=\varlimsup_\delta \varlimsup_\varepsilon \frac{\varepsilon}{2} 
\Big\langle\Big\langle}
\times
\frac{\sigma^\delta(u^{\varepsilon,\delta,N,L})}
              {\sigma(u^{\varepsilon,\delta,N,L})}
     \frac{D(v^{\varepsilon,\delta,L,+})}{\sigma(v^{\varepsilon,\delta,L,+})} 
      v^{\varepsilon,\delta,L,+}_x
  \Big\rangle\Big\rangle
\\
&&\qquad 
\le \varlimsup_\delta \varlimsup_\varepsilon \frac{\varepsilon}{2} 
\Big\langle\Big\langle
 D(v^{\varepsilon,\delta,L,+}) v^{\varepsilon,\delta,L,+}_x, 
   (\Xi^{N,L})^2 \frac{D(v^{\varepsilon,\delta,L,+})}
                     {\sigma(v^{\varepsilon,\delta,L,+})} 
      v^{\varepsilon,\delta,L,+}_x 
  \Big\rangle\Big\rangle
\\
&&\qquad 
\le H(\tilde{u})
\end{eqnarray*}
so that \eqref{e:ucost} holds.

\smallskip
\noindent\textit{Step 4.}  Since $\{H_\varepsilon\}$ is equicoercive on
$\mc X$ and \eqref{e:ucost} holds, there exist
$\delta_0,\,\varepsilon_0 \equiv \varepsilon_0(\delta_0)$ small enough
and a compact set $\mc K_0 \subset \mc X$ such that
$u^{\varepsilon,\delta,N,L} \in \mc K_0$ for each $\varepsilon
<\varepsilon_0$, $\delta<\delta_0$, $N \in \bb N$ and $L \ge 1$. In
this step we show that any limit point $u$ of
$\{u^{\varepsilon,\delta,N,L}\}$ coincide with $\tilde{u}$, provided
the limits in $\varepsilon$, $\delta$, $N$ and $L$ are taken in a
suitable order, see \eqref{e:uconv}. This will conclude the proof.

Let $z^{\varepsilon,\delta,N,L,\pm}: \tau^L \times \bb R \to [-1,1]$,
$z^{\varepsilon,\delta,N,L,\pm}:=u^{\varepsilon,\delta,N,L}
-v^{\varepsilon,\delta,L,\pm}$. By \eqref{e:3.6}, \eqref{e:ucost} and
\eqref{e:vest}, for each $N \in \bb N$
\begin{equation}
\label{e:zbound}
\varepsilon \int_{R^{N,L,\pm}}\!dt\,dx\, (z^{\varepsilon,\delta,N,L,\pm}_x )^2
\le \tilde{C}^{N,L}
\end{equation}
for some constant $\tilde{C}^{N,L} >0$ independent of $\varepsilon$ and
$\delta$.

Since we will first perform the limit $\varepsilon \to 0$, we now fix
$\delta$, $N$, $L$ as above, and we drop for a few lines these
indexes, thus writing $u^\varepsilon \equiv
u^{\varepsilon,\delta,N,L}$, $v^{\varepsilon,\pm} \equiv
v^{\varepsilon,\delta,L,\pm}$, $z^{\varepsilon,\pm} \equiv
z^{\varepsilon,\delta,N,L,\pm}$, $\Xi \equiv \Xi^{N,L}$. Recalling the
definition \eqref{e:fnl}, by \eqref{e:udef} and \eqref{e:v-}, we have
weakly on $P^{N,L,-}$
\begin{equation*}
z^{\varepsilon,-}_t+\big(f(u^\varepsilon)-f(v^{\varepsilon,-})\big)_x
=\frac{\varepsilon}{2}\big(D(u^\varepsilon) z^{\varepsilon,-}_x \big)_x
  +\frac{\varepsilon}{2} \big(\big[D(u^\varepsilon)-D(v^{\varepsilon,-})\big]
                                      v^{\varepsilon,-}_x \big)_x
\end{equation*}
Let now $l\in C^2([-1,1])$ and $\varphi \in C^\infty_{\mathrm{c}}
\big(P^{N,L,-}\big)$. It follows
\begin{eqnarray}
\label{e:lfuncext1}
\nonumber
&&
-\langle \langle l(z^{\varepsilon,-}),\varphi_t\rangle \rangle
-\langle \langle f(u^\varepsilon)-f(v^{\varepsilon,-}), l'(z^{\varepsilon,-})
                \varphi_x \rangle \rangle
\\
\nonumber
&&\qquad \phantom{=}
-\langle \langle f(u^\varepsilon)-f(v^{\varepsilon,-}), l''(z^{\varepsilon,-})
                z^{\varepsilon,-}_x \varphi \rangle \rangle
\\
\nonumber
& & 
\qquad =
-\frac{\varepsilon}{2} \langle \langle D(u^\varepsilon) z^{\varepsilon,-}_x,
l''(z^{\varepsilon,-}) z^{\varepsilon,-}_x \varphi \rangle \rangle
- \frac{\varepsilon}{2} \langle \langle D(u^\varepsilon) z^{\varepsilon,-}_x,
 l'(z^{\varepsilon,-}) \varphi_x \rangle \rangle
\\
\nonumber
& &
\qquad \phantom{=}
-\frac{\varepsilon}{2} \langle \langle
\big[D(u^\varepsilon)-D(v^{\varepsilon,-})\big] v^{\varepsilon,-}_x,
l''(z^{\varepsilon,-}) z^{\varepsilon,-}_x \varphi \rangle \rangle
\\
& &
\qquad \phantom{=}
-\frac{\varepsilon}{2} \langle \langle
\big[D(u^\varepsilon)-D(v^{\varepsilon,-})\big] v^{\varepsilon,-}_x, 
l'(z^{\varepsilon,-}) \varphi_x \rangle \rangle
\end{eqnarray}

In the same fashion, by \eqref{e:v+}, weakly on $P^{N,L,+}$
\begin{eqnarray*}
z^{\varepsilon,+}_t+\big(f(u^\varepsilon)-f(v^{\varepsilon,+})\big)_x
& = & \frac{\varepsilon}{2}\big(D(u^\varepsilon) z^{\varepsilon,+}_x \big)_x
  +\frac{\varepsilon}{2} \big(\big[D(u^\varepsilon)-D(v^{\varepsilon,+})\big]
                                      v^{\varepsilon,+}_x \big)_x
\\
& &
 +\varepsilon \Big( \Big[\sqrt{\sigma(v^{\varepsilon,+})} 
              - \sqrt{\sigma^{\delta}(u^\varepsilon)} \Big]
\frac{D(v^{\varepsilon,+})}{\sqrt{\sigma(v^{\varepsilon,+})}} 
                            v^{\varepsilon,+}_x \Big)_x
\end{eqnarray*}
Since $v^{\varepsilon,+}$ takes values in $[\delta,1-\delta]$, we have
$\sigma^{\delta}(v^{\varepsilon,+})=\sigma(v^{\varepsilon,+})$ and thus, in
the same fashion as above, for each $l \in C^2([-1,1])$ and $\varphi
\in C^\infty_{\mathrm{c}}\big(P^{N,L,+}\big)$
\begin{eqnarray}
\label{e:lfuncext2}
&&
\nonumber 
-\big\langle \big\langle l(z^{\varepsilon,+}),\varphi_t
\big\rangle \big\rangle 
-  \big\langle\big\langle 
f(u^\varepsilon)-f(v^{\varepsilon,+}), l'(z^{\varepsilon,+})
                    \varphi_x 
\big\rangle\big\rangle
\\
&&
\nonumber 
\qquad \phantom{=}
-\big\langle\big\langle
f(u^\varepsilon)-f(v^{\varepsilon,+}), l''(z^{\varepsilon,+})
                z^{\varepsilon,+}_x \varphi 
\big\rangle\big\rangle
\\
\nonumber  
& & \qquad =
-\frac{\varepsilon}{2} 
\big\langle\big\langle D(u^\varepsilon) z^{\varepsilon,+}_x,
l''(z^{\varepsilon,+}) z^{\varepsilon,+}_x \varphi \big\rangle\big\rangle
- \frac{\varepsilon}{2} 
\big\langle \big\langle D(u^\varepsilon) z^{\varepsilon,+}_x,
 l'(z^{\varepsilon,+}) \varphi_x \big\rangle \big\rangle
\\
\nonumber 
& &
\qquad \phantom{=}
-\frac{\varepsilon}{2} \big\langle \big\langle
\big[D(u^\varepsilon)-D(v^{\varepsilon,+})\big] v^{\varepsilon,+}_x,
l''(z^{\varepsilon,+}) z^{\varepsilon,+}_x \varphi \big\rangle \big\rangle
\\
& &
\nonumber 
\qquad \phantom{=}
-\frac{\varepsilon}{2} \big\langle \big\langle
\big[D(u^\varepsilon)-D(v^{\varepsilon,+})\big] v^{\varepsilon,+}_x, 
l'(z^{\varepsilon,+}) \varphi_x \big\rangle \big\rangle
\\
\nonumber 
& &
\qquad \phantom{=}
-\varepsilon \Big\langle\Big\langle 
  \big[\sqrt{\sigma^\delta(v^{\varepsilon,+})}
    -\sqrt{\sigma^\delta(u^\varepsilon)} \big]
   \frac{ D(v^{\varepsilon,+})} {\sqrt{\sigma(v^{\varepsilon,+})}}
 v^{\varepsilon,+}_x, 
    \varphi l''(z^{\varepsilon,+})z^{\varepsilon,+}_x \Big\rangle\Big\rangle
\\ 
& &
\qquad \phantom{=}
-\varepsilon \Big\langle \Big\langle 
   [\sqrt{\sigma^\delta(v^{\varepsilon,+})}  
-\sqrt{\sigma^\delta(u^\varepsilon)}]
   \frac{ D(v^{\varepsilon,+})} {\sqrt{\sigma(v^{\varepsilon,+})}}
 v^{\varepsilon,+}_x,
    l'(z^{\varepsilon,+}) \varphi_x \Big\rangle\Big\rangle
\qquad\quad
\end{eqnarray}
For $l$ convex and $\varphi$ nonnegative, the first term in the second
lines of \eqref{e:lfuncext1} and \eqref{e:lfuncext2} is
nonpositive. With these assumptions on $l$ and $\varphi$ we thus
define $B_l \equiv B_{l,\varphi}^{\varepsilon,\delta,N,L,\pm}:=\big[
\langle \langle D(u^\varepsilon) z^{\varepsilon,\pm}_x,
l''(z^{\varepsilon,\pm}) z^{\varepsilon,\pm}_x \varphi \rangle \rangle
\big]^{1/2}$ and let, for $F \in C([0,1])$,
\begin{equation*}
C_{F,l}^\delta:= \max \big\{ l''(z) |F(v+z)-F(v)|^2\,:\:v\in
[\delta,1-\delta],\,z\in[-1,1],\,v+z\in [0,1] \big\}
\end{equation*}
Since $v^{\varepsilon,\pm}_x,\,z^{\varepsilon,\pm}_x \in
L_{2,\mathrm{loc}}\big(P^{N,L,\pm}\big)$, by \eqref{e:zbound}, Cauchy-Schwarz
inequality and the fact that $D$ is uniformly positive, we have for each
nonnegative $\varphi^\pm \in C^\infty_{\mathrm{c}}\big(P^{N,L,\pm}\big)$, and
for some constant $C\equiv C^{\varepsilon,\delta,N,L}_{\varphi^\pm}$ independent
of $l$
\begin{eqnarray*}
&&
\big| 
   \langle \langle f(u^\varepsilon)-f(v^{\varepsilon,\pm}), l''(z^{\varepsilon,\pm})
           z^{\varepsilon,\pm}_x \varphi^\pm \rangle \rangle
\big|
\\
&& 
\qquad \phantom{\le}
+\big|
\frac{\varepsilon}{2} \langle \langle
\big[D(u^\varepsilon)-D(v^{\varepsilon,\pm})\big] v^{\varepsilon,\pm}_x,
l''(z^{\varepsilon,\pm}) z^{\varepsilon,\pm}_x \varphi^\pm \rangle \rangle
\big|
\\
&& 
\qquad \phantom{\le}
+ \big|
\varepsilon \langle \langle 
  \big[ \sqrt{\sigma^\delta(v^{\varepsilon,+})} 
        -\sqrt{\sigma^\delta(u^\varepsilon)} \big]
  \frac{ D(v^{\varepsilon,+})} {\sqrt{\sigma(v^{\varepsilon,+})}}
 v^{\varepsilon,+}_x, 
    \varphi^+ l''(z^{\varepsilon,+})z^{\varepsilon,+}_x \rangle \rangle
\big|
\\
& &\qquad
\le C\,\big[\sqrt{C_{f,l}^\delta}
    +  \sqrt{C_{D,l}^\delta}
+\sqrt{C_{\sqrt{\sigma^\delta},l}^\delta} \big]
           B_l
\end{eqnarray*}
We also let $C_l:=\max_{z\in [-1,1]} |l'(z)|$ and note that, in view
of \eqref{e:vest} and \eqref{e:zbound}, for any nonnegative
$\varphi^\pm \in C^\infty_{\mathrm{c}}\big(P^{N,L,\pm}\big)$ and for
some constant $\tilde{C}=\tilde{C}^{\delta,N,L}_{\varphi^\pm}$
independent of $\varepsilon$ and $l$
\begin{eqnarray*} 
&&
    \frac{\varepsilon}{2} 
 \Big| \big\langle\big\langle D(u^\varepsilon) z^{\varepsilon,\pm}_x,
     l'(z^{\varepsilon,\pm}) \varphi^\pm_x \big\rangle \big\rangle 
\Big|
+ \frac{\varepsilon}{2} 
\Big|\big\langle\big\langle
    \big[D(u^\varepsilon)-D(v^{\varepsilon,\pm})\big] v^{\varepsilon,\pm}_x, 
    l'(z^{\varepsilon,\pm}) \varphi^\pm_x \big\rangle \big\rangle 
\Big|
\\
&&\qquad \phantom{\le}
+ \varepsilon \Big|  \Big\langle\Big\langle 
     [\sqrt{\sigma(v^{\varepsilon,+})}  -\sqrt{\sigma(u^\varepsilon)}]
  \frac{ D(v^{\varepsilon,+})} {\sqrt{\sigma(v^{\varepsilon,+})}}
 v^{\varepsilon,+}_x, 
    l'(z^{\varepsilon,+}) \varphi^+_x \Big\rangle\Big\rangle 
\Big|
\\
&&\qquad
\le \tilde{C} \, C_l \, \sqrt{\varepsilon}
\end{eqnarray*}

Patching all together, for each nonnegative $\varphi^\pm \in
C^\infty_{\mathrm{c}}\big(P^{N,L,\pm}\big)$ we gather
\begin{eqnarray}
\label{e:zext1}
\nonumber
&& 
-\big\langle\big\langle l(z^{\varepsilon,\pm}),\varphi^\pm_t
\big\rangle\big\rangle
- \big\langle\big\langle f(u^\varepsilon)-f(v^{\varepsilon,\pm}),
          l'(z^{\varepsilon,\pm}) \varphi^\pm_x 
\big\rangle \big\rangle
\\
\nonumber
&&
\qquad \quad
\le -\frac{\varepsilon}{2} B_l^2
  +C\,\Big[\sqrt{C_{f,l}^\delta}+ \sqrt{C_{D,l}^\delta} 
+ \sqrt{C_{\sqrt{\sigma^\delta},l}^\delta} \Big] 
   B_l
+ \tilde{C} C_l \sqrt{\varepsilon}
\qquad \quad
\\
&&
\qquad \quad
\le \frac{3}{2\varepsilon} C^2 \big[
  C_{f,l}^\delta + C_{D,l}^\delta
  +  C_{\sqrt{\sigma^\delta},l}^\delta
     \big]
+ \tilde{C} C_l \sqrt{\varepsilon}
\end{eqnarray}
It is then easily seen that we can take a sequence of convex smooth
functions $\{l_n\} \subset C^2([-1,1])$ such that $|l'_n(z)|\le 1$,
$l_n(z)\to |z|$, $z l_n'(z) \to |z|$ uniformly on $[-1,1]$, and such
that, by the H\"older continuity hypotheses on $D$ and $\sigma$
\begin{equation*}
  \lim_{n \to \infty} \big( C_{f,l_n}^\delta + C_{D,l_n}^\delta
  +  C_{\sqrt{\sigma^\delta},l_n}^\delta \big)=0
\end{equation*}
Evaluating \eqref{e:zext1} for $l\equiv l_n$, taking the limit $n \to
\infty$, and recalling that we assumed $f'$ to be positive on $[0,1]$,
we gather for each nonnegative $\varphi^\pm \in
C^\infty_{\mathrm{c}}\big(P^{N,L,\pm}\big)$
\begin{equation}
\label{e:zext2}
-\langle \langle |u^\varepsilon -v^{\varepsilon,\pm}|,
                            \varphi^\pm_t\rangle \rangle
 -  \langle \langle \big|f(u^\varepsilon)-f(v^{\varepsilon,\pm})\big|,
                                        \varphi^\pm_x \rangle \rangle
\le \tilde{C} \sqrt{\varepsilon}
\end{equation}
We now reintroduce the dropped indexes $\delta,N,L$, and recall that
for $\delta \le \delta_0$, $\varepsilon \le \varepsilon_0(\delta_0)$,
$N \in \bb N$ and $L \ge 1$ we have $u^{\varepsilon,\delta,N,L} \in
\mc K_0$ for some compact $\mc K_0 \subset \mc X$. Let $u^{N,L} \in
\mc K_0$ be a generic limit point of $\{u^{\varepsilon,\delta,N,L}\}$
in $\mc X$ as $\varepsilon \to 0$ and successively $\delta \to 0$. By
\eqref{e:vconv} and \eqref{e:zext2}, for each nonnegative $\varphi \in
C^{\infty}_{\mathrm{c}}\big(P^{N,L,-}\cup P^{N,L,+}\big)$
\begin{equation}
\label{e:zext3}
    - \big\langle\big\langle |u^{N,L}-\tilde{u}|, 
                    \varphi_t \big\rangle\big\rangle
    - \big\langle\big\langle \big|f(u^{N,L}) -f(\tilde{u})\big|,
                    \varphi_x \big\rangle\big\rangle 
 \le 0
\end{equation}
Since $u^{N,L} \in \mc K_0$, there exist $u^L \in \mc X$ and a
subsequence $\{N_k\} \subset \bb N$ such that $u^{N_k,L} \to u^L$ in
$\mc X$ as $k \to +\infty$. By \eqref{e:xilarge} and \eqref{e:zext3},
it follows that for each nonnegative $\varphi \in
C^{\infty}_{\mathrm{c}}\big(R^{L,-} \cup R^{L,+} \big)$ 
\begin{equation}
\label{e:zext4}
   - \big\langle \big\langle |u^{L}-\tilde{u}|, 
                    \varphi_t \big\rangle \big\rangle
    - \big\langle\big\langle \big|f(u^{L}) -f(\tilde{u})\big|,
                    \varphi_x \big\rangle \big\rangle 
 \le 0
\end{equation}
Since $\tau^L$ is dense in $[0,T]$, by \eqref{e:rlarge} and
\eqref{e:nil} we have that, for $L \ge 1$, $R^{L,+} \cup R^{L,-}$ is
dense in $[0,T] \times \big[-L+\frac{1}{4L},L-\frac{1}{4L}\big]$.
Note also that $\tilde{u} \in \mc S_\sigma \subset \mc E$ by
hypotheses. Furthermore, since $u^{L}$ is a limit point of a sequence
with uniformly bounded $H_\varepsilon$-cost, we also have $u^{L} \in
\mc E$ by item (ii) in Theorem~\ref{t:ecne}, namely $\tilde{u}$ and
$u^L$ are entropy-measure solutions to \eqref{e:1.1}. By
Lemma~\ref{l:l1}, $\tilde{u},\,u^{L} \in
C\big([0,T];L_{1,\mathrm{loc}}(\bb R)\big)$. By the same
Lemma~\ref{l:l1} and the assumption $V_f^->0$ we have that the maps $x
\mapsto \tilde{u}(t,x)$ and $x \mapsto u^L(t,x)$ are continuous from
$\bb R$ to $L_1\big([0,T]\big)$.  Therefore, since the boundaries of
$R^{L,+}$ and $R^{L,-}\setminus R^{L,+}$ are countable unions of
segments parallel to the $x$ and $t$ axes, we have that
\eqref{e:zext4} holds for each nonnegative $\varphi \in
C^{\infty}_{\mathrm{c}} \big((0,T)\times
(-L+\frac{1}{4L},L-\frac{1}{4L})\big)$.

Recalling $\{u^L\} \subset \mc K_0$, let $u$ be a limit point of
$\{u^L\}$ along a subsequence $L_k \to \infty$. From \eqref{e:zext4}
we get for each nonnegative $\varphi \in C^\infty_{\mathrm{c}} \big(
(0,T)\times \bb R\big)$
\begin{equation}
\label{e:zext5}
 - \langle \langle |u-\tilde{u}|,
                      \varphi_t
   \rangle \rangle
- \langle \langle
    \big|f(u) -f(\tilde{u})\big|, \varphi_x
\rangle \rangle
\le 0
\end{equation}
Reasoning as above, we also have $u \in \mc E$, and thus setting
$z:=u-\tilde{u}$, by Lemma~\ref{l:l1}, $u,\,\tilde{u},\,z \in
C\big([0,T];L_{1,\mathrm{loc}}(\bb R) \big)$. By \eqref{e:zext5}, it
is then easily seen that for each bounded nonnegative Lipschitz
function $\varphi$ on $[0,T]\times \bb R$ such that
$\int\!dt\,dx\,[|\varphi|+ |\varphi_t|+|\varphi_x|]<+\infty$, and for
each $t \in [0,T]$
\begin{eqnarray}
\label{e:zext6}
\nonumber
&& 
\!\!
\langle |z(t)|, \varphi(t) \rangle - \langle |z(0)|, \varphi(0) \rangle
\\
&&
\;
- \int_{[0,t]}\!dr 
\big[ \langle |z|(r), \varphi_r(r) \rangle
+ \langle \big|f(\tilde{u}(r)+z(r)) -f(\tilde{u}(r))\big|, \varphi_x(r)
\rangle \big]
\le 0 \qquad\;\;\;
\end{eqnarray}
Fixed $L \ge 1$, we evaluate the inequality \eqref{e:zext6} for
$\varphi(t,x) \equiv \varphi^{L}(x)$ defined as
\begin{equation*}
\varphi^L(x):=
\begin{cases}
e^{-(L-x)} & \text{if $x<-L$}
\\
1 & \text{if $-L\le x \le L$}
\\
e^{-(x-L)} & \text{if $x>L$}
\end{cases}
\end{equation*}
so that setting $Z^L(t):= \langle |z(t)|,\varphi^L\rangle$ we have
\begin{equation*}
Z^L(t)-Z^L(0) \le V_f^+ \int_{[0,t]}\! dr\, 
\langle |z|(r), |\varphi^L_x|\rangle \le V_f^+ \int_{[0,t]}\!dr\,Z^L(r)
\end{equation*}
By Gronwall inequality, for each $L \ge 1$ and each $t \in [0,T]$, 
we have $Z^L(r) \le \exp[V_f^+\,t] Z^L(0)$. 
Note that $u(0,x)=\tilde{u}(0,x)$ by
\eqref{e:udl} and the definition of convergence in $\mc X$. Therefore
$Z^L(0)=0$, and thus $Z^L(t)=0$ for each $t \in [0,T]$ and
$L \ge 1$. Hence $u=\tilde{u}$.
\qed

\smallskip
\noindent
\textbf{Proof of Proposition~\ref{p:H}.}
In order to show that $H$ is lower semicontinuous, first note that the
set of weak solutions is closed in $\mc X$. Moreover for each entropy
sampler $\vartheta$ the map $\mc X \ni u \mapsto P_{\vartheta,u} \in \bb
R$ is continuous. On the other hand, if $u$ is a weak solution to
\eqref{e:1.1} then the equalities in \eqref{e:Hsup} holds, and thus
$H$ is a supremum of continuous maps.

Since $D(\cdot)/\sigma(\cdot)$ is uniformly positive on $[0,1]$,
$H(u)=0$ iff $u \in \mc E$ and $\varrho_u^+=0$, thus $u$ is entropic.
Conversely, entropic solutions $u$ are in $\mc E$ by item (i) in
Proposition~\ref{p:kin}, and the entropic condition is thus equivalent
to $\varrho_u^+=0$.

The coercivity of $H$ follows from the Tartar's method of compensated
compactness, that we already applied in the proof of
Theorem~\ref{t:ecne} item (ii). Suppose indeed that we are given a
sequence $\{u^n\} \subset \mc X$ such that $H(u^n) \le C_H <+\infty$
for each $n$. Then each $u^n$ is an entropy-measure solution to
\eqref{e:1.1} by the definition of $H$. For each entropy $\eta$, each
$n,\,L>0$, by the same bound in the proof of Proposition~\ref{p:kin},
$\|\wp_{\eta,u^n}\|_{{\mathrm{TV},L}} \le 2
\|\wp_{\eta,u^n}^+\|_{\mathrm{TV},L} +
2\,(\|\eta\|_{\infty}+\|q\|_{\infty})(2L+T)$. On the other hand, for
each $\eta \in C^2([0,1])$ such that $\sigma \eta'' \le D$,
$\|\wp_{\eta,u^n}^+\|_{\mathrm{TV},L} \le H(u^n)$ and therefore
$\|\wp_{\eta,u^n}\|_{{\mathrm{TV},L}}$ is bounded uniformly in $n$.
Since $\eta$ and $q$ are bounded, we have that
$\{\eta(u^n)_t+q(u^n)_x\}$ is precompact in
$H^{-1}_{\mathrm{loc}}([0,T]\times \bb R)$. As we already noted in the
proof of Theorem~\ref{t:ecne} item (ii), see \cite[Ch. 9]{Se}, this
yields the compactness of $\{u^n\}$ in $\mc X$.  \qed

\smallskip
\noindent
\textbf{Proof of Remark~\ref{r:bv}.}
  By well known properties of functions of locally bounded variation,
  for each entropy $\eta$ and $u \in \mc X \cap
  BV_{\mathrm{loc}}([0,T]\times \bb R)$ we have that $\wp_{\eta,u}$ is
  a Radon measure on $(0,T)\times \bb R$. If $u$ is a weak solution to
  \eqref{e:1.1}, by Vol'pert chain rule \cite{AFP}, the absolutely
  continuous and Cantor parts of $\wp_{\eta,u}$ w.r.t.\ the Lebesgue
  measure on $(0,T)\times \bb R$ vanish, and we get
\begin{equation*}
d\wp_{\eta,u}
          = \Big\{ \big[\eta(u^+)-\eta(u^-)\big]n^t
           +\big[q(u^+)-q(u^-)\big]n^x \Big\} \, d \mc H^{1} \res J_u
\end{equation*}
On the other hand the Rankine-Hugoniot condition $\big[u^+ - u^-
\big]n^t+\big[f(u^+)-f(u^-)\big]n^x=0$ holds. The statement of the
remark follows by direct calculations.
\qed

\smallskip
\noindent
\textbf{Proof of Remark~\ref{r:JV}.}
For $u \in \mc E$ we have
\begin{eqnarray*}
  H^\prime(u) &=& \sup  \big\{
  \wp_{\eta,u}(\varphi) \, , \: 
  \varphi \in C^\infty_{\mathrm{c}}\big((0,T)\times \bb R;[0,1] \big),
\\
&&\phantom{ \sup  \big\{  \wp_{\eta,u}(\varphi)\, , \:}
 \eta \in C^2([0,1])\,:\:0\le \sigma \eta'' \le D \big\}
\end{eqnarray*}
so that the inequality $H \ge H^\prime$ follows from the equalities in
\eqref{e:Hsup}. The same inequality yields $H(u)=H^\prime(u)$ if there
exists a set $E^+$ as in the statement of the remark. If $f$ is convex
or concave and $u$ has locally bounded variation, we can take
$E^+=\{(t,x) \in J_u\,:\:\exists v \in [0,1]\,:\:
\rho(v,u^+,u^-)>0\}$, where $J_u$, $u^\pm$ and $\rho$ are defined as
in Remark~\ref{r:bv}.

If $f$ is neither convex nor concave, then there exist
$u^-,\,u^+,\,v^\prime,\,v^{\prime \prime} \in (0,1)$ such that
$\rho(v^\prime,u^+,u^-)>0$ and $\rho(v^{\prime \prime},u^+,u^-)<0$,
where $\rho$ is defined as in Remark~\ref{r:bv}. Let
$V:=\frac{f(u^+)-f(u^-)}{u^+-u^-}$, and define $u:[0,T]\times \bb R
\to [0,1]$ by
\begin{equation*}
u(t,x):=
\begin{cases}
u^+ & \text{for $x< V\,t$}
\\
u^- & \text{for $x>V\,t$}
\end{cases}
\end{equation*}
Then $u \in \mc E$ and by a direct computation $H(u) > H^\prime(u)$.
\qed

\appendix
\section{$\mc I$-approximation of atomic Young measures}
\label{s:A}

Here we prove the claims stated in the proof of Theorem~\ref{t:lsce},
\textit{Step 1}, where the sets $\mc M^1_n$, $\widetilde{\mc M}^1_n$,
$\upbar{\mc M}^1_n$ are defined.

\smallskip
\noindent \textit{Claim 1: $\widetilde{\mc M}_1^n$ is $\mc I$-dense in
  $\upbar{\mc M}_1^n$.} For $n\ge 1$, let $\mu \in \upbar{\mc
  M}_1^n$, let $G^\mu$ be defined as in
Lemma~\ref{l:riesz2}. Let also $r,\,\alpha^i,\,u^i$ be as in the
definition of $\upbar{\mc M}_1^n$ and $L,\,\mu_\infty$ be as in the
definition \eqref{e:mg} of $\mc M_g$. With no loss of generality, we
can assume that $u^{i+1}\ge u^i$, $i=1,\ldots,n-1$, since we can
reorder the $u^i(t,x)$ for all $(t,x)$ preserving continuity of the
$u^i$ and measurability of the $\alpha^i$. Analogously it is not
restrictive to assume, for $|x|>L$, $u^i(t,x)=u^i_\infty$,
$\alpha^i(t,x)=\alpha^i_\infty$ for some constants
$u^i_\infty,\,\alpha^i_\infty \in (0,1]$; in particular
$\mu_\infty=\sum_i \alpha^i_\infty \delta_{u^i_\infty}$.

Let now $\{\jmath^k\} \subset C^\infty_{\mathrm{c}}(\bb R \times \bb
R)$ be a sequence of smooth mollifiers supported by $[-T/k,T/k]\times
\bb [-1,1]$, and recall the definition \eqref{e:b} of $b^k$.  For
$i=1,\ldots,n$ and $h,k \ge 1$ define $\alpha^{i;k} \in C^1
\big([0,T]\times \bb R;[r,1] \big)$, and $u^{i;h,k} \in
C^1\big([0,T]\times \bb R;[h^{-1},1-h^{-1}]\big)$ by
\begin{eqnarray}
\label{e:alphaik}
\nonumber
&&
\alpha^{i;k}(t,x)  :=  
      \int\!dy\,ds\,
        \jmath^k(t-s,x-y) \alpha^i\big(b^k(s,y)\big)
\\
\nonumber
&& u^{i;h,k}(t,x) := h^{-1} \Big[1 +
                  \frac{i}{n \sum_{i^\prime} i^{\prime} 
                                   \alpha^{i^\prime;k}(t,x)} \Big]
\\
&& \qquad\quad
+\frac{1-3\,h^{-1}}{\alpha^{i;k}(t,x)} 
              \int\!dy\, ds\, \jmath^k(t-s,x-y)\,
              \alpha^i\big(b^k(s,y)\big)\,u^i\big(b^k(s,y)\big)
\qquad\quad
\end{eqnarray}
Clearly $\alpha^{i;k}$ and $u^{i;h,k}$ are smooth, with $\alpha^{i;k}
\ge r$, $\sum_i \alpha^{i;k}=1$, and $\alpha^{i;k}$, $u^{i;h,k}$ are
constant for $|x|>L+1$. Furthermore for $i=1,\ldots,n-1$ and $(t,x)
\in [0,T] \times \bb R$
\begin{eqnarray*}
& &
  \lim_{k\to \infty}
  \big[u^{i+1;h,k}(t,x) - u^{i;h,k}(t,x) -\frac{h^{-1}}{n^2}\big]
\\
&& \qquad\quad
\ge \lim_{k \to \infty} 
 \Big[u^{i+1;h,k}(t,x) - u^{i;h,k}(t,x) 
        -\frac{h^{-1}}{n \sum_{i^\prime} i^{\prime} 
         \alpha^{i^\prime;k}(t,x) }\Big]
\\ &&
\qquad\quad
     = \big[1-3h^{-1}\big] \big[u^{i+1}(t,x)-u^{i}(t,x) \big]
\end{eqnarray*}
Since the $u^i$ are continuous, it is not difficult to see that
convergence in the last line above is uniform on compact subsets of
$[0,T]\times \bb R$. On the other hand, since the $u^{i}$ and
$u^{i;h,k}$ are constant for $|x|>L+1$, we have that convergence is
indeed uniform on $[0,T]\times \bb R$. It follows that for each $h>1$
there exists $K^h\ge 1$ such that $u^{i+1;h,k} \ge u^{i;h,k} +
h^{-1} n^{-2}/2$ for each $k \ge K^h$. Therefore, defining $\mu^{h,k}
\in \mc M$ by
\begin{equation*}
\mu^{h,k}_{t,x}:=\sum_{i=1}^n \alpha^{i;k}(t,x)
                      \delta_{u^{i;h,k}(t,x)}
\end{equation*}
we get, for $k \ge K^h$, $\mu^{h,k} \in \widetilde{\mc M}_1^n$
provided $\mc I(\mu^{h,k})<+\infty$. Recalling Lem\-ma~\ref{l:riesz2},
this follows by the existence of $G^{\mu^{h,k}} \in L_2([0,T]\times
\bb R)$ satisfying weakly on $(0,T) \times \bb R$:
\begin{equation*}
\mu^{h,k}(\imath)_t+\mu^{h,k}(f)_x =-G^{\mu^{h,k}}_x
\end{equation*}
Indeed $G^{\mu^{h,k}}$ can be computed
explicitly as
\begin{eqnarray*}
G^{\mu^{h,k}}(t,x) & := & \big(1-3h^{-1}\big)
                        \int \!ds\, dy\,
                          \jmath^k(t-s,x-y) G^\mu \big(b^k(s,y)\big)
\\ & &
 +  \big(1-3h^{-1}\big) \int \!ds \, dy\,
 \jmath^k(x-y,t-s)\,\mu_{b^k(s,y)}(f)
\\ & &
 -  \mu^{h,k}(f)
 -(1-3h^{-1}) \mu_\infty(f) 
 + \mu_\infty^{h,k}(f)
\end{eqnarray*}
where
\begin{equation*}
\mu_\infty^{h,k}(f):= \sum_i \alpha^i_\infty f\Big(h^{-1} 
         +\frac{i\,h^{-1}}{n \sum_{i^\prime} i^{\prime} 
                                 \alpha^{i^\prime}_\infty} 
        +\big(1-3h^{-1}\big)u^i_\infty  \Big)
\end{equation*}
It immediately follows that $\lim_{h \to \infty} \lim_{k \to \infty}
\|G^{\mu^{h,k}}-G^{\mu}\|_{L_2([0,T]\times \bb R)}=0$, and it is also
straightforward to see that, for each $F \in C([0,1])$
\begin{equation*}
  \lim_{h\to \infty} \lim_{k\to \infty} \mu^{h,k}(F)=\mu(F) 
  \qquad
  \text{strongly in $L_{1,\mathrm{loc}}([0,T]\times \bb R)$}
\end{equation*}
By Remark~\ref{l:iconv}, we can extract a subsequence $\{\mu^k\}$ from
$\{\mu ^{h,k}\}$ that $\mc I$-converges to $\mu$.

\smallskip
\noindent \textit{Claim 2: $\upbar{\mc M}_1^n$ is $\mc I$-dense in
  $\mc M_1^n$.} For $n \ge 1$, let $\mu \in \mc M_1^n$. Let also
$\alpha^i,\,u^i$ and $L$ be as in the definition of $\mc M_1^n$ and
$\mc M_g$. With no loss of generality, we can assume that
$\alpha^i>0$, since we do not require the $u^i$ to be distinct. As in
\textit{Claim 1} above, we can also assume that, for $|x|>L$,
$u^i(t,x)=u^i_\infty$, $\alpha^i(t,x)=\alpha^i_\infty$ for some
constants $u^i_\infty,\,\alpha^i_\infty \in [0,1]$.

With these assumptions, for $h,\,k \ge 1$ and $i=1,\ldots,n$, let us
define $\alpha^{i;k}$ as in \eqref{e:alphaik}, and 
$u^{i;k}$ by
\begin{equation*}
\begin{array}{lcl}
u^{i;k}(t,x)  := \frac{1}{\alpha^{i;k}} 
    \int\!dy\, ds\,
      \jmath^k(t-s,x-y) \alpha^i\big(b^k(s,y)\big) 
         u^i\big(b^k(s,y) \big)
\end{array}
\end{equation*}
Letting
\begin{equation*}
\mu^{k}_{t,x}:=\sum_{i=1}^n \alpha^{i;k}(t,x)
                      \delta_{u^{i;k}(t,x)}
\end{equation*}
we gather $\mu^k \in \upbar{\mc M}_1^n$. A computation similar to the
one carried out in \textit{Claim 1} shows that $\mu^k$ $\mc
I$-converges to $\mu$ as $k\to \infty$.

\smallskip
\noindent \textit{Claim 3: $\mc M_1^n$ is $\mc I$-dense in
  $\widetilde{\mc M}_1^{n+1}$.} This is the key step in the proof of
Theorem~\ref{t:lsce}. For $n \ge 1$, let $\mu \in \widetilde{\mc
  M}_1^{n+1}$, and let $G^\mu$ be defined as in Lemma~\ref{l:riesz2}.
Let also $r,\,\alpha^i,\,u^i$ be as in the definition of
$\widetilde{\mc M}_1^{n+1}$, and $L,\,\mu_\infty$ as in the definition
\eqref{e:mg} of $\mc M_g$. Note that for $|x|>L$,
$\alpha^{i}(t,x)=\alpha^{i}_\infty$ and $u^i(t,x)=u^i_\infty$ for some
constants $\alpha^i_\infty \in [r,1-r]$, $u^i_\infty \in [r,1-r]$,
with $u^{i+1}_\infty \ge u^i_\infty + r$, $i=1,\ldots,n$.

Let us define the Young measures $\nu^1,\,\nu^0 \in \widetilde{\mc
M}_1^n$ by
\begin{equation*}
  \nu^1_{t,x}:=\delta_{u^{n+1}(t,x)} \qquad 
  \nu^0_{t,x}:=\sum_{i=1}^n \frac{\alpha^i(t,x)}{1-\alpha^{n+1}(t,x)} 
                                 \delta_{u^i(t,x)}
\end{equation*}
so that, letting $\beta(t,x):= \alpha^{n+1}(t,x) \le 1-r$
\begin{equation*}
\mu_{t,x}=\beta(t,x)\nu^1_{t,x}+\big(1-\beta(t,x)\big) \nu^0_{t,x}
\end{equation*}
The basic idea is to build up a sequence $\{\mu^k\}$ $\mc I$-converging to
$\mu$, as follows: we first slice up $[0,T]\times \bb R$ in small
strips, alternating a strip of width $\beta\,k^{-1}$ with a strip of
width $(1-\beta) k^{-1}$; we then set $\mu^k_{t,x}=\nu^1_{t,x}$ for
$(t,x)$ in the first family of strips, and $\mu^k_{t,x}=\nu^0_{t,x}$
for $(t,x)$ in the second family of strips. As we let $k\to \infty$,
we easily get $\mu^k \to \mu$; however, to get also $\mc I(\mu^k) \to
\mc I(\mu)$, we will have to carefully define these strips.

For $j \in \bb Z$ and $k \in \bb N$, let us consider the maps
$\gamma^k_j:[0,T]\to \bb R$ solutions to
\begin{equation}
\label{e:gamma}
\begin{cases}
{\displaystyle
  \dot{\gamma}=\frac{\nu^1_{t,\gamma}(f)-\nu^0_{t,\gamma}(f)}
  {\nu^1_{t,\gamma}(\imath)-\nu^0_{t,\gamma}(\imath)}
}
\\
{\displaystyle \vphantom{\big\{^{\Big\{}}
\gamma(0)=\frac{j}{k}
}
\end{cases}
\end{equation}
These equations are well-posed since $\nu^1(f)$, $\nu^0(f)$,
$\nu^1(\imath)$, $\nu^0(\imath)$ are Lipschitz functions in the
$(t,x)$ variables, and $\nu^1(\imath)-\nu^0(\imath) \ge r$, by the
definition of $\widetilde{\mc M}_1^{n+1}$. Furthermore, by standard
theory for \eqref{e:gamma}, $\gamma^k_j \in C^0([0,T])\cap
C^1((0,T))$; $|\dot{\gamma}^k_j| \le 2 r^{-1} \max_{v \in
  [0,1]}|f(v)|$; $\gamma^k_{j+1}>\gamma^k_j$; and
$\gamma^k_{j+1}(t)-\gamma^k_{j}(t) \le C k^{-1}$ for some constant $C$
independent of $k,\,j$ and $t$.

We next define the maps $\beta^k_j:[0,T] \to \bb R$ by
\begin{equation}
\label{e:beta}
\int_{\gamma^k_j(t)}^{\gamma^k_j(t)+\beta^k_j(t)} \!dx\,
             \big[ \nu^1_{t,x}(\imath) - \nu^0_{t,x}(\imath)\big] = 
\int_{\gamma^k_j(t)}^{\gamma^k_{j+1}(t)} \!dx\, \beta(t,x) 
             \big[ \nu^1_{t,x}(\imath) - \nu^0_{t,x}(\imath)\big]
\end{equation}
Since $\nu^1_{t,x}(\imath) - \nu^0_{t,x}(\imath) \ge r>0$, for any
fixed $t \in [0,T]$ the l.h.s.\ of this equation is strictly increasing in
$\beta^k_j(t)$. Since it vanishes for $\beta^k_j(t)=0$ and it is
larger than the r.h.s.\ for
$\beta^k_j(t)=\gamma^k_{j+1}(t)-\gamma^k_{j}(t)$ (recall $\beta(t,x)
\in [r, 1-r]$), there exists a unique
$0<\beta^k_j(t)<\gamma^k_{j+1}(t)-\gamma^k_j(t)$ satisfying
\eqref{e:beta}. Furthermore, since $\beta$ and the $\gamma^k_j$ are
smooth, we have $\beta^k_j \in C^0([0,T])\cap C^1((0,T))$. The mean
value theorem then implies
\begin{equation}
\label{e:alphaconv}
\Big|\beta^k_j(t)
   - \int_{\gamma^k_j(t)}^{\gamma^k_{j+1}(t)}\!dx\, \beta(t,x)\Big| 
\le C \big[\gamma^k_{j+1}(t)-\gamma^k_j(t) \big]^2
\le C^\prime k^{-2}
\end{equation}
for suitable constants $C,\,C^\prime$. For $h$ and $k$ two positive
integers, we next define the Young measure $\mu^{h,k} \in \mc M$ by
\begin{equation*}
\mu^{h,k}_{t,x}:=
\begin{cases}
\nu^0_{t,x} & \text{if $\exists j \in \bb Z$, $|j|\le h\,k$ such
            that $\gamma^k_j(t) +\beta^k_j(t)<x < \gamma^k_{j+1}(t)$}
\\
\nu^1_{t,x} & \text{otherwise}
\end{cases}
\end{equation*}
Since $\nu^1_{t,x}$ is constant for $|x|$ sufficiently large, we have
$\mu^{h,k} \in \mc M_g$ for $h$ large enough. Furthermore, since
convergence in $\mc M$ is \emph{local}, \eqref{e:alphaconv} yields
$\lim_{h\to \infty} \lim_{k\to \infty} \mu^{h,k}=\mu$ in $\mc M$, and
for each $F \in C([0,1])$
\begin{equation*}
  \lim_{h\to \infty} \lim_{k\to \infty} \mu^{h,k}(F)=\mu(F) 
   \qquad \text{strongly in $L_{1,\mathrm{loc}}([0,T]\times \bb R)$}
\end{equation*}
We next prove that $\mc I(\mu^{h,k}) <+\infty$ and $\lim_h \lim_k
G^{\mu^{h,k}}=G$ in $L_2([0,T]\times \bb R)$; so that, reasoning as in
the proof of \textit{Claim 1}, by Remark~\ref{l:iconv} we get the
existence of a subsequence $\{\mu^k\}$ $\mc I$-converging to $\mu$.
For each $F \in C([0,1])$, $(t,x) \mapsto \mu_{t,x}^{h,k}(F)$ is
smooth outside the graph of the curves $\gamma^k_j$. Therefore by
Lemma~\ref{l:ruho} there exists $G^{h,k} \in
L_{2,\mathrm{loc}}([0,T]\times \bb R)$, such that $\mu^{h,k}(\imath)_t
+\mu^{h,k}(f)_x = -G^{h,k}_x$ holds weakly. First we show that we can
choose $G^{h,k}$ to be compactly supported, so that $G^{h,k}\in
L_2([0,T]\times \bb R)$, and thus $\mc I(\mu^{h,k})<+\infty$ with
$G^{h,k}=G^{\mu^{h,k}}$ according to the definition given in
Lemma~\ref{l:riesz2}.

Since $G^{h,k}$ is defined up to a measurable function of $t$,
and $G^{h,k}_x(t,x)=0$ for $x<\gamma^k_{-hk}(t)$ (we are
considering $h$ large enough as above), we can assume
$G^{h,k}(t,x)=G^{\mu}(t,x)=0$ for $x<\gamma^k_{-hk}(t)$.
Furthermore, by \eqref{e:ruho} and \eqref{e:gamma}, for each $j \in
\bb Z$, $G^{h,k}$ is continuous in the regions $\{(t,x)\,:\:
\gamma^k_j(t)+\beta^k_j(t)<x<\gamma^k_{j+1}(t)+\beta^k_{j+1}(t)\}$.
Let now $j \in \bb Z$ with $|j|\le hk$, and $t\in[0,T]$; by
\eqref{e:ruho} and \eqref{e:gamma}
\begin{eqnarray*}
&&-\Big[G^{h,k}\big(t,[\gamma^k_j(t)+\beta^k_j(t)]^- \big) 
    - G^{h,k}\big(t,\gamma^k_j(t)\big)
     \Big]
\\
&&\qquad\quad
= \;
\nu^1_{t,\gamma^k_j(t)+\beta^k_j(t)}(f)-\nu^1_{t,\gamma^k_j(t)}(f)
 + \int_{\gamma^k_j(t)}^{\gamma^k_j(t)+\beta^k_j(t)}\!dx\,
   \big[\nu^1_{t,x}(\imath)\big]_t
\end{eqnarray*}
and 
\begin{eqnarray*}
&&
-\big[G^{h,k}\big(t,\gamma^k_{j+1}(t) \big) 
      -G^{h,k}\big(t,[\gamma^k_j(t)+\beta^k_j(t)]^-\big)
     \big]
\\
&&\qquad\quad
= \;
\nu^0_{t,\gamma^k_{j+1}(t)}(f)-\nu^0_{t,\gamma^k_j(t)+\beta^k_j(t)}(f)
\\
& &\qquad\quad
\phantom{ = \;}
 + \int_{\gamma^k_j(t)+\beta^k_j(t)}^{\gamma^k_{j+1}(t)}\!dx\,
     \big[\nu^0_{t,x}(\imath)\big]_t
+ \big[\nu^0_{t,\gamma^k_j(t)+\beta^k_j(t)}(f)
     - \nu^1_{t,\gamma^k_j(t)+\beta^k_j(t)}(f) \big]
\\ 
&&\qquad\quad
\phantom{ = \;} 
- \big[\nu^0_{t,\gamma^k_j(t)+\beta^k_j(t)}(\imath) 
     - \nu^1_{t,\gamma^k_j(t)+\beta^k_j(t)}(\imath) \big]
         \big[\dot{\gamma}^k_j(t)+\dot{\beta}^k_j(t)\big]
\end{eqnarray*}
By \eqref{e:gamma}, \eqref{e:beta} and simple algebraic manipulations
\begin{eqnarray*}
&&
G^{h,k}\big(t,\gamma^k_{j+1}(t) \big) 
      -G^{h,k}\big(t,\gamma^k_j(t) \big)
\\
&&\qquad\quad
= 
-\big[ \mu_{t,\gamma^k_{j+1}(t)}(f)-\mu_{t,\gamma^k_j(t)}(f) \big]
- \int_{\gamma^k_j(t)}^{\gamma^k_{j+1}(t)}\!dx\,
     \big[\mu_{t,x}(\imath)\big]_t
\\
&&
\qquad\quad
=  G^{\mu}\big(t,\gamma^k_{j+1}(t) \big) 
      -G^{\mu}\big(t,\gamma^k_j(t) \big)
\end{eqnarray*}
 Since $G^{h,k}(t,\gamma^k_{-hk}(t))=G^{\mu}(t,\gamma^k_{-hk}(t))=0$,
 we deduce that for each $j \in \bb Z$. 
 we have $G^{h,k}(t,\gamma^k_{j}(t))=G^{\mu}(t,\gamma^k_{j}(t))$.
 In particular, since
 $G^{\mu}(t,\gamma^k_{hk}(t))=0$ and $G^{h,k}_x(t,x)=G^{\mu}_x(t,x)=0$
 for $x>\gamma^k_{hk}(t)$, we have $G^{h,k}(t,x)=G^{\mu}(t,x)=0$ for
 $x>\gamma^k_{hk}(t)$ and $x<\gamma^k_{-hk}(t)$. That is, $G^{h,k}$
 and $G^{\mu}$ are compactly supported. Thus $\mc
 I(\mu^{h,k})<+\infty$ and $G^{h,k}=G^{\mu^{h,k}}$.

 Finally, by the definition of $G^{\mu}$ and $G^{\mu^{h,k}}$,
 recalling $G^{h,k}(t,\gamma^k_{j}(t))=G^{\mu}(t,\gamma^k_{j}(t))$ we
 have
\begin{eqnarray*}
&&
\!\!
\big\|G^{\mu^{h,k}}-G^{\mu}\big\|^2_{L_2([0,T]\times \bb R)}
= 
\sum_{j=-hk}^{hk}\! \int_{[0,T]} \!dt
\int_{\gamma^k_j(t)}^{\gamma^k_{j+1}(t)} \!\!\!\!\!\!\!\!dx\,
\big(G^{\mu^{h,k}}(t,x)-G^{\mu}(t,x)\big)^2
\\ &&
\;
=\sum_{j=-hk}^{hk} \int_{[0,T]} \!dt \:
\Bigg\{
\\ &&
\quad
    \int_{\gamma^k_j(t)}^{\gamma^k_{j}(t)+\beta^k_j(t)}
   \!\!\!\!\!\!\!\!dx
        \Big[\int_{\gamma^k_j(t)}^x \!dy
             [\nu^1_{t,y}(\imath)-\mu_{t,y}(\imath)]_t
           + [\nu^1_{t,y}(f)-\mu_{t,y}(f)]_y
        \Big]^2
\\
&&
\quad
   +\int_{\gamma^k_{j}(t)+\beta^k_j(t)}^{\gamma^k_{j+1}(t)} \!dx\,
        \Big[\int_{\gamma^k_{j+1}(t)}^x \!dy\,
             [\nu^0_{t,y}(\imath)-\mu_{t,y}(\imath)]_t
           + [\nu^0_{t,y}(f)-\mu_{t,y}(f)]_y
        \Big]^2
\Bigg\}
\end{eqnarray*}
Since all the integrands in the are last two lines of this formula are
bounded uniformly in $h$ and $k$, each term of the sum is bounded by
$C\,k^{-3}$ for some constant $C>0$. Therefore the sum itself is
bounded by $2\,C\,h\,k^{-2}$, and we get $\lim_{h \to \infty} \lim_{k
  \to \infty} \mc I(\mu^{h,k}) = \mc I(\mu)$.

\section{$\Gamma$-viscosity cost 
for scalar Hamilton-Jacobi equations}
\label{s:B}
In this appendix we establish a $\Gamma$-convergence result for a
sequence of functionals associated with the Hamilton-Jacobi equation
\eqref{e:1.1}
\begin{equation}
\label{e:HJ}
b_t+f(b_x)=0
\end{equation}
which is related to \eqref{e:1.1} via the transformation $u=b_x$. In
\eqref{e:HJ} we understand $(t,x) \in [0,T]\times \bb R$ and $b(t,x)
\in \bb R$.  As usual, we assume $f$ to be a Lipschitz function on
$[0,1]$, $D$ and $\sigma$ continuous functions on $[0,1]$, with $D$
uniformly positive and $\sigma$ strictly positive on $(0,1)$. We will
just sketch most of the proofs, since they are similar to the proofs
of the corresponding statements for \eqref{e:1.1}.

We introduce the equivalence $\sim$ on
$C\big([0,T];L_{2,\mathrm{loc}}(\bb R \big)$ by setting $b^1 \sim b^2$
iff $b^1-b^2$ is constant in $[0,T]\times \bb R$. We let $\mc B$ be
the set of functions $b \in C\big([0,T];L_{2,\mathrm{loc}}(\bb R
\big)/\sim$ such that $b_x \in \mc U$. The requirement $b_x \in \mc U$
is clearly compatible with $\sim$, so that $\mc B$ is well defined. We
equip $\mc B$ with the metric
\begin{equation}
\label{e:dB}
d_{\mc B}(b^1,b^2):= d_{\mc U}(b^1_x,b^2_x) 
   + \inf_{c \in \bb R} \sup_{t \in [0,T]} 
       \sum_{N=1}^\infty \frac{1}{2^N}
       \|b^1(t,\cdot)-b^2(t,\cdot)+c\|_{L_2\big([-N,N] \big)}
\end{equation}
Note that the second term in the r.h.s. of \eqref{e:dB} is the
projection of the $C\big([0,T];L_{2,\mathrm{loc}}(\bb
R)\big)$-distance w.r.t.\ the $\sim$ equivalence. $(\mc B,d_{\mc B})$
is a complete separable metric space.

For $b \in \mc B$ such that $b_{x x} \in
L_{2,\mathrm{loc}}([0,T]\times \bb R) $ and $\varepsilon>0$ we next
define the linear functional $a^b_\varepsilon$ on
$C^{\infty}_{\mathrm{c}}\big((0,T)\times \bb R\big)$ by
\begin{equation}
\label{e:a}
a^{b}_{\varepsilon}(\varphi):= -\langle \langle b,\varphi_t \rangle \rangle
      +\langle \langle f(b_x),\varphi \rangle \rangle
      -\frac{\varepsilon}{2} 
       \langle \langle D(b_x)b_{x x},\varphi \rangle \rangle
\end{equation}
and the functional $J_\varepsilon:\mc B \mapsto [0,+\infty]$ 
as follows. If $b_{x x} \in L_{2,\mathrm{loc}}([0,T]\times \bb R)$ we
set 
\begin{equation}
\label{e:J}
J_\varepsilon(b):=
 \sup_{\varphi\in C^\infty_{\mathrm{c}}((0,T)\times \bb R)}
  \Big[ 
   a_\varepsilon^b(\varphi)
   - \frac 12 \langle\langle \sigma(b_x)\varphi , \varphi \rangle\rangle \Big]
\end{equation}
letting $J_\varepsilon(b)=+\infty$ otherwise.
We want to study the $\Gamma$-convergence of $\{J_\varepsilon\}$. As shown
below, this problem is strictly related to the $\Gamma$-convergence of
$\{I_\varepsilon\}$ defined in \eqref{e:2.6}.

We introduce the set $\mc A:=\big\{(b,\mu) \in \mc B \times \mc M\,:\:
b_x = \mu(\imath) \big\}$ which we equip with the metric
\begin{equation}
\label{e:dA}
d_{\mc A}\big((b^1,\mu^1),(b^2,\mu^2)\big):= d_{\mc B}(b^1,b^2)+d_{\mc
  M}(\mu^1,\mu^2)
\end{equation}
We say that $(b,\mu) \in \mc A$ is a measure-valued solution to
\eqref{e:HJ} iff $b_t+\mu(f)=0$ weakly in $(0,T)\times \bb R$. We lift
$J_\varepsilon$ to a functional $\mc J_\varepsilon: \mc A \to [0,+\infty]$
by setting
\begin{equation}
\label{e:mcJe}
   \mc J_\varepsilon(b,\mu):=
   \begin{cases}
     J_\varepsilon(b) & \text{ if } \quad
     \mu_{t,x} = \delta_{b_x(t,x)} \\
     +\infty  & \text{ otherwise }
   \end{cases}
\end{equation}

\begin{theorem}
\label{t:HJ1} 
The sequence $\{\mc J_\varepsilon\}$ is equicoercive on $\mc A$ and
$\Gamma$-converges to
\begin{equation*}
  \mc J \big((b,\mu) \big) :=
 \sup_{\varphi \in C^\infty_{\mathrm{c}}((0,T) \times\bb R) } \Big\{ 
   - \langle\langle b, \varphi_t \rangle\rangle 
   + \langle\langle \mu (f), \varphi \rangle\rangle 
   - \frac 12\, \langle\langle \mu (\sigma) 
                  \varphi ,\varphi \rangle\rangle 
     \Big\}
\end{equation*}
\end{theorem}

Note that $\mc J\big((b,\mu)\big)=0$ iff $(b,\mu)$ is a measure-valued
solution to \eqref{e:HJ}.

On the set $\mc B$ we next introduce the metric $d_{\mc Y}$
\begin{equation*}
d_{\mc Y}(b^1,b^2):= d_{\mc X}(b^1_x,b^2_x) 
   + \inf_{c \in \bb R} \sup_{t \in [0,T]} 
       \sum_{N=1}^\infty \frac{1}{2^N}
       \|b^1(t)-b^2(t)+c\|_{L_2\big([-N,N] \big)}
\end{equation*}
and denote by $(\mc Y,d_{\mc Y})$ the complete separable metric space
consisting of the same set $\mc B$ equipped with the distance $d_{\mc
  Y}$. We say that $b \in \mc Y$ is a weak solution to \eqref{e:HJ}
iff $-\langle \langle b,\varphi_t \rangle \rangle +\langle \langle
f(b_x),\varphi \rangle \rangle=0$ for each $\varphi \in
C^\infty_{\mathrm{c}}\big((0,T)\times \bb R \big)$. We denote by $\mc
W \subset \mc Y$ the set of weak solutions to \eqref{e:HJ}. We rescale
the functional $J_\varepsilon$ defining $K_\varepsilon:\mc Y \to [0,+\infty]$ as
\begin{equation}
\label{e:Keps}
K_\varepsilon:= \varepsilon^{-1} J_\varepsilon
\end{equation}

\begin{theorem}
\label{t:HJ2}
Let $K_\varepsilon$ be the functional on $\mc W$ as defined in
\eqref{e:Keps} and \eqref{e:J}.
\begin{itemize}
\item[{\rm (i)}] The sequence of functionals $\{K_\varepsilon\}$
  satisfies the $\Gamma$-liminf inequality
\begin{equation*}
\big(\Gliminf_{\varepsilon\to 0} K_\varepsilon\big)(b) \ge 
\begin{cases}
H(b_x) & \text{if $b \in \mc W$}
\\
+\infty & \text{otherwise}
\end{cases}
\end{equation*}
\item[{\rm (ii)}] Assume there is no interval where $f$ is affine.
Then the sequence $\{K_\varepsilon\}$ is equicoercive on
$\mc Y$.
\item[{\rm (iii)}] Suppose furthermore $f \in C^2([0,1])$ and
  $D,\sigma \in C^{\alpha}([0,1])$ for some $\alpha>1/2$. Then
\begin{equation*}
\big(\Glimsup_{\varepsilon\to 0} K_\varepsilon\big)(b) \le 
\begin{cases}
\upbar{H}(b_x) & \text{if $b \in \mc W$}
\\
+\infty & \text{otherwise}
\end{cases}
\end{equation*}
\end{itemize}
\end{theorem}

Since $b(0,\cdot)$ is bounded and Lipschitz, by a well known
connection between entropic solutions to \eqref{e:1.1} and viscosity
solutions to \eqref{e:HJ}, see e.g.\ \cite[Theorem~1.1]{KR}, we gather
$(\Gliminf_\varepsilon K_\varepsilon)(b)=0$ iff $b$ is a viscosity solutions
to \eqref{e:HJ}. It follows that if $b^\varepsilon$ satisfies the
equation
\begin{equation}
\label{e:HJforced}
b_t+f(b_x)=\frac{\varepsilon}{2} D(b_x)b_{xx} - \sigma(b_x) E^\varepsilon
\end{equation}
for some $E^\varepsilon$
such that 
$\lim_\varepsilon \varepsilon
\|E^\varepsilon\|_{ L_2([0,T]\times \bb R, \sigma(b_x)dt\,dx)}^2
=0$, then limit points of $\{b^\varepsilon\}$ are viscosity solutions to
\eqref{e:HJ}. Analogously if $b^\varepsilon$ solves
\eqref{e:HJforced} for some $E^\varepsilon$ with $\varepsilon
\|E^\varepsilon\|_{ L_2([0,T]\times \bb R, \sigma(b_x)dt\,dx)}^2$
uniformly bounded, then limit points $b$ of $\{b^\varepsilon\}$ are such
that $b_x \in \mc E$.

In order to prove Theorem~\ref{t:HJ1} and Theorem~\ref{t:HJ2} we first
establish some preliminary results. Given a measurable map
$a:[0,T]\times \bb R \to [0,+\infty]$, we let $\mc L_a$ be the Hilbert
space obtained by identifying and completing the set $\big\{\varphi
\in C^\infty_{\mathrm{c}}\big((0,T)\times \bb R \big)\,:\: \langle
\langle a \varphi,\varphi\rangle \rangle<+\infty \big\}$ w.r.t.\ the
seminorm $\langle \langle a \varphi,\varphi\rangle \rangle$.

\begin{lemma2}
\label{l:Jriesz}
Let $\varepsilon>0$ and $b \in \mc B$ be such that
$J_\varepsilon(b)<+\infty$. Then there exists $E^{\varepsilon,b} \in
\mc L_{\sigma(b_x)}$ such that
\begin{equation}
\label{e:E}
b_t+f(b_x)=\frac{\varepsilon}{2} D(b_x)\,b_{x x} -\sigma(b_x)\,E^{\varepsilon,b}
\end{equation}
holds weakly on $(0,T) \times \bb R$ and $J_\varepsilon(b)=\frac 12
\|E^{\varepsilon,b}\|_{\mc L_{\sigma(b_x)}}^2$. Furthermore
$I_\varepsilon(b_x)<+\infty$ and there exists $\gamma^{\varepsilon,b} \in
\mc L_{\sigma(b_x)^{-1}}$ such that $\gamma^{\varepsilon,b}_x=0$ and
$\sigma(b_x) E^{\varepsilon,b}=
\sigma(b_x)\Psi^{\varepsilon,b_x}+\gamma^{\varepsilon,b}$, where
$\Psi^{\varepsilon,b_x}$ is defined as in Lemma~\ref{l:riesz}. 
In particular
\begin{equation}
\label{e:costdec}
J_\varepsilon(b)
  = \frac 12 \|\Psi^{\varepsilon,b_x}\|_{\mc D^1_{\sigma(b_x)}}^2
            + \frac 12 \|\gamma^{\varepsilon,b}\|_{\mc L_{\sigma(b_x)}}^2
  =I_\varepsilon(b_x)+\frac 12 \langle \langle \sigma(b_x)^{-1}
                   \gamma^{\varepsilon,b},\gamma^{\varepsilon,b} \rangle \rangle 
\end{equation}
\end{lemma2}

\proof
  The existence of $E^{\varepsilon,b}$, \eqref{e:E} and the equality
  $J_\varepsilon(b)=\frac 12 \|E^{\varepsilon,b}\|_{\mc
    L_{\sigma(b_x)}}^2$ are achieved as in Lemma~\ref{l:riesz}. We
  also have
\begin{eqnarray*}
J_\varepsilon(b) & = & 
  \sup_{\varphi \in C^\infty_{\mathrm{c}}((0,T)\times  \bb R)} 
   \Big\{
       a_\varepsilon^b(\varphi)
           -\frac 12 \langle \langle \sigma(b_x) \varphi,
                                 \varphi \rangle \rangle \Big\}
\\
&\ge &
   \sup_{\phi \in C^\infty_{\mathrm{c}}((0,T)\times \bb R)} 
      \Big\{
         a_\varepsilon^b(\phi_x) 
            -\frac 12 \langle \langle \sigma(b_x) \phi_x, 
                                 \phi_x \rangle \rangle  \Big\}
\\ 
& = &  
\sup_{\phi \in C^\infty_{\mathrm{c}}((0,T)\times \bb R)} 
        \Big\{ 
            \ell^\varepsilon_{b_x}(\phi)-\frac 12 \langle \langle
\sigma(b_x) \phi_x,\phi_x \rangle \rangle \Big\} =I_\varepsilon(b_x)
\end{eqnarray*}
By \eqref{e:3.1} and \eqref{e:E} there exists $\Psi^{\varepsilon,b_x} \in
\mc D^{1}_{\sigma(b_x)}$ such that
$\big(\sigma(b_x)E^{\varepsilon,b}\big)_x= \big(\sigma(b_x)
\Psi^{\varepsilon,b_x}_x \big)_x$, namely $\sigma(b_x)
E^{\varepsilon,b}=\sigma(b_x) \Psi^{\varepsilon,b_x}_x +
\gamma^{\varepsilon,b}(t)$ for some measurable map
$\gamma^{\varepsilon,b}:[0,T]\to \bb [-\infty,+\infty]$. It is then easy
to check \eqref{e:costdec}.
\qed

\smallskip
The following lemma is proven analogously.
\begin{lemma2}
\label{l:Jriesz2}
Let $(b,\mu) \in \mc A$ be such that $\mc J\big((b,\mu)
\big)<+\infty$. Then there exists $E^{(b,\mu)} \in \mc
L_{\mu(\sigma)}$ such that
\begin{equation*}
b_t+\mu(f)=-\mu(\sigma) E^{(b,\mu)}
\end{equation*}
and $\mc J\big((b,\mu)\big)=\frac 12 \|E^{(b,\mu)}\|_{\mc
  L_{\mu(\sigma)}}^2$. Furthermore $\mc I(\mu)<+\infty$ and there
exists $\gamma^{(b,\mu)} \in \mc L_{\mu(\sigma)^{-1}}$ such that
$\gamma^{(b,\mu)}_x=0$ and
\begin{equation}
\mc J\big((b,\mu)\big)
  = \frac 12 \|\Psi^{\mu}\|_{\mc D^1_{\mu(\sigma)}}^2
            + \frac 12 \|\gamma^{(b,\mu)}\|_{\mc L_{\mu(\sigma)}}^2
  =\mc I(\mu)+\frac 12 \langle \langle \mu(\sigma)^{-1}
                   \gamma^{(b,\mu)},\gamma^{(b,\mu)} \rangle \rangle 
\end{equation}
where $\Psi^\mu$ is defined as in Lemma~\ref{l:riesz2}.
\end{lemma2}

\begin{lemma2}
\label{l:Jequic}
The sequence of functional $\{J_\varepsilon\}$ is equicoercive on $(\mc
B,d_{\mc B})$.
\end{lemma2}

\proof 
Let $\{b^\varepsilon\} \subset \mc B$ be such that
$J_\varepsilon(b^\varepsilon) \le C_J$ for some $C_J<+\infty$. By
\eqref{e:costdec} $I_\varepsilon(b^\varepsilon_x) \le C_J$, and thus
$\{b^\varepsilon_x\}$ is precompact in $\mc U$ by Lemma~\ref{l:equic}.
We are left with the proof of the compactness of $\{b^\varepsilon\}$
w.r.t.\ the second term on the r.h.s.\ of \eqref{e:dB}. By
\eqref{e:costdec} and \eqref{e:3.6} we have that for any $N>0$,
$\varepsilon^2 \int_{[0,T]\times [-N,N]}\!dt\,dx \,(b^\varepsilon_{x
  x})^2 \le C( C_J+\varepsilon\,N+1)$ for some constant $C>0$
depending only on $f$ and $D$. It then follows by \eqref{e:E} that for
each $N>0$, $\|b^\varepsilon_t\|_{L_2([0,T]\times [-N,N]}$ is bounded
uniformly in $\varepsilon$. Since $b_x \in \mc U$ for each $b \in \mc
B$, we also have $0\le b_x^\varepsilon \le 1$. Recalling that elements
in $b$ are defined up to a constant, the conclusion follows by these
bounds on $b^\varepsilon_t$, $b^\varepsilon_x$ and compact Sobolev
embedding.
\qed

\smallskip
The following remark follows by Proposition~\ref{l:equic} and
Lemma~\ref{l:Jriesz} and the definition \eqref{e:dB} of $d_{\mc B}$.

\begin{remark2}
\label{e:Jlsc}
For each $\varepsilon>0$, $J^\varepsilon$ is lower semicontinuous on $(\mc
B,d_{\mc B})$.
\end{remark2}

\begin{lemma2} 
\label{l:bu}
For each $u \in \mc U$ such that $I_\varepsilon(u)<+\infty$ there
exists $b^{\varepsilon,u} \in \mc B$ such that $b^{\varepsilon,u}_x=u$
and $J_\varepsilon(b^{\varepsilon,u})=I_\varepsilon(u)$. Furthermore
if $b \in \mc B$ is such that $b_x=u$ and $J_\varepsilon(b)<+\infty$,
then $b_t=b^{\varepsilon,u}_t+\gamma^{\varepsilon,b}$, where
$\gamma^{\varepsilon,b} \in \mc L_{\sigma(u)^{-1}}$ is defined as in
Lemma~\ref{l:Jriesz}. Conversely, given $\gamma \in \mc
L_{\sigma(u)^{-1}}$ with $\gamma_x=0$, there exists a unique $b \in
\mc B$ such that $b_x=u$ and $b_t=b^{\varepsilon,u}_t+\gamma$.
\end{lemma2}

\proof
From the definitions \eqref{e:2.6} and \eqref{e:J}, it is not
difficult to gather
\begin{equation*}
I_\varepsilon(u)= \inf_{b \in \mc B\,:\:b_x=u} J_\varepsilon(b)
\end{equation*}
Since $J_\varepsilon$ is coercive and lower semicontinuous on $\mc B$,
and $\{b \in \mc B\,:\:b_x=u\}$ is a closed subset of $\mc B$, there
exists a $b^u$ on which the infimum is attained.

If $b$ is such that $b_x=u$ and $J_\varepsilon(b)<+\infty$, then by the
decomposition of $E^{\varepsilon,\cdot}$ in Lemma~\ref{l:Jriesz} we have
$(b-b^u)_t =\sigma(u) (E^{\varepsilon,b}-E^{\varepsilon,b^u})=
\gamma^{\varepsilon,b}$. The converse statement follows by choosing
$b(t,x)=b^u(t,x)+\int^t\!ds \gamma(s)$, which identifies a unique $b
\in \mc B$.
\qed

\smallskip
\noindent\textbf{Proof of Theorem~\ref{t:HJ1}.}
Equicoercivity follows by \eqref{e:costdec}, the equicoercivity
statement in Theorem~\ref{e:2.1} and Lemma~\ref{l:equic}.

  In order to prove the $\Gamma$-liminf inequality, let
  $\{(b^\varepsilon,\mu^\varepsilon)\} \subset \mc A$ converge to some
  $(b,\mu) \in \mc A$. It is not restrictive to assume
  $J_\varepsilon(b^\varepsilon)<+\infty$, and thus $b^\varepsilon_{x x} \in
  L_{2,\mathrm{loc}}\big([0,T]\times \bb R \big)$ and
  $\mu^\varepsilon=\delta_{b^\varepsilon_x}$. Then for each $\varphi \in
  C^{\infty}_{\mathrm{c}} \big((0,T)\times \bb R\big)$
\begin{eqnarray*}
&&\mc J_\varepsilon\big((b^\varepsilon,\mu^\varepsilon)\big)
   = J_\varepsilon(b^\varepsilon) 
\\
&&\qquad\;
\ge -\langle \langle b^\varepsilon,\varphi_t \rangle \rangle
       +\langle \langle f(b^\varepsilon_x),\varphi \rangle \rangle
       -\frac{\varepsilon}{2} 
            \langle \langle D(b^\varepsilon_x)b^\varepsilon_{x x},
                                           \varphi \rangle \rangle
       - \frac 12 \langle\langle \varphi , 
           \sigma(b^\varepsilon_x) \varphi \rangle\rangle
\end{eqnarray*} 
As in the proof of the $\Gamma$-liminf inequality in
Theorem~\ref{t:2.1}, an integration by parts shows that the third term
in the l.h.s. vanishes as $\varepsilon \to 0$. Hence
\begin{equation*}
\varliminf_{\varepsilon\to 0} 
     \mc J_\varepsilon\big((b^\varepsilon,\mu^\varepsilon)\big) \ge
            -\langle \langle b,\varphi_t \rangle \rangle
            +\langle \langle \mu(f),\varphi \rangle \rangle
            - \frac 12 \langle\langle \mu(\sigma)\varphi , 
                         \varphi \rangle\rangle
\end{equation*}
and the $\Gamma$-liminf inequality is achieved by optimizing over
$\varphi$.

Let $(b,\mu) \in \mc A$ be such that $\mc
J\big((b,\mu)\big)<+\infty$. By Lemma~\ref{l:Jriesz2} $\mc
I(\mu)<+\infty$ and by the $\Gamma$-limsup inequality in
Theorem~\ref{t:2.1} there exists a sequence $\{u^\varepsilon\} \subset
\mc U$ such that $\delta_{u^\varepsilon} \to \mu$ in $\mc M$ and
$\varlimsup_\varepsilon I_\varepsilon(u^\varepsilon)=\varlimsup_\varepsilon  \mc
I_\varepsilon(\delta_{u^\varepsilon}) \le \mc I(\mu)$. By
Corollary~\ref{l:bu} there exists $b^{\varepsilon,u^\varepsilon} \in
\mc B$ such that $b^{\varepsilon,u^\varepsilon}_x =u^\varepsilon$ and
$J_\varepsilon(b^{\varepsilon,u^\varepsilon})=I_\varepsilon(u^\varepsilon)$. Letting
$\gamma^{(b,\mu)}$ be defined as in Lemma~\ref{l:Jriesz2}, it is also
easily seen that there exists a sequence $\gamma^\varepsilon \in \mc
L_{\sigma(u^\varepsilon)^{-1}}$ such that $\gamma^\varepsilon_x=0$,
$\gamma^{\varepsilon} \to \gamma^{(b,\mu)}$ weakly in $L_{2}([0,T])$,
and $\|\gamma^\varepsilon\|_{\mc L_{\sigma(u^\varepsilon)^{-1}}} \to
\|\gamma^{(b,\mu)}\|_{\mc L_{\mu(\sigma)^{-1}}}$. Recalling
Corollary~\ref{l:bu}, we define the sequence $b^\varepsilon$ by the
requirements $b^{\varepsilon}_x=u^\varepsilon$ and $b^{\varepsilon}_t
=b^{\varepsilon,u^\varepsilon}_t +\gamma^\varepsilon$. We have
\begin{eqnarray*}
\varlimsup_{\varepsilon\to 0} 
      \mc J_\varepsilon\big(b^\varepsilon,\delta_{b^\varepsilon_x}) \big)
         &=& \varlimsup_{\varepsilon\to 0} J_{\varepsilon}(b^{\varepsilon,u^\varepsilon}) 
          + \frac 12 \langle\langle \sigma(u^\varepsilon) \gamma^\varepsilon, 
                                     \gamma^\varepsilon \rangle \rangle
\\
&\le& \mc I(\mu)  
           + \frac 12 \langle \langle \mu(\sigma) \gamma^{(b,\mu)}, 
                            \gamma^{(b,\mu)} \rangle \rangle
         = \mc J\big((b,\mu)\big)
\end{eqnarray*}
On the other hand $\delta_{b^\varepsilon_x} \to \mu$ in $\mc M$, and
it is not difficult to check $b^\varepsilon_t \to b_t$
weakly. Therefore any limit point in $\mc A$ of
$\{(b^\varepsilon,\delta_{b^\varepsilon_x})\}$ coincides with
$(b,\mu)$.
\qed

\smallskip
\noindent\textbf{Proof of Theorem~\ref{t:HJ2}.}
  If $b \in \mc Y$ is such that $(b,\delta_{b_x})$ is a measure-valued
  solution to \eqref{e:HJ}, then $b \in \mc W$. By the $\Gamma$-liminf
  inequality in Theorem~\ref{t:HJ1} we thus obtain
  $(\Gliminf_\varepsilon K_\varepsilon)(b)=+\infty$ if $b \not \in \mc
  W$. The $\Gamma$-liminf inequality on $\mc W$ follows immediately by
  (i) in Theorem~\ref{t:ecne} and \eqref{e:costdec}.

  Equicoercivity is a consequence of (ii) in Theorem~\ref{t:ecne} and
  Lemma~\ref{l:Jequic}.

  In order to prove the $\Gamma$-limsup inequality, let $b \in \mc W$
  be such that $\upbar{H}(b_x)<+\infty$. By (iii) in
  Theorem~\ref{t:ecne} there exists a sequence $\{u^\varepsilon\}
  \subset \mc X$ converging to $u:=b_x$ in $\mc X$ and such that
  $\varlimsup_\varepsilon H_\varepsilon(u^\varepsilon) \le \upbar{H}(u)$. Let
  $b^\varepsilon:=b^{\varepsilon,u^{\varepsilon}}$; by
  Corollary~\ref{l:bu} $\varlimsup_\varepsilon 
  K_\varepsilon(b^{\varepsilon,u^\varepsilon})\le K(b)$. Furthermore,
  by (i) and (ii) proven above, $\{b^\varepsilon\}$ is precompact in
  $\mc Y$ and its limit points are in $\mc W$. Let $\tilde{b} \in \mc
  W$ be a limit point of $\{b^\varepsilon\}$. Then $\tilde{b}_x=b_x$,
  since $b^\varepsilon_x =u^\varepsilon \to u=b_x$ in $\mc X$; on the
  other hand $b_t+f(b_x)=0=\tilde{b}_t+f(\tilde{b}_x)$, so that we
  also gather $b_t=\tilde{b}_t$. It follows $\tilde{b}=b$.
\qed

\textit{Acknowledgment.} 
We thank Stefano Bianchini and Corrado Mascia
for fruitful discussions. The authors gratefully acknowledge the
hospitality and the support of the Centro De Giorgi of the Scuola
Normale Superiore di Pisa, where this paper was partially written.

\end{document}